\magnification1200

\centerline{\bf Reflection trees of graphs as boundaries of Coxeter groups}
\medskip
\centerline{{\bf Jacek \'Swi\c atkowski}
\footnote{*}{The author was partially supported by  (Polish) Narodowe Centrum Nauki, grant no UMO-2017/25/B/ST1/01335.}}
\medskip
\centerline{Instytut Matematyczny, Uniwersytet Wroc\l awski}
\centerline{pl. Grunwaldzki 2/4, 50-384 Wroc\l aw, Poland}
\centerline{e-mail: {\tt swiatkow@math.uni.wroc.pl}}

\bigskip
\hfill {\it I dedicate this paper to the memory of my Parents}

\bigskip
{\bf Abstract.} 
To any finite graph $X$ (viewed as a topological space) we assosiate
some explicit compact metric space ${\cal X}^r(X)$ which we call 
{\it the reflection tree of graphs $X$}. This space is of topological 
dimension $\le1$ and its connected components are locally connected.
We show that if $X$ is appropriately triangulated (as a simplicial graph
$\Gamma$ for which $X$ is the geometric realization) then the visual boundary
$\partial_\infty(W,S)$ of the right angled Coxeter system $(W,S)$
with the nerve
isomorphic to $\Gamma$ is homeomorphic to ${\cal X}^r(X)$.
For each $X$, this yields in particular many word hyperbolic groups
with Gromov boundary homeomorphic to the space ${\cal X}^r(X)$.

\bigskip\noindent
{\bf 1. Introduction.}
\medskip

In [KK] the authors show that if the Gromov boundary $\partial G$
of a word-hyperbolic
group $G$
is connected and has no local cut point (which corresponds to the fact
that $G$ has no virtual splitting over a finite or 2-ended subgroup),
and if this boundary has topological dimension 1,
then it is homeomorphic to the Sierpi\'nski curve or to the Menger curve.  
In this paper we introduce a big family of 1-dimensional topological
spaces, each of which is either disconnected or has a local cut point, 
and each of which appears as
Gromov boundary of many right angled hyperbolic Coxeter groups.
We call these spaces {\it the reflection trees of graphs}.

The description of the above mentioned spaces is contained in Section 2,
and here we only mention that to any finite graph $X$
(viewed as a topological space and thus called a topological graph) 
there is associated one such space,
denoted ${\cal X}^r(X)$, and called {\it the reflection tree of graphs $X$}.
The space ${\cal X}^r(X)$ is connected if and only if $X$ is connected and
has no separating point, and if this is the case, the space ${\cal X}^r(X)$
is also locally connected. 
${\cal X}^r(X)$ has topological dimension 1 if and only if $X$
contains a cycle (equivalently, iff $X$ is not a tree or a disjoint union of trees).
Otherwise, except for the trivial cases when $X$ is a singleton
or doubleton, the space ${\cal X}^r(X)$ is homeomorphic to the Cantor set.
For each $X$ there are many distinct topological graphs $Y$ such that the
spaces ${\cal X}^r(X)$ and ${\cal X}^r(Y)$ are homeomorphic.
The complete topological classification of the reflection trees of graphs is not known, and in this paper we make only few initial observations
in this direction.
One such observation consists of referring to the paper [DT] by Dani
and Thomas, in which the authors study some aspects of the topology 
of boundaries of a class of hyperbolic Coxeter groups, including those
which (in view of the present paper) 
yield connected reflection trees of graphs as boundaries.
More precisely, Dani and Thomas completely determine the structure of
local cut points (encoded in an object called the Bowditch's JSJ tree)
of the studied class of boundaries.
In the light of the present paper, this gives in particular a complete 
description of the structure of local cut points of connected
reflection trees of graphs, and this allows to distinguish many of the
latter spaces up to homeomorphism.
More detailed comment concerning the relationship between
the results of [DT] and of the present paper is given in Remarks 
2.5.(4) and 2.5.(5).

\medskip
To state the main result of the paper, Theorem 1.1 below, 
we need some preparations. Recall that any right angled Coxeter system
is uniquely determined by an object called {\it the nerve}
of this system, which is an arbitrary finite flag simplicial complex. 
Moreover, to each Coxeter system $(W,S)$ there is canonically associated  
a compact metrisable space called
{\it the visual boundary} (or {\it boundary at infinity}) 
of $(W,S)$, denoted $\partial_\infty(W,S)$. 
In case when $W$ is word-hyperbolic, this boundary coincides
(up to homeomorphism) with the Gromov boundary $\partial W$ of $W$. 
We refer the reader to [D] for a detailed introduction to right angled
Coxeter systems and their boundaries.

An {\it essential vertex} of a finite topological graph $X$ is any point which
locally splits $X$ into different than 2 number of connected components.
The set of essential vertices of $X$ is obviously finite. An (open)
{\it essential edge} of $X$ is any connected component of the complement
in $X$ of the set of all essential vertices of $X$.
Note that loops and multiple edges typically occur as essential edges.
We also view as an essential edge falling into the category of loop edges 
each connected component of $X$
homeomorphic to the circle.

The main result of this paper is the following (compare Theorem 7.10 in the text,
where a slightly more general statement is given).

\medskip\noindent
{\bf 1.1 Theorem.}
{\it Let $\Gamma$ be a finite simplicial graph which is flag, and denote by $|\Gamma|$
the underlying topological graph. Suppose that each essential edge of $|\Gamma|$
which is not a loop is the union of at least 3 edges of $\Gamma$.
Let $(W_\Gamma,S_\Gamma)$ be the right angled Coxeter system with nerve $\Gamma$.
Then the visual boundary $\partial_\infty(W_\Gamma,S_\Gamma)$ is homeomorphic
to the reflection tree of graphs $|\Gamma|$, i.e. $\partial_\infty(W_\Gamma,S_\Gamma)\cong
{\cal X}^r(|\Gamma|)$.}

\medskip\noindent
{\bf Remarks.}
\item{(1)}
Note that flagness assumption for $\Gamma$ implies that any essential loop edge
of $|\Gamma|$ is the union of at least 4 edges of $\Gamma$. Note also that the
no-empty square condition, which guarantees word-hyperbolicity
of the group $W_\Gamma$ (see Theorem 12.2.1 in [D]), 
holds for $\Gamma$ satisfying the assumptions of Theorem 1.1 
exacly when each essential loop
edge of $|\Gamma|$ is the union of at least 5 edges of $\Gamma$. In particular,
any finite topological graph $X$ admits many triangulations $\Gamma$
which satisfy all the assumptions of Theorem 1.1, and for which
the group $W_\Gamma$ is word-hyperbolic. Consequently, each space
${\cal X}^r(X)$ appears as the Gromov boundary of many right angled
hyperbolic Coxeter groups. 

\item{(2)} Note that if a graph $\Gamma$ is the simplicial suspension of a finite
set consisting of at least 3 points (so that each essential edge of $|\Gamma|$
consists of two edges of $\Gamma$), then $\partial_\infty(W_\Gamma,S_\Gamma)$ 
is homeomorphic to the suspension of the Cantor set, and this space
is not homeomorphic to ${\cal X}^r(|\Gamma|)$ 
(e.g. because it is not locally connected). This example shows
that the assumption that each essential edge of $|\Gamma|$ consists of at least
3 edges of $\Gamma$ cannot be omitted (even though it certainly can be weakened).

\break

\medskip
We finish the introduction with a brief description of the organization
of the paper. 
In Section 2 we introduce {\it reflection trees of graphs},
the class of topological spaces mentioned in the title of the paper.
In the same section we derive also the most basic properties of those spaces.
In Section 3 we make a step towards the topological classification
of reflection trees of graphs, by showing that the disconnected ones
among them
are fully understood in terms of connected ones,
by means of 
the operation of the dense amalgam (as introduced in [Amalgam]).
This observation, apart from being interesting on its own, is also an ingredient in
the proof of the main result of the paper - Theorem 7.12
(which is a slightly stronger version of Theorem 1.1 stated above).
In Section 4 we observe that reflection trees of graphs allow certain 
much more flexible description, in terms of so called {\it violated reflection
inverse systems}. Section 5 deals with yet another relaxation in the description
of reflection trees of graphs. 
We describe them as limits of 
the so called {\it inverse sequences of $X$-graphs
and $X$-blow-ups}. The latter description is powerful enough to allow in Sections
7 and 8 identification of boundaries of some Coxeter groups as appropriate
reflection trees of graphs (up to homeomorphism). Before turning to
the proof of this identification, we recall in Section 6 the concept of
the Coxeter-Davis complex of a Coxeter system.
This is some CAT(0) geodesic metric space canincally associated to a Coxeter system, 
the visual boundary of which is, by definition, the boundary of the corresponding
system. In the same section, we derive
some geometric properties of such complexes, for the cases under our interest
in this paper. In Section 7, we formulate (as Theorem 7.12) the main result
of this paper in its full generality. We also state a crucial technical result,
called Approximation Lemma, which relates geometry of the
considered Coxeter-Davis
complexes with the corresponding
inverse sequences of $X$-graphs and $X$-blow-ups (introduced in Section 5).
We show in the same Section 7 how the main result follows 
from Approximation Lemma. 
Finally, in the last Section 8, we give a proof of Approximation Lemma.


\magnification1200

\bigskip
\noindent
{\bf 2. Reflection trees of graphs ${\cal X}^r( X)$}. 

\medskip
In this paper, by a {\it graph} (or {\it topological graph}) we mean
the underlying topological space $X=|\Gamma|$ of a finite simplicial graph $\Gamma$.
The {\it natural cell structure} of a graph $X$ is the coarsest cell structure
on $X$. The vertices of this structure (also called the {\it essential vertices}
of $X$) are these points $x\in X$ which locally split $X$ into different
than 2 number of connected components. The open edges of this structure
(called {\it essential open edges}) are the connected components of the
complement in $X$ of the set of all essential vertices. Note that loops
(loop edges) and multiple edges typically occur in the natural cell structure
of a graph $X$. We also view as an essential loop edge each connected
component of $X$ homeomorphic to the circle.
Obviously, the natural cell structure of each graph $X$ consists of
finitely many vertices and edges.

\medskip
In this section,
for any graph $X$, we describe a compact metric space
${\cal X}^r(X)$, of topological dimension $\le1$, called 
{\it the reflection tree of graphs} $X$. 
To start describing this space, 
we need an auxiliary object,
namely a countable dense subset $D\subset X$ containing all essential
vertices. Such a subset is topologically unique, in the precise sense provided by
the following easy observation.

\medskip\noindent
{\bf 2.1 Lemma.}
{\it Let $D_1,D_2$ be two countable dense subsets of $X$, both
containing all essential vertices of $X$. Then there is a homeomorphism
$h:X\to X$ such that $h(D_1)=D_2$. Moreover, $h$ can be chosen
to preserve all essential vertices and edges of $X$.}

\medskip
Given $X$ and $D\subset X$ as above, we describe now
an inverse system ${\cal S}(X,D)$ whose limit is, by definition,
the reflection tree of graphs $X$.

\medskip\noindent
{\it The underlying poset of ${\cal S}( X,D)$.}

Let $T$ be a tree each vertex of which has valence equal to the cardinality
of the subset $D$.
Denote by $V_T,E_T$ the sets of vertices and (unoriented) edges of $T$,
respectively. Let $\cal F$ be the poset of all finite subtrees of $T$, ordered
by inclusion. The inverse system $S( X,D)$ will have a form
$$
{\cal S}( X,D)=\big( \{ X_F:F\in{\cal F} \} , \{  \pi_{F',F}:F\subset F' \}
 \big),
$$
where each $X_F$ is a compact metric space (some graph formed of
few copies of $ X$), and each $\pi_{F',F}:X_{F'}\to X_F$
is a continuous map, as described below. 

\medskip\noindent
{\it Blow-up of $\, X$ at a point or at a finite set.}

For any point $x\in X$, denote by $ X^\#(x)$ and call
{\it the blow-up of $ X$ at $x$} the graph obtained by attaching
to $ X\setminus\{x\}$ as many points (which all become vertices 
of valence 1 in
$ X^\#(x)$) as the number of components into which $x$ splits
its any normal neighbourhood in $ X$. Denote the set of all 
those attached vertices by $P_x$. 
Denote also by 
$\rho_x:  X^\#(x)\to  X$ the {\it blow down map}
which projects $P_x$ to $x$ and which is identical on the remaining part
$ X\setminus\{x\}$.

For any finite subset $J\subset X$, denote by $ X^\#(J)$
the graph obtained from $ X$ by performing blow-ups at all points
$x\in J$ (the result does not depend on the order). 
Given any finite subsets $J,K$ of $ X$ such that $J\subset K$, denote by 
$\rho_{K,J}: X^\#(K)\to X^\#(J)$ the blow down map
which shrinks each of the subsets $P_x:x\in K\setminus J$ to the corresponding
point $x$, and which is identical on the remaining parts of the graphs.

\medskip\noindent
{\it The spaces $X_F$.}

To each vertex $t\in V_T$ associate a copy of the graph $ X$,
and denote it $X_t$.
Equip also $T$ with a labelling $\lambda:E_T\to D$ such that for any
$t\in V_T$, denoting by $A_tT$ the set of edges in $T$ adjacent to $t$,
the restriction of $\lambda$ to $A_tT$ is a bijection on $D$.
Such a labelling clearly exists, and is unique up to an automorphism of $T$.

Intuitively, for each edge $e=[t_1,t_2]\in E_T$ the label $\lambda(e)$
will represent (in the construction of the spaces $X_F$ given formally in the next
paragraph) an operation of ``doubling'' $ X$ at the point $\lambda(e)$, 
consisting of blowing up the 
copies $X_{t_1},X_{t_2}$ at the point $\lambda(e)$ in each copy,
and of gluing to each other the sets $P_{\lambda(e)}$ in the blown up copies
through the identity map.

Fix a subtree $F\in{\cal F}$, and denote by $V_F, E_F$ its vertex and edge set,
respectively. For any $t\in V_F$, denote by $A_tF$ the set of edges in $F$
adjacent to $t$. Put
$$
X_F=\bigsqcup_{t\in V_F}X_t^\#(\lambda(A_tF))\slash\sim,
$$
where $\bigsqcup$ denotes the disjoint union, and where 
the equivalence relation $\sim$
is induced by the following equivalences: for each $e=[t_1,t_2]\in E_F$
and each $p\in P_{\lambda(e)}$, identify 
$p\in P_{\lambda(e)}\subset X_{t_1}^\#(\lambda(A_{t_1}F))$ with
$p\in P_{\lambda(e)}\subset X_{t_2}^\#(\lambda(A_{t_2}F))$.

Note that for each $t\in V_F$ 
the set $X_t^\#(\lambda(A_tF))$ 
is naturally a subspace in $X_F$.

\medskip\noindent
{\it The maps $\pi_{F',F}$.}

Given any finite subtrees $F, F'$ of $T$ such that $F\subset F'$,
define $\pi_{F',F}:X_{F'}\to X_F$ as follows.
For $t\in V_F$, the restriction of $\pi_{F',F}$ to
$X_t^\#(\lambda(A_tF'))\subset X_{F'}$ coincides with the blow down map
$\rho_{\lambda(A_tF'),\lambda(A_tF)}$
(whose image $X_t^\#(\lambda(A_tF))$ is viewed as a subset in $X_F$).
For $t\in V_{F'}\setminus V_F$, the subset $X_t^\#(\lambda(A_tF'))$
is mapped by $\pi_{F',F}$ to the point $\lambda(e)\in X_s^\#(\lambda(A_sF))$,
where $s\in V_F$ is the closest to $t$ (in the polygonal metric in $T$)
vertex of the subtree $F$, and where $e$ is the first edge on the path from $s$
to $t$ in $T$. 
Obviously, $X_s^\#(\lambda(A_sF))$ is viewed in the previous sentence
as a subset in $X_F$.

\medskip
The above description of the inverse system ${\cal S}( X,D)$,
together with Lemma 2.1, immediately yield the following.

\medskip\noindent
{\bf 2.2 Lemma.}
{\it  Let $D_1,D_2$ be two countable dense subsets of $ X$, both
containing all essential vertices of $ X$. Then the inverse systems
${\cal S}( X,D_1)$ and ${\cal S}( X,D_2)$ are isomorphic.
Consequently, the inverse limits $\lim_\leftarrow {\cal S}( X,D_1)$
and  $\lim_\leftarrow {\cal S}( X,D_2)$ are homeomorphic.}

\medskip
In view of Lemma 2.2, we will denote any inverse system
${\cal S}( X,D)$ simply by ${\cal S}_ X$, and call it
{\it the standard reflection inverse system for $ X$}.

\medskip\noindent
{\bf 2.3 Definition.}
The {\it reflection tree of graphs $ X$}, denoted ${\cal X}^r( X)$,
is the space obtained as inverse limit $\lim_\leftarrow{\cal S}_ X$ 
of the standard reflection inverse system for $ X$.

\medskip
At the end of this section, in Remark 2.6, we indicate an alternative 
way of describing reflection trees of graphs. Although we do not use this
alternative description in the present paper, it exhibits connections
of reflection trees of graphs with some other classes of spaces
studied in the context of boundaries of groups, notably with 
{\it trees of manifolds}.

Next result collects some basic properties of the reflection trees of graphs.

\medskip\noindent
{\bf 2.4 Lemma.}
\item{(1)} {\it Reflection tree of a singleton graph is the empty space,
and reflection tree of a doubleton is a doubleton.}
\item{(2)} {\it If $ X$ is a tree or a disjoint union of trees, and if 
it is not a singleton and not a doubleton, 
then the space ${\cal X}^r( X)$
is homeomorphic to the Cantor set.}
\item{(3)} {\it If $ X$ is homeomorphic to the circle $S^1$ then 
${\cal X}^r( X)\cong S^1$.}

\item{(4)} {\it For any graph $ X$ the topological dimension
of the space ${\cal X}^r( X)$ is $\le1$.}
\item{(5)} {\it If $ Y$ is a subgraph of a graph $ X$ 
(with respect to the natural cell structure in $X$),
then the space ${\cal X}^r( X)$ contains an
embedded copy of ${\cal X}^r(Y)$. Moreover, if $Y$ is a proper subgraph
not reduced to a single vertex then there are countably infinitely many
natural embeddings of ${\cal X}^r(Y)$ in ${\cal X}^r(X)$,
with pairwise disjoint images.}

\item{(6)} {\it We have $\dim {\cal X}^r( X)=1$ iff $ X$ is not a tree
or a disjoint union of trees (i.e. iff $ X$ contains a cycle).}

\item{(7)} {\it The space ${\cal X}^r( X)$ is connected iff $\, X$
is connected and has no separating vertex.}

\item{(8)} {\it If ${\cal X}^r( X)$ is connected, then it is also locally
connected and contains local cut points.}

\medskip
Before giving the proof of the above lemma, we make few comments.

\break

\medskip\noindent
{\bf 2.5 Remarks.} 
\item{(1)} Note that, due to part (8) of Lemma 2.4, 
the following 1-dimensional compacta  
do not appear among reflection trees of graphs:
\itemitem{$\bullet$} suspension of the Cantor set and the join of two
Cantor sets (these spaces are not locally connected);
\itemitem{$\bullet$} Sierpi\'nski curve and Menger curve
(these spaces have no local cut points).

\item{(2)}
Further topological properties of the spaces ${\cal X}^r( X)$ wait to be studied. In particular, it is an open problem to classify the spaces ${\cal X}^r( X)$
up to homeomorphism.

\item{(3)} In Section 3 we present a result (Proposition 3.8)
which gives a full description of disconnected reflection trees of graphs,
in terms of connected ones. Thus, the topological classification
of the reflection trees of graphs reduces to the classification
of the connected ones (which is anyway the most interesting case).

\item{(4)} Some further topological properties of reflection trees
of graphs can be deduced from the results of Dani and Thomas in [DT]. 
In that paper the authors give, among others, a complete
description of the structure of local cut points (the so called Bowditch's
JSJ tree) for Gromov boundaries of a class of right angled hyperbolic
Coxeter groups, including all hyperbolic Coxeter groups appearing
in Theorem 1.1 of the present paper. This gives the description
of the structure of local cut points for all connected reflection trees of graphs.
This structure (the Bowditch's JSJ tree) is a powerfull topological invariant,
containing full information about degrees of local cut points, as well as of their
mutual position inside the space.
In our context, this invariant allows to distinguish many 
of the connected reflection
trees of graphs, up to homeomorphism. 

\item{(5)}
As a complement to the previous remark (4), the reader should keep in mind
that continua with the same Bowditch's JSJ trees are not necessarily homeomorphic.
Hence, Bowditch's JSJ trees cannot be used to deduce that the boundaries
$\partial_\infty(W_\Gamma,S_\Gamma)$ as in Theorem 1.1, for distinct $\Gamma$ having the
same underlying topological graph $|\Gamma|$, are homeomorphic.

\medskip\noindent
{\bf Proof of Lemma 2.4:}

The proof of part (1) is straightforward, and the proof of
part (2)
is an exercise on checking the conditions from the
following well known
characterization of the Cantor set:  it is the unique compact
metric space which is totally disconnected and has no isolated points.
We skip the details.

To see part (3), note first that if $X$ is a circle then all of the spaces
$X_F$ in the system ${\cal S}_X$ are also homeomorphic to a circle.
Moreover, the maps $\rho_{F',F}$ all have the following form: a circle
is mapped to a circle by means of shrinking few pairwise disjoint
subsegments in the source circle to points. Obviously, each such map
$\rho_{F',F}$ can be approximated  by a homeomorphism of the
involved circles, i.e. it is a near-homeomorphism. By a result of Brown [Br],
the inverse limit of a system with bonding maps which are near-homeomorphisms coincides with the spaces in the system,
so in our case it is a circle, as required.

To get part (4), note that  since each space $X_F$ is a graph, 
the well known estimate for the dimension of an inverse limit
yields in our case the required inequality, as follows:
$$
\dim{\cal X}^r( X)=
\dim\lim_{\longleftarrow}(\{ X_F \},\{ \pi_{F',F} \})\le\sup\{ \dim X_F:F\in{\cal F} \}
\le1.
$$

The proof of part (5) requires some preparations concerning morphisms
of inverse systems. As a reference for this subject  we use
Section 2.5 in the book [En]. 
Recall that given inverse systems
$$
{\cal S}_1=(\{ Y_\lambda:\lambda\in\Lambda \}, \{ \pi^1_{\lambda',\lambda} \}) \hbox{\quad and \quad}  
{\cal S}_2=(\{ X_\mu:\mu\in M \}, \{ \pi^2_{\mu',\mu} \}),
$$
with bonding maps $\pi^1_{\lambda',\lambda}:X_{\lambda'}\to X_\lambda$
and  $\pi^2_{\mu',\mu}:X_{\mu'}\to X_\mu$,
a {\it morphism} ${\bf f}:{\cal S}_1\to{\cal S}_2$ is a pair 
${\bf f}=(\phi,\{ f_\mu:\mu\in M \})$ such that:
\item{$\bullet$} $\phi:M\to\Lambda$ is an order preserving map of
underlying posets such that the image $\phi(M)$ is a cofinal subset
of $\Lambda$;
\item{$\bullet$} $f_\mu:Y_{\phi(\mu)}\to X_\mu$ are continuous maps,
one for each $\mu\in M$,
such that whenever $\mu\le\mu'$ then
$$
f_\mu\pi^1_{\phi(\mu'),\phi(\mu)}=\pi^2_{\mu',\mu}f_{\mu'}.
$$

\noindent
A morphism $\bf f$ as above induces a continuous map 
$f:\lim_\leftarrow{\cal S}_1\to\lim_\leftarrow{\cal S}_2$ such that if
$p\in\lim_\leftarrow{\cal S}_1$ 
is represented by a thread $\{y_\lambda\}_{\lambda\in\Lambda}$  then
its image $f(p)$ is represented by 
$$
\{ f_\mu(y_{\phi(\mu)}) \}_{\mu\in M}.
$$
Moreover, we have the following result (see Lemma 2.5.9 in [En]).

\medskip\noindent
{\bf Fact A.}
{\it If all the maps $f_\mu$ of a morphism ${\bf f}=(\phi,\{ f_\mu \})$
are injective, then the induced map $f$ is also injective.}

\medskip
To proceed, under assumptions on graphs $X$ and $Y$ as in the statement
of part (5), we introduce a uniform notation that indicates dependence
of the involved objects on a graph $Z$ from
the set $\{ X,Y \}$. $D_Z$ denotes a countable dense subset of $Z$
containing all essential vertices of $Z$. $T_Z$ is a tree whose valence at
every vertex is equal to the cardinality of the set $D_Z$, and whose vertex
and edge sets are denoted $V_{T_Z},E_{T_Z}$, respectively.
$T_Z$ is equipped with a labelling $\lambda_Z:E_{T_Z}\to D_Z$
such that its restriction to any set of edges adjacent to a fixed vertex
is a bijection on $D_Z$. We denote by ${\cal F}_Z$ the poset of all finite
subtrees of $T_Z$, and by ${\cal S}_Z={\cal S}(Z, D_Z)$ the corresponding
inverse system as described above in this section, 
consisting of the spaces 
$Z_F:F\in{\cal F}_Z$ and the maps $\pi^Z_{F',F}:Z_{F'}\to Z_F$,
and such that $\lim_\leftarrow{\cal S}_Z={\cal X}^r(Z)$.
Furthermore, we assume (without loss of generality) that $D_Y=Y\cap D_X$.

Let $S$ be any maximal subtree of $T_X$ with the property that 
for each edge $e$ of $S$ we have $\lambda_X(e)\in D_Y$.
Note that if $Y$ is a proper subgraph of $X$ then there are countably 
infinitely many subtrees $S$ as above, and the vertex sets of these subtrees
form a partition of the vertex set $V_{T_X}$.
For each such $S$ there is a label preserving isomorphism 
$\psi:T_Y\to S$.

After fixing any $S$ and $\psi$ as above, 
denote by ${\cal F}_X^S$ the subposet
of ${\cal F}_X$ consisting of those finite subtrees $F$ of $T_X$
which have nonempty intersection with $S$.
Note that ${\cal F}_X^S$ is cofinal in ${\cal F}_X$, and thus the
restricted inverse system ${\cal S}_X^S:={\cal S}_X|_{{\cal F}_X^S}$
has the same limit as ${\cal S}_X$.
We describe a morphism ${\bf f}^\psi=(\phi,\{ f_F^\psi \}):{\cal S}_Y
\to{\cal S_X^S}$ as follows.
Given $F\in{\cal F}_X^S$, we put $\phi(F):=\psi^{-1}(F\cap S)$
and we note that this is a nonempty subtree of $T_Y$.
Obviously, $\phi$ is order preserving and surjective, so its image
is cofinal in ${\cal F}_Y$.
For any $F\in{\cal F}_X^S$ define 
$f_F^\psi:Y_{\psi^{-1}(F\cap S)}\to X_F$
as the injective map induced by the natural inclusions 
$Y_s\to X_{\psi(s)}$, for all vertices $s$ of the subtree $\psi^{-1}(F\cap S)$
We skip the straightforward details concerning the description of these maps,
the verification that they are injective, and that they satisfy the appropriate
commutativity equations, thus forming a morphism.

Let $f^\psi:\lim_\leftarrow{\cal S}_Y\to\lim_\leftarrow{\cal S}_X^S$
be the map induced by the morphism ${\bf f}^\psi$. Due to the above recalled
Fact A, this map is injective. Since 
$\lim_\leftarrow{\cal S}_Y\cong{\cal X}^r(Y)$ and 
$\lim_\leftarrow{\cal S}_X^S\cong{\cal X}^r(X)$, the map $f^\psi$
embeds a copy of ${\cal X}^r(Y)$ in ${\cal X}^r(X)$,
as required in the first part of the assertion.
To prove the second assertion of part (5), it is sufficient to show that
the images of the embeddings $f^\psi$ corresponding
to distinct subtrees $S$ are pairwise disjoint.

Fix two distinct subtrees $S,S'$ of $T_X$ equipped with label preserving
isomorphisms $\psi:T_Y\to S$ and $\psi':T_Y\to S'$. Consider any subtree
$F\in{\cal F}_X$ which has nonempty intersection with both $S$ and $S'$.
Note that, since $S\cap S'=\emptyset$, we have
$$
f_F^\psi(Y_{\psi^{-1}(F\cap S)})\cap f_F^{\psi'}(Y_{\psi^{-1}(F\cap S')})=
\emptyset,
$$
and this obviously implies that the images of $f^\psi$ and $f^{\psi'}$
are disjoint, which completes the proof of part (5).

\medskip

One implication in part (6) follows by observing that if $ X$ contains
a cycle then, due to parts (3) and (5), ${\cal X}^r( X)$ contains
an embedded copy of $S^1$, and due to part (4) this means that
$\dim{\cal X}^r( X)=1$.  The remaining implication follows
from part (1) and (2), in view of the well known fact that
the dimension of the Cantor set is 0.

To prove (7), suppose first that $ X$ either is not connected, or has
a separating vertex. It is not hard to see that then some space $X_F$ 
from the inverse system ${\cal S}_ X$ is not connected.
Since the bonding maps are all surjections (which is generally true in the 
inverse systems ${\cal S}_ X$, unless $ X$ is a singleton), 
we get disconnectedness
of the inverse limit (i.e. disconnectedness of the space ${\cal X}^r( X)$).
To get the converse, note that if $ X$ is connected and has no 
separating vertex then all the
spaces $X_F$ are connected. Since in such a case the inverse limit is connected too,
this completes the proof of part (7).

To get the local connectedness assertion in part (8), note that if ${\cal X}^r( X)$
is connected then (due to part (7)) $ X$ is connected and has no
separating vertex. As it was already observed,  each space $X_F$ is then connected.
For a similar reason, point preimages for all bonding maps
of the system ${\cal S}_ X$ are also connected. 
Consequently, the assertion follows by the general fact that  limit of a monotone inverse
system of locally connected continua is locally connected,
see [Capel].

To get existence of local cut points (i.e. the second assertion of part (8)),
note that under the assumption of this part the graph $X$ has at least one edge,
and hence $D$ contains points which are distinct from essential vertices.
Consider a point $x\in D$ not being an essential vertex of $ X$,
and choose an edge $e=[t_1,t_2]$ in $T$ labelled with $x$ (i.e. such that
$\lambda(e)=x$). Let $F$ be the subtree of $T$ spanned on the vertices
$t_1,t_2$ (coinciding with the edge $e$). Then any point 
$p\in P_x$ (viewed as a point of $X_F$ under the canonical embedding
of $X^\#_{t_1}(x)$ or $X^\#_{t_2}(x)$ in $X_F$) is a local cut point of $X_F$.
Moreover, it is not hard to realize that the preimage of $p$
in the inverse limit $\lim_\leftarrow{\cal S}( X,D)={\cal X}^r( X)$
is a singleton, and we denote the unique point in this preimage by $\tilde p$.
Now, let $U$ be any connected neighbourhood in $X_F$ which is disconnected
by $p$, and let $A,B$ be some nonempty open and closed subsets 
of $U\setminus\{p\}$  forming
a partition of this subset. Denoting by 
$\pi_F:\lim_\leftarrow{\cal S}( X,D)\to X_F$ the canonical projection
associated to ${\cal S}(X,D)$, we claim that
\item{(1)} $\pi_F^{-1}(U)$ is a connected neighbourhood of $\tilde p$, and
\item{(2)} the sets $\pi_F^{-1}(A), \pi_F^{-1}(B)$ form a nontrivial
open and closed partition of $\pi_F^{-1}(U)\setminus\{ \tilde p \}$.

\noindent
Property (1) above follows by observing (similarly as in the proof of part (7))
that for each $F'\supset F$ the preimage $\pi_{F',F}^{-1}(U)$ is connected.
Property (2) follows from the fact that $\pi_F^{-1}(p)=\{ \tilde p \}$,
and from surjectivity of $\pi_F$ (which in turn follows from surjectivity
of the bonding maps in the system ${\cal S}(X,D)$).
Properties (1) and (2) mean that 
$\tilde p$ is a local cut point in ${\cal X}^r( X)$,
hence the assertion.

\medskip\noindent
{\bf 2.6 Remark.}
Each of the spaces ${\cal X}^r( X)$ can be alternatively described using
the setting of the paper [Tr-metr], as limit of some
tree system $\Theta^r( X)$ associated to a graph $ X$ in a natural way.
We explain the details below, 
referring the reader to [Tr-metr] for the general introduction
concerning tree systems, as well as details of terminology and notation.
However, it should be made clear that in the present paper we make
no essential use of this alternative description of reflection trees of graphs.

Given any $ X$,
consider a family $\cal D$ of (closed)
normal neighbourhoods $\Delta(x)$ of points $x$ in $ X$  
such that:

\itemitem{(b1)} the sets from $\cal D$ are pairwise disjoint, and their
union is dense in $ X$;

\itemitem{(b2)} for each essential vertex $v$ in $ X$ there is
a set of form $\Delta(v)$ in $\cal D$.

\noindent

\noindent
It is not hard to observe that
a family $\cal D$ as above is unique
up to a homeomorphism of $ X$ preserving
all essential vertices and edges.

The {\it reflection tree system of graphs $ X$} is the
tree system $\Theta^r( X)$ described uniquely up
to an isomorphism of tree systems by the following conditions:
\itemitem{(r1)} each vertex space is homeomorphic to the space
$ X\setminus\cup_{\Delta\in{\cal D}}\,\hbox{int}(\Delta)$,
where ``int'' denotes here the topological interior in $ X$;
\itemitem{(r2)} for each vertex space 
$ X\setminus\cup_{\,\Delta\in{\cal D}}\,\hbox{int}(\Delta)$
the family of peripheral subsets in this space coincides
with the family $\{ \hbox{bd}(\Delta):\Delta\in{\cal D} \}$,
where ``bd'' denotes here the topological boundary in $ X$;
\itemitem{(r3)} the glueing maps in the system are the identity maps
$\hbox{bd}(\Delta)\to\hbox{bd}(\Delta)$ between the copies of the same
peripheral set $\hbox{bd}(\Delta)$ in the two adjacent vertex spaces.

Note that, after shrinking to points the peripheral subsets 
$\hbox{bd}(\Delta):\Delta\in{\cal D}$ of the space 
$ X\setminus\cup_{\,\Delta\in{\cal D}}\,\hbox{int}(\Delta)$,
we get a space homeomorphic to $ X$, and the set of points
obtained by shrinking the subsets $\hbox{bd}(\Delta)$ is a countable dense
subset in the quotient, containing its all essential vertices.
Using this fact, it is straightforward to verify that 
the inverse system
${\cal S}_ X$ is isomorphic to the system ${\cal S}_{\Theta^r( X)}$
associated to the tree system $\Theta^r( X)$ as in Section 1.D
of [Tr-metr]. (The latter system is called in [Tr-metr] the standard inverse
system associated to the tree system $\Theta^r( X)$.)
It follows then from Proposition 1.D.1 in [Tr-metr] that the reflection tree
of graphs ${\cal X}^r( X)$ is homeomorphic to the limit
$\lim\Theta^r( X)$ of the tree system $\Theta^r( X)$.

The above observation allows to make
the following comment, which puts the spaces ${\cal X}^r( X)$
in a wider perspective. The pair of data
$$
( X\setminus\cup_{\,\Delta\in{\cal D}}\,\hbox{int}(\Delta), 
\{ \hbox{bd}(\Delta):\Delta\in{\cal D} \} ),
$$
which could be called
{\it the densely punctured graph $ X$}, is fairly analogous to
an object $M^\circ$ described in section 1.E.2 of [Tr-metr], called there
the densely punctured manifold $M$. The latter object leads (in a way similar
as in the above description of the tree system $\Theta^r(X)$)
to a tree system ${\cal M}(M)$ called the dense tree system of manifolds $M$.
(In fact, if $ X=S^1$ then the tree systems $\Theta^r( X)$
and ${\cal M}(S^1)$ are just isomorphic.)
Consequently, reflection trees of graphs ${\cal X}^r( X)$
are analogous to the spaces ${\cal X}(M)$,
called {\it trees of manifolds}, obtained as limits of the systems ${\cal M}(M)$
(or alternatively as inverse limits of the standard inverse systems
$\Theta_{{\cal M}(M)}$ associated to ${\cal M}(M)$). 
The latter spaces are known in topology for
many years (see e.g. [AS], [J], [St]), and more recently they have been studied in the context of boundaries of groups (see [Fi], [PS] and [Tr-mfld]).
In this context, our main result (Theorem 1.1) may be viewed as analogon
of Fisher's result [Fi] saying that Coxeter groups with PL manifold nerves
have boundaries homeomorphic to the appropriate trees of manifolds
(see Theorem 3.A.3 in [Tr-metr] for correct statement of this result).


\magnification1200

\bigskip\noindent
{\bf 3. Disconnected spaces ${\cal X}^r(X)$ as dense amalgams.}

\medskip
Recall that the space ${\cal X}^r(X)$ is connected if and only if 
$X$ is connected and 
has no separating essential vertex. If this is not the case, $X$ decomposes naturally
into a family of subgraphs, called {\it blocks} of $X$, which are connected and
which have no separating essential vertices 
(we describe this decomposition in detail below). 
A block is {\it nontrivial} if it is not a single (isolated) vertex of $X$
and not a single non-loop edge of $X$.

The aim of this section is to show that 
if the space ${\cal X}^r(X)$ is not connected,
then it is homeomorphic to 
some uniquely determined combination of the connected spaces
${\cal X}^r(Y)$ (whose copies appear as the nontrivial connected components
in ${\cal X}^r(X)$), where $Y$ runs through the family of nontrivial blocks
of $X$. The above mentioned ``combination of spaces'' corresponds to the operation
of the dense amalgam, introduced in [Amalgam] and recalled below in
Definition 3.6. See Proposition 3.8 below
for a precise statement of the above announced main result of the section.

We start with recalling the concepts of block decomposition of a graph
and dense amalgam of a family of compact metric spaces.

\medskip\noindent
{\it Block decomposition of a graph $X$.}

\medskip\noindent
{\bf 3.1 Definition.}
We say that a topological graph $X$ is {\it non-separable} if it is connected
and has no
separating essential vertex. (In particular, the circle is non-separable.)
A {\it block} of $X$
is its any subgraph (for the natural cell structure consisting
of essential vertices and edges)
which is maximal for the inclusion in the family of all non-separable
subgraphs of $X$, and which is not homeomorphic to a segment, 
and not a single isolated vertex of $X$. 



\medskip
We now describe the set of blocks of a graph $X$ in terms
of a sequence of operations called splittings.
We will need this description
in the arguments both in this section, and in Section 7 (in the proof
of Theorem 7.12).

\medskip\noindent
{\bf 3.2 Definition.} 
A {\it natural subgraph} of a graph $X$ is its any subgraph
for the natural cell structure (consisting of  essential edges and essential
vertices).
A {\it splitting} of a connected graph $X$ is a pair of connected
natural subgraphs $X_1,X_2$, both distinct from a singleton, 
such that $X_1\cup X_2=X$ and $X_1\cap X_2$
is a single essential vertex of $X$. The subgraphs $X_1$ and $X_2$, viewed as topological
graphs (equipped with their new natural cell structures), are then called
{\it the parts} of this splitting.

\medskip
Observe that a connected graph has no splitting if and only if it is non-separable.

\medskip\noindent
{\bf 3.3 Definition.}
A {\it split decomposition} of a connected graph $X$ is any
sequence of splittings having the following recursive description:
\item{$\bullet$} the empty sequence of splittings forms the trivial 
split decomposition of $X$, and the set of factors of this decomposition
is $\{ X \}$;
\item{$\bullet$} a single splitting of $X$ is a split decomposition, and its set
of factors is the set consisting of 
the
two parts of the splitting; 
\item{$\bullet$} if some sequence of splittings is a split decomposition of $X$,
and if a collection $\{ X_1,\dots,X_m \}$ of natural subgraphs of $X$ 
is the set of factors of this split decomposition, then adding to this 
sequence a splitting of one of those factors, say $X_m$, we also get a split
decomposition of $X$; moreover, if $X'_m,X''_m$ are the parts of the
above splitting of $X_m$, the set of factors of the new split decomposition
is $\{ X_1,\dots,X_{m-1},X'_m,X''_m \}$. 

\noindent
A split decomposition of $X$ is {\it terminal} if its every factor is non-separable.
A factor of a split decomposition is {\it trivial} if it is a singleton or 
if it is homeomorphic to a segment. Otherwise a factor is {\it nontrivial}.

\medskip
Every finite connected topological graph admits a terminal split decomposition.
This follows by observing that, after adding a new splitting to a split
decomposition, the total number of essential edges in the factors
does not increase (actually, it may happen that this number strictly decreases,
since the natural cell structure of a subgraph might be coarser than the
structure induced from the initial graph).
Moreover, each factor of any split decomposition of $X$ is a natural
subgraph of $X$. 

It is not hard to observe that any block of a connected graph $X$ is contained
in precisely one factor of any split decomposition of $X$.
As an easy consequence, the following result holds.

\medskip\noindent
{\bf 3.4 Lemma.}
{\it Let $X$ be a finite connected topological graph.
The set of nontrivial factors of any terminal split decomposition of $X$ coincides with the set of blocks of $X$.}

\medskip
Lemma 3.4 obviously implies the following.

\medskip\noindent
{\bf 3.5 Corollary.}
{\it Let $X$ be any finite topological graph (not necessarily connected).
the set of blocks of $X$ coincides with the union of the sets of nontrivial factors
of any terminal split decompositions of all connected components of $X$.}

\bigskip\noindent
{\it The dense amalgam of metric compacta.}

\medskip
We recall the concept of the dense amalgam of metric compacta,
as introduced in [Amalgam].

\medskip\noindent
{\bf 3.6 Definition.} The {\it dense amalgam} is
an operation $\widetilde\sqcup$ which to any finite tuple $B_1,\dots,B_k$ of nonempty metric compacta associates a  
metric compactum $$A=\widetilde\sqcup(B_1,\dots,B_k)$$ 
determined uniquely up to homeomorphism by the following.
The space
$A$ can be equipped with a countable infinite family $\cal A$ of subsets,
partitioned as ${\cal A}={\cal A}_1\sqcup\dots\sqcup{\cal A}_k$,
such that:

\itemitem{
(a1)} the subsets in $\cal A$ are pairwise disjoint and for each
$i\in\{ 1,\dots,k \}$ the family ${\cal A}_i$ consists of embedded copies of the space $B_i$;

\itemitem{
(a2)} the family $\cal A$  
is {\it null}, i.e. for any metric on $A$ compatible with the topology the diameters of sets in $\cal A$ converge to 0;

\itemitem{
(a3)} each subset from the family ${\cal A}$ is a boundary subset of $A$ (i.e. its complement is dense in $A$);

\itemitem{
(a4)} for each $i$, the union of the family ${\cal A}_i$ is dense in $A$;

\itemitem{
(a5)} any two points of $A$ which do not belong to the same subset
of $\cal A$ can be separated from each other by an open and closed subset $H\i A$ which is $\cal A$-saturated (i.e. such that any element of $\cal A$ is either contained in or disjoint with $H$).

\medskip
The uniqueness claim appearing in the above definition is proved in
[Amalgam], as Theorem 0.2. Section 1 of the same paper contains an effective
construction of the space $\widetilde\sqcup(B_1,\dots,B_k)$, for any
given $B_1,\dots,B_k$.

We call the result $\widetilde\sqcup(B_1,\dots,B_k)$ of the above described operation 
{\it the dense amalgam of the spaces $B_1,\dots,B_k$}.  
Obviously, the dense amalgam of any family $B_1,\dots,B_k$
of spaces is a disconnected perfect compact metric space.
Moreover, it follows easily from the condition (a5) above that 
if the spaces $B_1,\dots,B_k$ are connected, then the connected
components of their dense amalgam are precisely the subsets from the family
$\cal A$ and the singletons from the complement of the union $\cup{\cal A}$.

The following result concerning ``algebraic'' properties of the operation
of dense amalgam 
is proved as Proposition 0.1 in [Amalgam].

\medskip\noindent
{\bf 3.7 Proposition.} 
{\it For any nonempty metric compacta $B,B_1,\dots,B_k$ the following
equalities hold, up to homeomorphism:}

\item{(1)} {\it $\widetilde\sqcup(B_1,\dots,B_k)=
\widetilde\sqcup(B_1\sqcup\dots\sqcup B_k)$;} 

\item{(2)} {\it $\widetilde\sqcup(B_1,\dots,B_k)=\widetilde\sqcup(B_1,\dots,B_{i-1},\widetilde\sqcup(B_i,\dots, B_k))$ for any $k\ge1$ and any $1\le i\le k$;}

\item{(3)} {\it $\widetilde\sqcup(B,B_1,\dots,B_k)=\widetilde\sqcup(B,B,B_1,\dots,B_k)$ for any $k\ge0$;}

\item{(4)} {\it for any totally disconnected nonempty compact metric space $Q$, and any $k\ge1$, we have}
$$
\widetilde\sqcup(B_1,\dots,B_k,Q)=\widetilde\sqcup(B_1,\dots,B_k);
$$

\item{(5)} {\it for any totally disconnected space $Q$ (including the case when
$Q$ is a singleton) we have $\widetilde\sqcup(Q)=C$, where $C$ is the Cantor set.}

\medskip
In consistency with the properties from Proposition 3.7, we use the convention that $\widetilde\sqcup(\emptyset):=C$ and $\widetilde\sqcup(\emptyset,X_1,\dots,X_k):=\widetilde\sqcup(X_1,\dots,X_k)$.
Moreover, if $[\emptyset]$ denotes the empty family of metric compacta,
we set $\widetilde\sqcup([\emptyset]):=C$.

\bigskip\noindent
{\it Main result of the section, and its proof.}

\medskip
Now we are ready to state the main result of the section.

\medskip\noindent
{\bf 3.8 Proposition.}
{\it Let $X$ be a finite topological graph distinct from a singleton and from a doubleton,
and let $Y_1,\dots,Y_m$ be
the (possibly empty) 
family of homeomorphism types of all blocks of $X$.
Then
$$ 
{\cal X}^r(X)\cong\widetilde\sqcup({\cal X}^r(Y_1),\dots,{\cal X}^r(Y_m)).
$$}

Our proof of Proposition 3.8 is based on the following.

\medskip\noindent
{\bf 3.9 Lemma.}
{\it Suppose that for a finite topological graph $X$ one of the following
conditions holds:}
\item{(a)} {\it $X$ is the disjoint union of its
nonempty subgraphs $X_1$ and $X_2$, at least one of which is not a singleton;}
\item{(b)} {\it $X$ is connected and 
has a splitting with parts $X_1$ and $X_2$.}

\noindent
{\it Then}
$$
{\cal X}^r(X)\cong\widetilde\sqcup({\cal X}^r(X_1),{\cal X}^r(X_2)).
$$

We now give a proof of Proposition 3.8, assuming Lemma 3.9,
and we postpone the proof of Lemma 3.9 until the last part of the section.

\medskip
\noindent
{\bf Proof of Proposition 3.8 (assuming Lemma 3.9):}
Suppose that $Z$ is a cycle in $X$, i.e. a subgraph homeomorphic to the circle,
and let $X_0$ be the connected component of $X$ containing $Z$.
Observe that for any split decomposition of $X_0$, 
$Z$ is contained in precisely one of the
factors of this splitting. 
By Lemma 3.4,
$Z$ is then contained in a block of $X_0$, which is also a block
of $X$. In particular, this shows that
$X$ contains some block if and only if it contains a cycle.

Consider first the less interesting case when the family of blocks
of $X$ is empty. By the discussion of the previous paragraph, this happens
exactly when $X$ contains no cycle. By Lemma 2.4(2), the reflection
tree ${\cal X}^r(X)$ is then homeomorphic to the Cantor set $C$.
Since, by our convention, $C$ is the dense amalgam of the empty family
of spaces, the proposition follows in this case for ``conventional'' reasons.

We now pass to the real case of interest, when $X$ has at least one 
block, i.e. the parameter $m$ from the statement satisfies $m\ge1$.
Let $E_1,\dots,E_k$ be the family of all trivial factors of some terminal
split decompositions of all componets of $X$. We then have
$k\ge0$ (i.e. the family may be empty), and each $E_j$ is homeomorphic
either to a segment or to a singleton (the latter appear only when $X$ contains
a connected component which is a singleton). 
Let $X_1,\dots,X_p$ be the family of all connected components
of $X$. By applying recursively Lemma 3.9(a) 
together with Proposition 3.7(2)
(and using parts (1) and (2) of Proposition 2.4 and our conventions
concerning the empty set in the case when all $X_i$ are singletons),
we get 
$$
{\cal X}^r(X)\cong\widetilde\sqcup({\cal X}^r(X_1),\dots,{\cal X}^r(X_p)).
$$
By referring to the fact that blocks of $X$ are exactly the nontrivial factors in terminal
split decompositions of the components $X_1,\dots,X_p$ (Corollary  3.5), 
and applying again Proposition 3.7(2),
and then Proposition 3.7(3) (to exdclude repetitions of homeomorphism
types among the amalgamated nontrivial blocks),
we get
$$
{\cal X}^r(X)\cong\widetilde\sqcup({\cal X}^r(Y_1),\dots,{\cal X}^r(Y_m),
{\cal X}^r(E_1),\dots,{\cal X}^r(E_k)).
$$
Now, each of the spaces ${\cal X}^r(E_j)$ is either homeomorphic
to the Cantor set $C$ (by Lemma 2.4(2)) or is empty (by Lemma 2.4(1)). Then, applying Proposition 3.7(4) and
one of our conventions concerning the empty set, we get
$$
{\cal X}^r(X)\cong\widetilde\sqcup({\cal X}^r(Y_1),\dots,{\cal X}^r(Y_m)),
$$
as required. This finishes the proof.

\medskip\noindent
{\bf Proof of Lemma 3.9:}
The proof needs to be split according to the following two subcases:
\item{(i)} either $X=X_1\sqcup X_2$ and no one of $X_i$ is a singleton,
or $X$ is connected and has a splitting with parts $X_1,X_2$;
\item{(ii)} $X=X_1\sqcup X_2$ and one of the subgraphs $X_i$, say $X_2$, is a singleton.

In the subcase (i) both spaces ${\cal X}^r(X_i)$ are nonempty
and thus we need to check that the space
${\cal X}^r(X)$ can be equipped with a family 
${\cal A}={\cal A}_1\sqcup{\cal A}_2$ of subsets satisfying all the conditions
of Definition 3.6, for $k=2$, with ${\cal X}^r(X_i)$ substituted for $B_i$.
We deal with this subcase carrefully below.
In the cubcase (ii) we have ${\cal X}^r(X_2)=\emptyset$, so that our
assertion reads
$$
{\cal X}^r(X)\cong\widetilde\sqcup({\cal X}^r(X_1),\emptyset)=
\widetilde\sqcup({\cal X}^r(X_1)),
$$
and thus
we need to equip
${\cal X}^r(X)$ with a family 
${\cal A}={\cal A}_1$ of subsets satisfying all the conditions
of Definition 3.6, for $k=1$, with ${\cal X}^r(X_1)$ substituted for $B_1$.
The arguments necessary to deal with this subcase are either the same or very similar to those corresponding to the subcase (i), and we omit them.

We turn to the setting of the subcase (i).
Put ${\cal A}_i$ to be the family of images of the natural embeddings
of ${\cal X}^r(X_i)$ in ${\cal X}^r(X)$ described in the proof of part (5)
of Lemma 2.4 (with $X_i$ substituted for $Y$).
Since each $X_i$ is a proper subgraph of $X$ not reduced to a vertex,
both families ${\cal A}_i$ are countable infinite. We need to check
conditions (a1)-(a5) of Definition 3.6.

\medskip\noindent
{\it Condition} (a1).

In view of the above description of the families ${\cal A}_i$ and the second
assertion of Lemma 2.4(5), to verify this condition, it remains to
show that for any $E_1\in{\cal A}_1$ and $E_2\in{\cal A}_2$ we get
$E_1\cap E_2=\emptyset$. To see this, recall that both $E_i$ have the following
description. There is an appropriate subtree $S_i\subset T_X$, and a label
preserving isomorphism $\psi_i:T_{X_i}\to S_i$, which determine the morphism
${\bf f}^{\psi_i}:{\cal S}_{X_i}\to {\cal S}_X^{S_i}$ such that
for the induced map $f^{\psi_i}$ we have 
$E_i=f^{\psi_i}(\lim_\leftarrow{\cal S}_{X_i})$.

Consider first the case when either assumption (a) of the lemma holds
(i.e. $X=X_1\sqcup X_2$) or the subtrees $S_1,S_2$ are disjoint.
Let $F$ be any finite subtree of $T_X$ which has nonempty intersection with both $S_i$. It is not hard to see that then 
$$
f_F^{\psi_1}((X_1)_{\psi_1^{-1}(F\cap S_1)}) \, \cap \,
f_F^{\psi_2}((X_2)_{\psi_2^{-1}(F\cap S_2)}) = \emptyset,
$$
and hence the image sets $E_i$ of the maps $f^{\psi_i}$ are disjoint too,
as required.

We are now left with the case when assumption (b) of the lemma holds and
the subtrees $S_i$ intersect. Denote by $v$ the common vertex of $X_1$
and $X_2$, which is an essential vertex of $X$. Since we assume,
without loss of generality, that $D_{X_i}=X_i\cap D_X$ for both $i$,
we obviously get that $D_{X_1}\cap D_{X_2}=\{v\}$.
The intersection of the subtrees $S_i$ consists then of a single edge $e$ such that $\lambda_X(e)=v$. We denote this intersection by $F_e$, and we observe
that
$$
X_{F_e}=f_{F_e}^{\psi_1}((X_1)_{\psi_1^{-1}(F_e)}) \, \sqcup \,
f_{F_e}^{\psi_2}((X_2)_{\psi_2^{-1}(F_e)}),
$$
where $\sqcup$ denotes the disjoint union.
As in the previous case, this implies that
 the image sets $E_i$ of the maps $f^{\psi_i}$ are disjoint,
 which completes the verification of condition (a1).

\medskip\noindent
{\it Condition} (a2).

We will show that for $i=1,2$ the family ${\cal A}_i$ of subsets is null,
which obviously implies the condition.
We refer to the description of the subsets $E_i\in{\cal A}_i$ as recalled above,
in the part of proof concerning condition (a1).

Given any subtree $F\in{\cal F}_X$, for a fixed $i\in\{1,2\}$ only finitely many
of the subtrees $S_i$ (out of the countable infinite family of such subtrees)
intersect $F$. Consequently, denoting by 
$\pi^X_F:\lim_\leftarrow{\cal S}_X \to X_F$ the projection canonically
associated to the system, we get that only finitely many of the images
$\pi_F^X(E_i)$ are not singletons. This easily implies that the family
${\cal A}_i=\{ E_i \}$ is null, as required.

\medskip\noindent
{\it Condition} (a3).

For any $i\in\{ 1,2 \}$,
fix a subset $E_i\in{\cal A}_i$, described as 
$E_i=f^{\psi_i}(\lim_\leftarrow{S_{X_i}})$, as explained in the
earlier parts of this proof. We need to show that any open set in
${\cal X}^r(X)=\lim_\leftarrow{\cal S}_X$ contains a point from the
complement of $E_i$. Obviously, it is sufficient to show this for subsets $U$
from some basis of topology, e.g. for subsets of form
$$
U=(\pi_F^X)^{-1}(W):F\in{\cal F}_X,
\hbox{ $W$ is a nonempty open subset of $X_F$}
$$
where $\pi_F^X:\lim_\leftarrow{\cal S}_X\to X_F$ are the canonical projections
associated to the system.

For any $F$ and $W\subset X_F$ as above, there is a vertex $t$ of $F$ such that,
if we view the blow-up space $(X_t)^\#(\lambda_X(A_tF))$
naturally as a subspace in $X_F$, then $W$ has nonempty intersection 
with this subspace. Consequently, viewing  
$D_X\setminus\lambda_X(A_tF)$ as a subset in $(X_t)^\#(\lambda_X(A_tF))$,
there is $x\in D_X\setminus\lambda_X(A_tF)$ which belongs to $W$.
Denote by $F'$ the subtree of $T_X$ equal to the union of $F$ and the edge $e$
adjacent to $t$ and labelled with $x$. Denote by $s$ the vertex of $e$ other 
than $t$, and choose any point $x'$ in the space $X_s$ distinct from $x$ and
not belonging to $X_i$
(here we canonically view $X_i$ as a subset of $X_s$, 
since the latter is a copy of $X$). Now, viewing $x'$ as a point of $X_{F'}$,
and $x$ as a point of $X_F$, we make the following easy observations.
First, since the bonding maps in the system ${\cal S}_X$ are all surjective,
it follows that the preimage $(\pi_{F'}^X)^{-1}(x')$ is nonempty. 
Since $x'\notin X_i$,
we get that this preimage is disjoint with $E_i$. Finally, since 
$\pi_{F',F}(x')=x$, we get that 
$(\pi_{F'}^X)^{-1}(x')\subset(\pi_F^X)^{-1}(x)\subset
(\pi_F^X)^{-1}(W)=U$. Consequently, $U$ contains points from
the complement of $E_i$, which completes the verification of condition (a3).

\medskip\noindent
{\it Condition} (a4).

We need to show that in any open set $U=(\pi_F^X)^{-1}(W)$ from the
basis of topology (as described in the previous part of the proof) there is
a point $p$ belonging to the union $\bigcup{\cal A}_i$.

Similarly as in the previous part of the proof, consider a vertex $t$ 
of the subtree $F$
such that $W\cap(X_t)^\#(A_tF)\ne\emptyset$, and then consider a point
$x$ from this intersection which is also a point of 
$D_X\setminus\lambda_X(A_tF)$. Again, denote by $e$ this edge of $T_X$
adjacent to $t$ for which $\lambda_X(e)=x$, and denote by $s$ the vertex
of $e$ other than $t$. Denote by $F'$ the subtree of $T_X$ being the union
of $F$ and $e$. Moreover, let $x'$ be a point of $X_s$ which, under the 
identification of $X_s$ with $X$, is distinct from $x$ and belongs to $X_1$.
We then view this $x'$ as a point of $X_{F'}$.

Consider the set $E_i\in{\cal A}_i$ which corresponds to this subtree
$S_i$ of $T_X$ which contains $s$. Consider also a label preserving
isomorphism $\psi_i:T_{X_i}\to S_i$. Note that $x'$ is in the image
of the map $f_{F'}^{\psi_i}$ from the morphism 
${\bf f}^{\psi_i}$, so that for some $x''\in(X_i)_{\psi_i^{-1}(F'\cap S_i)}$
we have $f_{F'}^{\psi_i}(x'')=x'$. 
By surjectivity of the bonding maps, there is 
$q\in(\pi^{X_i}_{\psi_i^{-1}(F'\cap S_i)})^{-1}(x'')$, and then
$p:=f^{\psi_i}(q)$ belongs to both $(\pi_{F'}^X)^{-1}(x')$ and to $E_i$.
Since $\pi_{F',F}^X(x')=x$, we get
$$
p\in(\pi_{F'}^X)^{-1}(x') \subset (\pi_F^X)^{-1}(x) \subset
(\pi_F^X)^{-1}(W)=U,
$$
which shows that $E_i$ intersects $U$, as required.

\medskip\noindent
{\it Condition} (a5).

We start with describing a family of open, closed and $\cal A$-saturated
subsets of the space ${\cal X}^r(X)=\lim_\leftarrow{\cal S}_X$.

Suppose first that $X_1,X_2$ are the factors of a splitting of $X$,
and let $v$ be the vertex of $X$ at the intersection of these factors.
Let $e$ be any edge of the tree $T_X$ labelled with $v$, and denote by $F_e$
the subtree of $T_X$ coinciding with $e$. It is not hard to see that
the space $X_{F_e}$ splits as the disjoint union of two subgraphs, which can
be canonically identified as $(X_i)_{F_e}$ (if we view $e$ as 
an edge in the labelled trees $T_{X_i}$, respectively).
Denoting by $\pi^X_{F_e}:\lim_\leftarrow{\cal S}_X\to X_{F_e}$
the map canonically associated to the inverse limit, we put
$H_i^e:=(\pi^X_{F_e})^{-1}((X_i)_{F_e})$, thus getting a partition of
${\cal X}^r(X)=\lim_\leftarrow{\cal S}_X$ into two open and closed subsets
$H_1^e,H_2^e$.

If $X=X_1\sqcup X_2$, with both summands nonempty and not singletons,
fix any point $v\in D_{X_1}$, and consider any edge $e$ and the 
related subtree $F_e$ as in the previous paragraph.
Denote by $s,s'$ the vertices of the edge $e$.
  Then the space
$X_{F_e}$ splits as the disjoint union of two subgraphs, which can
be canonically identified as $(X_1)_{F_e}$ and as
$(X_2)_s\sqcup (X_2)_{s'}$
(if we view $e$ as 
an edge in the labelled trees $T_{X_i}$, respectively). We then put 
$H_1^e:=(\pi^X_{F_e})^{-1}((X_1)_{F_e})$ and 
$H_2^e:=(\pi^X_{F_e})^{-1}((X_2)_s\sqcup(X_2)_{s'})$,
again getting a partition of
${\cal X}^r(X)=\lim_\leftarrow{\cal S}_X$ into two open and closed subsets
$H_1^e,H_2^e$.

\smallskip
We next show that the subsets $H_i^e$ as above are $\cal A$-saturated.
Actually, it is obviously sufficient to show this for the subsets $H_1^e$. 
Fix one of the subsets $H_1^e$, and 
let $E$ be an arbitrary set from the family $\cal A$. Then
$E=f^\psi(\lim_\leftarrow{\cal S}_{X_i})$ for appropriate $i\in\{ 1,2 \}$
and appropriate label preserving isomorphism $\psi:T_{X_i}\to S$,
where $S$ is a subtree of $T_X$. 
Consider first the case when $e$ is disjoint with $S$. Let $e'$ be the first
edge on the shortest path from $e$ to $S$ in $T_X$.
Obviously, $\lambda_X(e')\ne\lambda_X(e)=v$. It follows that
if $\lambda_X(e')\in X_1$ then $E\subset H_1^e$, and if 
$\lambda_X(e')\in X_2$ then $E\cap H_1^e=\emptyset$.
In the remaining case, when $e$ intersects $S$, $e$ is in fact contained in $S$.
Thus, if $i=1$, we easily see that $E\subset H_1^e$, and if $i=2$
then $E\cap H_1^e=\emptyset$. This shows that $H_1^e$ is
$\cal A$-saturated, as required.

\smallskip
We now turn to the required separation property.
Let $p=(p_F)$ and $q=(q_F)$ be two distinct points of 
${\cal X}^r(X)=\lim_\leftarrow{\cal S}_X$. Put 
$E_X^v=\{ e\in E_X:\lambda_X(e)=v \}$.
Observe that if for some $e\in E_X^v$ exactly one of the points 
$p_{F_e},q_{F_e}$ belongs to the subgraph $(X_1)_{F_e}\subset X_{F_e}$,
then the set $H_1^e$ separates $p$ from $q$. Thus, we need to prove
the following.

\smallskip\noindent
{\bf Claim.} 
{\it If for each $e\in E_X^v$ the alternative of the following two 
conditions holds:}
\item{(a)} $p_{F_e},q_{F_e}\in (X_1)_{F_e}$, {\it or}
\item{(b)} $p_{F_e}\notin (X_1)_{F_e}$ {\it and} $q_{F_e}\notin(X_1)_{F_e}$,

\noindent
{\it then there is $E\in{\cal A}$ which contains both $p$ and $q$.}

\smallskip
To prove the claim, we first make the following observation.
In its statement, we view singletons $\{t\}$ for $t\in V_X$ as
subtrees of $T_X$.

\smallskip\noindent
{\bf Subclaim 1.}
{\it There is $t\in V_X$ such that $p_{\{t\}}\ne q_{\{t\}}$.}

\smallskip
To justify the subclaim, suppose that its assertion is not true.
It is not hard to prove inductively (using the arguments similar
as in the more difficult proof of Subclaim 2 below), 
that for each $F\in{\cal F}_X$
we have then $p_F=q_F$, which yields $p=q$, 
a contradiction. 

\smallskip
Coming back to the proof of Claim, let $t$ be as in the assertion of
Subclaim 1, and let $e$ be the edge from $E_X^v$ adjacent to $t$.
We then obviously have $p_{F_e}\ne q_{F_e}$.
Suppose that condition (a) as in the claim holds, i.e. 
$p_{F_e},q_{F_e}\in (X_1)_{F_e}$. Denote by $S$ this maximal subtree
of $T_X$ with all edge labels in $D_{X_1}$, which contains $e$.
Denote by $\psi:T_{X_1}\to S$ any label preserving isomorphism.
Put $E=f^\psi(\lim_\leftarrow{\cal S}_{X_1})$ and note that 
$E\in{\cal A}_1\subset{\cal A}$.

\smallskip\noindent
{\bf Subclaim 2.}
{\it In the setting as above, we have $p,q\in E$}.

\smallskip
To get Subclaim 2, by cofinality of the subposet of ${\cal F}_X$ appearing
below, it is sufficient to show that for those $F\in{\cal F}_X$ which
contain $e$ we have
$$
p_F,q_F\in f^\psi_F((X_1)_{\psi^{-1}(S\cap F)}).\leqno{(+)}
$$
We will do this by induction. Note that 
due to condition (a) of the claim (under assumption of which we work),
assertion $(+)$ holds true for $F=F_e$. Suppose it holds for some 
$F_0$ containing $e$, and let $F_0'$ be a subtree containing $F_0$
and having one more vertex than $F_0$. We will show that $(+)$
holds true for $F=F_0'$. 

Denote by $t_0$ the vertex of $F_0'$ not contained in $F_0$, and by 
$e_0=[t_0,s_0]$
the edge connecting $t_0$ to $F_0$. Put also $x_0=\lambda_X(e_0)$.
By our assumptions, $p_{F_0}$ and $q_{F_0}$ are two distinct points
belonging to the subset 
$f^\psi_{F_0}((X_1)_{\psi^{-1}(S\cap F_0)})\subset X_{F_0}$.
As a consequence, at least one of the points $p_{F_0'}$ and $q_{F_0'}$,
say $p_{F_0'}$, is not contained in the part $X_{t_0}^\#(x_0)$
of $X_{F_0'}$. Consequently, identifying canonically $X_{t_0}$ with $X$,
we get $p_{\{t_0\}}=x_0$.

If $x_0$ happens to belong to $(X_1)_{t_0}\subset X_{t_0}$,
by the assumption of the claim we get that $q_{\{ t_0 \}}\in(X_1)_{t_0}$
as well. This obviously implies assertion $(+)$ for $F=F_0'$.
If $x_0=\lambda_X(e_0)\notin(X_1)_{t_0}$, we get that $p_{F_0'}$
cannot be contained in $X_{t_0}^\#(x_0)$, since otherwise
$p_{F_0}=\pi^X_{F_0',F_0}(p_{F_0'})$ would be equal to $x_0$
viewed as belonging to the part $X_{s_0}^\#(A_{s_0}F_0)$
of $X_{F_0}$, contradicting the assumption that
$p_{F_0}\in f_{F_0}^\psi((X_1)_{\psi^{-1}(F_0\cap S)})$.
Thus, both $p_{F_0'}$ and $q_{F_0'}$ are then not contained
in the part $X_{t_0}^\#(x_0)$, and hence they must belong to the subset
$f_{F_0'}^\psi((X_1)_{\psi^{-1}(F_0'\cap S)})\subset X_{F_0'}$,
as required. This completes the proof of Subclaim 2.

To conclude the proof of Claim, it remains to consider the case
when condition (b) of the claim holds. In this case, an argument similar 
as above shows that for some $E\in{\cal A}_2$ we have $p,q\in E$.
We omit further details, hence finishing the proof of Claim,
the verification of condition (a5), and the proof of Lemma 3.9.


\magnification1200

\bigskip
\noindent
{\bf 4. Violated reflection trees of graphs.} 

\medskip
In this section we show that reflection trees of graphs allow a bit more
flexible description than the one given in Section 2. We will need this
description in the later arguments in the paper, notably in the proof
of the main result of the next section (Proposition 5.5) which is then used
in the proof of the main result of the paper - Theorem 7.12.

Let $X$ be a finite topological graph. A {\it violated reflection inverse system
for $X$} is any inverse system ${\cal S}^v_X$ described as follows.
Fix a countable dense subset $D$ of $X$ containing all essential vertices.
Let $T$, $V_T$, $E_T$, $\cal F$ and $\lambda:E_T\to D$ be as described
in Section 2. To each vertex $t\in V_T$ associate a copy $X_t$
of the graph $X$. Denote by ${\cal E}_T$ the set of all oriented edges
of $T$, and for each $\epsilon\in{\cal E}_T$ denote by $\alpha(\epsilon),
\omega(\epsilon)$ the initial and the terminal vertex of $\epsilon$,
respectively, and by $|\epsilon|$ the underlying unoriented edge of $T$.
Fix any family $\beta_\epsilon:\epsilon\in{\cal E}_T$ of maps
$P_{\lambda(|\epsilon|)}\to P_{\lambda(|\epsilon|)}$
(where the source and the target  blow-up locus 
$P_{\lambda(|\epsilon|)}$ is viewed  
as a subset in any appropriate blow-ups of the copies $X_{\alpha(\epsilon)}$ 
and $X_{\omega(\epsilon)}$, respectively) such that

\item{(i)} $\beta_\epsilon$ is the identity if $\lambda(|\epsilon|)$
is an essential vertex of $X$, and it is either the identity or the transposition
of the doubleton $P_{\lambda(|\epsilon|)}$ otherwise;
\item{(ii)} for each $\epsilon\in{\cal E}_T$, 
if $\bar\epsilon$ denotes the oppositely oriented edge $\epsilon$,
we have $\beta_{\bar\epsilon}=\beta_\epsilon^{-1}$.
\smallskip
\noindent
The underlying poset of the now being described system ${\cal S}^v_X$
is $\cal F$, and its spaces $Z_F:F\in{\cal F}$ are as follows.
Fix a subtree $F\in{\cal F}$, and use the notation $V_F$, $E_F$ and
$A_tF:t\in V_F$ as in Section 2. Denote also by ${\cal E}_F$
the set of all oriented adges of $F$. Put
$$
Z_F:=\bigsqcup_{t\in V_T}X_t^\#(\lambda(A_tF))/\approx,
$$
where the equivalence relation $\approx$ is induced by the following equivalences:
for each $\epsilon\in{\cal E}_F$ and each $p\in P_{\lambda(|\epsilon|)}$
identify $p\in P_{\lambda(|\epsilon|)}\subset 
X^{\#}_{\alpha(\epsilon)}(\lambda(A_{\alpha(\epsilon)}F))$ with
$\beta_\epsilon(p)\in P_{\lambda(|\epsilon|)}\subset 
X^{\#}_{\omega(\epsilon)}(\lambda(A_{\omega(\epsilon)}F))$.
The maps $\pi^v_{F',F}:Z_{F'}\to Z_F$ of the system ${\cal S}^v_X$
are then defined in the same way as the corresponding maps $\pi_{F',F}$
of the system ${\cal S}_X$, see Section 2.

We now state, as a proposition below, the main result of this section.
The assumptions which we put on $X$ are not essential,
but they simplify the argument, and correspond to the case that we need
for our later purposes in the paper. 

\medskip\noindent
{\bf 4.1 Proposition.}
{\it Let $X$ be a finite connected topological graph without
essential loop edges, and 
let ${\cal S}^v_X$ be any violated reflection inverse system for $X$.
Then $\lim_\leftarrow{\cal S}^v_X\cong{\cal X}^r(X)$.}

\medskip
The remaining part of the section is devoted to the proof of Proposition 4.1.
Choose a function $\kappa$ which to any essential edge $a$ of $X$
associates an arbitrarily selected essential vertex adjacent to $a$.
Denote by $X^2_a$ the graph obtained by gluing the two copies of
the blow-up $X^\#(\kappa(a))$ through the identity map along the 
blow-up locus $P_{\kappa(a)}$. Denote also by $D^2_a\subset X^2_a$
the union of two copies of the set 
$D\setminus\{ \kappa(a) \}\subset X^\#(\kappa(a))$.
Denote by $\tilde e$ this essential edge of $X^2_a$ which is made
of the two copies of this edge of $X^\#(\kappa(a))$ which 
corresponds to the edge $a$ of $X$ (we denote this edge of 
$X^\#(\kappa(a))$ also by $a$). To prove Proposition 4.1, we will need
the following easy to observe result, whose proof we omit.

\medskip\noindent
{\bf 4.2 Lemma.}
{\it For any points $x_1,x_2\in D^2_a\cap\tilde e$ which are not 
essential vertices
of $X^2_a$ (i.e. no one of them is an endpoint of $\tilde e$), 
including the case $x_1=x_2$, consider the blow-ups $(X^2_a)^\#(x_i)$
with their corresponding blow-up loci $P_{x_i}$, which are doubletons.
Then for any bijection $\alpha:P_{x_1}\to P_{x_2}$
there is a homeomorphism
$h:(X^2_a)^\#(x_1)\to (X^2_a)^\#(x_2)$ such that:}
\item{(1)} {\it $h$ restricted to the complement of 
the union of two edges in 
$(X^2_a)^\#(x_1)$ intersecting the blow-up locus  $P_{x_1}$ 
coincides either with the identity or with the relevant restriction of the
canonical involution of $X^2_a$
which exchanges identically the two copies of $X^\#(\kappa(a))$ contained in $X^2_a$;}
\item{(2)} {\it $h$ restricted to $P_{x_1}$ coincides with $\alpha$;}
\item{(3)} {\it $h(D^2_a\setminus\{ x_1 \})=D^2_a\setminus\{ x_2 \}$.} 


\medskip
We now turn to the essential part of the proof of Proposition 4.1.
We will show that for some appropriately chosen cofinal subposets
${\cal G},{\cal G}'\subset{\cal F}$ the restricted inverse systems 
${\cal S}_X|_{\cal G}$ and ${\cal S}^v_X|_{{\cal G}'}$ are isomorphic,
thus yielding 
$$
\lim_\leftarrow{\cal S}^v_X=\lim_\leftarrow{\cal S}^v_X|_{{\cal G}'}\cong
\lim_\leftarrow{\cal S}_X|_{{\cal G}}=\lim_\leftarrow{\cal S}_X=
{\cal X}^r(X).
 $$
Actually, we will find appropriate $\cal G$ and ${\cal G}'$ in the form of 
increasing sequences $F_n:n\ge 0$ and $F_n':n\ge0$
of finite subtrees of $T$ such that 
$\bigcup F_n=\bigcup F_n'=T$. Thus, the restricted to $\cal G$
and ${\cal G}'$ 
inverse systems mentioned above
will be in fact inverse sequences.   

We start with describing an appropriate sequence $F_n$.
Put $F_0$ to be any single vertex of $T$. Recursively, having described $F_n$,
we choose a vertex $w\in V_T$ not contained in $F_n$ and such that
it is connected to $F_n$ by a single edge of $T$, which we denote $e$.
If $\lambda(e)$ is an essential vertex of $X$, we put $F_{n+1}$ to be
the subtree of $T$ spanned on $F_n\cup\{ w \}$. If $\lambda(e)$
is not an essential vertex of $X$, let $a_e$ be the essential edge of $X$
containing $\lambda(e)$ in its interior. Let $e'$ be the edge of $T$
adjacent to $w$ and such that $\lambda(e')=\kappa(a_e)$, and let
$z$ be the vertex of $e'$ other than $w$. 
We then put $F_{n+1}$ to be the subtree of $T$ spanned on
$F_n\cup\{ w,z \}$.
It is clear that choices of the vertices $w$ at corresponding steps of the
above construction can be made so that $\bigcup F_n=T$, and
after making these choices in such a way we get a sequence
${\cal G}=(F_n)_{n\ge0}$ as required.

We now turn to describing inductively subtrees $F_n'$ and homeomorphisms
$g_n:X_{F_n}\to Z_{F_n'}$ which commute with the bonding maps
$\pi_{F_i,F_j}$ and $\pi^v_{F_i',F_j'}$ in both inverse sequences
${\cal S}_X|_{\cal G}$ and ${\cal S}^v_X|_{{\cal G}'}$, respectively.
We start with putting $F_0'=F_0$ and, under canonical identification
of both spaces $X_{F_0}$ and $Z_{F_0'}$ with $X$, we take $g_0$
to be the identity of $X$.

To proceed with the induction, suppose that for some $n\ge0$ and for all
$i=0,\dots,n$ subtrees $F_i'$ and homeomorphisms 
$g_i:X_{F_i}\to Z_{F_i'}$ have been already defined, and they satisfy the following properties:
\item{(a)} $F_0'\subset\dots\subset F_n'$;
\item{(b)} the commutativity equation $g_j\pi_{F_i,F_j}=\pi^v_{F_i',F_j'}g_i$
holds for all $0\le j<i\le n$;
\item{(c)} for $x\in X_{F_n}$ we have that    
$x\in D\setminus\lambda(A_tF_n)\subset X_t^\#(\lambda(A_tF_n))
\subset X_{F_n}$ for some $t\in V_{F_n}$ if and only if 
there is $s\in V_{F_n'}$ such that
$g_n(x)\in D\setminus\lambda(A_sF_n')\subset X_s^\#(\lambda(A_sF_n'))
\subset Z_{F_n'}$;
\item{(d)} if $x$ as  in (c) is an essential vertex of $X$ then $x$ and $g_n(x)$,
viewed as points of $X$, coincide;  
\item{(e)} if $x$ as in (c) is not an essential vertex of $X$ then $g_n(x)$
is also not an essential vertex of $X$; moreover, $x$ and $g_n(x)$
(viewed as points of $X$) belong to the same essential edge of $X$.

\noindent
Note that for $n=0$ (and thus for the already defined subtree $F_0'$
and map $g_0$) properties (a)-(e) are indeed satisfied.
To describe the subtree $F_{n+1}'$ and the map 
$g_{n+1}:X_{F_{n+1}}\to Z_{F_{n+1}'}$, suppose first that,
in the description of the sequence $F_k:k\ge0$, 
the subtree $F_{n+1}$
is spanned on $F_n\cup\{ w \}$, i.e. $\lambda(e)$ is an essential vertex 
of $X$.
Denote by $t$ the vertex of $F_n$ adjacent to $e$ and view $x=\lambda(e)$
as a point in $X_t$, and hence also in $X_{F_n}$.
Let $s\in V_{F_n'}$ be this vertex for which, due to property (c) above,
$g_n(x)\in D\setminus\lambda(A_sF_n')\subset X_s^\#(\lambda(A_sF_n'))
\subset Z_{F_n'}$. Let $d$ be this edge of $T$ adjacent to $s$
for which $\lambda(d)=g_n(x)$ (where $g_n(x)$ is viewed here as
an element of $D\subset X$). Note that, by property (d), we have
$x=g_n(x)$ (as points of $X$), and hence $\lambda(e)=\lambda(d)$
(as points of $D$). Let $u$ be the vertex of $d$ other than $s$.
Note that then $u\notin V_{F_n'}$ and put $F_{n+1}'$ to be the subtree
of $T$ spanned on $F_n'\cup\{ u \}$. Define a map
$g_{n+1}:X_{F_{n+1}}\to Z_{F_{n+1}'}$ as follows.
Viewing $\lambda(e)$ again as a point of $X_{F_n}$ (as above),
and at the same time as a point of $X_w$, we have 
$X_{F_{n+1}}=X_{F_n}^\#(\lambda(e))\cup X_w^\#(\lambda(e))$,
where the intersection corresponds to the blow-up loci $P_{\lambda(e)}$
in both summands, which are identified with each other by the identity map.
Similarly, we have 
$Z_{F_{n+1}'}=Z_{F_n'}^\#(\lambda(d))\cup X_u^\#(\lambda(d))$,
with the subsets $P_{\lambda(d)}$ in both summands glued to each other
also by the identity map (due to condition (i) in the description of 
a violated reflection inverse system for $X$). 
Let $\hat g_n:X_{F_n}^\#(\lambda(e))\to Z_{F_n'}^\#(\lambda(d))$
be the map induced by $g_n$ 
(which makes sense since $g_n(\lambda(e))=\lambda(d)$).
Put $g_{n+1}$ to be equal to $\hat g_n$ on the part $X_{F_n}^\#(\lambda(e))$,
and put it to be equal to the identity map 
$X_w^\#(\lambda(e))\to X_u^\#(\lambda(d))$
(after identifying canonically both $X_w$ and $X_u$ with $X$,
and remembering that $\lambda(e)=\lambda(d)$).

We now describe the subtree $F_{n+1}'$ and the map 
$g_{n+1}:X_{F_{n+1}}\to Z_{F_{n+1}'}$ in the second case,
when $F_{n+1}$ is spanned on $F_n\cup\{ w,z \}$, i.e. $\lambda(e)$
is not an essential vertex of $X$.
As in the previous case, denote by $t$ the vertex 
of $F_n$ adjacent to $e$ and view $x_1=\lambda(e)$
as a point in $X_t$, and hence also in $X_{F_n}$.
Let $s\in V_{F_n'}$ be this vertex for which, due to property (c) above,
$x_2:=g_n(x_1)\in D\setminus\lambda(A_sF_n')\subset X_s^\#(\lambda(A_sF_n'))
\subset Z_{F_n'}$. Let $d$ be this edge of $T$ adjacent to $s$
for which $\lambda(d)=x_2$ (where $x_2$ is viewed here as
an element of $D\subset X$). Note that, by property (e), we have that
the points $x_1,x_2$ (viewed as points of $X$) belong to the
interior of the same essential edge $a_e$ of $X$. 
Let $u$ be the vertex of $d$ other than $s$.
Note that then $u\notin V_{F_n'}$.
Let $d'$ be the edge of $T$ adjacent to $u$ and such that 
$\lambda(d')=\kappa(a_e)$, and let $v$ be the vertex of $d'$
other than $u$. 
Put $F_{n+1}'$ to be the subtree
of $T$ spanned on $F_n'\cup\{ u,v \}$.
Define a map
$g_{n+1}:X_{F_{n+1}}\to Z_{F_{n+1}'}$ as follows.
Recall that we denote by $e'$ the edge connecting the vertices $w$ and $z$
of the subtree $F_{n+1}$. Denote by $F_{e'}$ and $F_{d'}$ the subtrees
of $T$ consisting of single edges $e'$ and $d'$, respectively.
View $x_1=\lambda(e)$ again as a point of $X_{F_n}$ (as above),
and at the same time as a point of $X_w$, and hence also of $X_{F_{e'}}$.
We then have 
$X_{F_{n+1}}=X_{F_n}^\#(x_1)\cup X_{F_{e'}}^\#(x_1)$,
where the intersection corresponds to the blow-up loci $P_{x_1}$
in both summands, which are identified with each other by the identity map.
Similarly, we have 
$Z_{F_{n+1}'}=Z_{F_n'}^\#(x_2)\cup Z_{F_{d'}}^\#(x_2)$,
with the subsets $P_{x_2}$ in both summands glued to each other
by a map $\beta$ which is either the identity or the transposition (due to condition (ii) in the description of
a violated reflection inverse system for $X$). 
Note that since both $x_1,x_2$ (viewed as points of $X$) belong
to the interior of the same essential edge $a_e$ of $X$, and since 
$\lambda(e')=\lambda(d')=\kappa(a_e)$, we can canonically identify
the spaces $X_{F_{e'}}^\#(x_1)$ and $Z_{F_{d'}}^\#(x_2)$
with the spaces $(X^2_{a_e})^\#(x_1)$ and $(X^2_{a_e})^\#(x_2)$
(in the notation used in Lemma 4.2), respectively. 
Let $h:X^\#_{F_{e'}}(x_1)\to Z^\#_{F_{d'}}(x_2)$ be a homeomorphism
as in Lemma 4.2, for  $\alpha=\beta\hat g_n|_{P_{x_1}}$
(under the above mentioned identifications of the source and the target with the appropriate
blow-ups of $X^2_{a_e}$).
Let $\hat g_n:X_{F_n}^\#(x_1)\to Z_{F_n'}^\#(x_2)$ be the map 
induced by $g_n$ (which makes sense since $g_n(x_1)=x_2$).
Put $g_{n+1}$ to be equal to $\hat g_n$ on the part $X_{F_n}^\#(x_1)$,
and put it to be equal to $h$ on the remaining part $X^\#_{F_{e'}}(x_1)$.

We skip a straightforward verification that in both cases above
conditions (a)-(e) are still satisfied, for $n$ replaced with $n+1$.
We also note that, due to condition (c) kept during the inductive construction,
we have $\bigcup F_n'=T$, and hence the subposet 
${\cal G}'=(F_n')_{n\ge0}$ is cofinal in $\cal F$. Moreover, the commuting
sequence of homeomorphisms $(g_n)$ obviously yields an isomorphism
of the restricted inverse systems ${\cal S}_X|_{\cal G}$ and
${\cal S}^v_X|_{{\cal G}'}$, as required.
This finishes the proof of Proposition 4.1.


\magnification1200

\bigskip\noindent
{\bf 5. Inverse sequences of $ X$-graphs and $ X$-blow-ups.}

\medskip
In this section, for any topological graph $X$ without loop edges
and without vertices of degree 1,
we describe a usefull class of inverse sequences of topological graphs,
called {\it null and dense inverse sequences of $X$-graphs
and $X$-blow-ups}, such that for any inverse sequence $\cal R$ in this class
we have $\lim_{\leftarrow}{\cal R}\cong {\cal X}^r(X)$. One may think
of these sequences as of another (quite flexible and convenient) way
of expressing the spaces ${\cal X}^r(X)$. We will use such sequences as a tool
in the proof of the main
result of this paper (Theorem 7.12),
to recognize visual boundaries of some Coxeter systems as (being
homeomorhic to) the trees of graphs ${\cal X}^r(X)$.

\smallskip
Let $ X$ be a (finite) topological graph, equipped with its natural stratification
into essential vertices and edges.
For all of this section we assume that $X$ contains no essential loop edges
(but it may have multiple edges), and no vertices of degree 1.
Let $V,E$ be the sets of (essential) vertices and edges of $ X$, respectively.
For $v\in V$, let $ X_v$ be the {\it link of $ X$ at $v$},
i.e. the set of edges issuing from $v$. 
For each $v\in V$ consider 
the induced labelling of the elements of the link $ X_v$ with
the elements of $E$.

\medskip\noindent
{\bf 5.1 Definition.} Given a graph $ X$ as above,
with the sets $V,E$ of vertices and edges,
an $ X$-{\it graph} is a finite topological graph $X'$
having no esential loop edges, 
equipped with labellings of its vertices and edges with
elements of the sets $V,E$ respectively, such that
if $v'$ is a vertex
of $ X'$ labelled with $v\in V$ then there is a bijection
$\beta_{v'}: X'_{v'}\to X_v$ 
preserving the induced labels from the set $E$.
Clearly, $ X$ itself is tautologically an $ X$-graph.
Note that bijections $\{ \beta_{v'} \}$ are in fact unique,
and we view them as part of the structure of an $ X$-graph.

\medskip
We now describe, in Definitions 5.2 and 5.3, 
a class of maps between $ X$-graphs
that we call {\it $ X$-blow-ups}.

\medskip\noindent
{\bf 5.2 Definition.} Let $ X'$ be an $ X$-graph and $v'$
a vertex of $ X'$ labelled with $v\in V$. {\it The $X$-blow-up at $v'$}
is a pair $( X'',f)$ consisting of an $ X$-graph $ X''$
and a map $f: X''\to X'$ defined as follows:
\item{$\bullet$} delete $v'$ from $ X'$ and complete the resulting
space naturally so that to distinct edges issuing from $v'$
distinct points are attached in place of $v'$; label the attached points $p$
naturally with elements of the link $ X'_{v'}$ and denote these labels by 
$\lambda(p)$ (this labelling plays only auxiliary role);
delete from $ X$ the open star of the vertex $v$, i.e. $v$
and the interiors of the edges issuing from $v$; 
glue the above modified $X'$ to the above modified $X$ through the map,
which to any point $p$ attached to $ X'\setminus\{ v' \}$ associates the terminal
vertex (other than $v$) of this edge deleted from $X$ which coincides with
the label $\lambda(p)$ (this gluing map is not necessarily a bijection); view the obtained
space as a topological graph (with natural stratification), 
and denote it by $ X''$;
consider also the
obvious labelling of vertices and edges of $ X''$, induced
from the corresponding labellings in $ X'$ and $ X$, 
and note that with this labelling $ X''$ is
an $ X$-graph;
\item{$\bullet$} define  $f: X''\to X'$ to be
the identity on the part of $ X''$ corresponding to
$ X'\setminus\{ v' \}$, and to be the constant map with value $v'$
on the remaining part of $ X''$. 

\noindent
If $f:X''\to X'$ is an $X$-blow-up at a vertex $v'$ of $X'$, we call $v'$
{\it the blow-up point} of $f$.

\medskip\noindent
{\bf 5.3 Definition.}  Let $ X'$ be an $ X$-graph,
$e'$ an edge of $ X'$ labelled with $e\in E$, and let $L$
be a (closed) segment contained in the interior of $e'$. 
The {\it $X$-blow-up at $L$}
is a pair $( X'',f)$ consisting of an $ X$-graph $ X''$
and a map $f: X''\to X'$ defined as follows:
\item{$\bullet$} delete the interior of $L$ 
from $ X'$ and the interior of $e$ from $ X$; 
choose any bijection between the sets of endpoints of $L$
and of $e$ and glue the remaining parts of $ X'$ and $ X$
accordingly with this bijection; view the obtained
space as a topological graph (with natural stratification), 
and denote it by $ X''$;
consider the obvious labelling of vertices and edges of $ X''$ 
induced from those in $ X'$ and $ X$ and note that 
with this labelling $ X''$ is an $ X$-graph;

\item{$\bullet$} define  $f: X''\to X'$ to be
the identity on the part of $ X''$ corresponding to
$ X'\setminus\hbox{int}(L)$;
if $\tilde e$ is an edge of $ X$ distinct
from $e$ and connecting the endpoints of $e$, let $f$ map
$\tilde e$ homeomorphically on $L$, consistently with the identification 
of the endpoints;
to define  $f$ on the remaining part of $ X''$, choose a point $p$
in the interior of $L$ and consider two half-segments into which
$p$ splits $L$; 
for any edge $\varepsilon$ issuing from
an endpoint of $e$ and terminating not at the other endpoint of $e$,  
let $f$ map $\varepsilon$ homeomorphically on the appropriate
half-segment of $L$;
finally, let $f$ map the remaining part of $ X''$ to $p$.  

\noindent
If $f:X''\to X'$ is an $X$-blow-up at a segment $L$, as above, we call $L$
{\it the blow-up segment} of $f$.

\smallskip
Note that in fact we have two possible operations of $ X$-blow-up 
at a segment $L$,
related to the two choices of the gluing bijection in the above description.
These two operations are in general essentially distinct.

Note also that, by the assumption that $X$ has no vertex of degree 1,
the $X$-blow-up maps at segments are surjective, and thus all
$X$-blow-up maps are surjective.

\medskip
To unify terminology, we use the term {\it blow-up-locus} for both blow-up 
points and blow-up segments.
If $f: X''\to X'$ is an $ X$-blow-up with blow-up locus $\Lambda$, 
we view all points in $ X'\setminus\Lambda$
as points of $ X''$, and we say that these are
the points of $ X'$ {\it unaffected by $f$}.

\medskip
In the main result of this section, Proposition 5.5, we deal
with inverse sequences 
$$
{X_1}\longleftarrow\hskip-14pt^{\pi_1}X_2
\longleftarrow\hskip-14pt^{\pi_2}\dots
$$
of the form presented in the following definition.

\medskip\noindent
{\bf 5.4 Definition.}
We call an inverse sequence ${\cal R}=(\{ X_n \}_{n\ge1},\{ \pi_n \}_{n\ge1})$
a {\it null and dense inverse sequence of $ X$-graphs and $ X$-blow-ups}
if it satisfies the following conditions:
\smallskip
\itemitem{(i1)} $X_1= X$, each $X_i$ is an $ X$-graph, and
each map $\pi_{i}:X_{i+1}\to X_i$ is an $ X$-blow-up
(either at some vertex of $X_i$, or at some segment
$L$ contained in the interior of some edge of $X_i$); 
\smallskip
\itemitem{(i2)} for each $i$ every vertex $v'$ in $X_i$ eventually blows-up;
more precisely, there is
$k\ge i$ such that $v'$ is unaffected by the maps 
$\pi_{i},\pi_{i+1},\dots,\pi_{k-1}$ (i.e. it is naturally a vertex
of $X_k$) and the map
$\pi_{k}:X_{k+1}\to X_k$ is the $ X$-blow-up at $v'$.
\smallskip
\itemitem{(i3)} for each $i$ the family of the images in $X_i$ 
(through appropriate compositions of the maps $\pi_k$)
of all blow-up segments of the maps $\pi_m$ 
with $m>i$, is  null, i.e. the 
diameters of those images converge to 0; more precisely, after choosing any
metric in $X_i$, for any $\epsilon>0$ the set of indices $m$ such that
$m>i$, $\pi_m$ is a blow-up at a segment, 
and the image of this blow-up segment of $\pi_m$
in $X_i$ has diameter greater than $\epsilon$, is finite;
\smallskip
\itemitem{(i4)} for each $i$ the union of the images in $X_i$ of all 
blow-up segments of the maps $\pi_m$ with $m>i$ is dense in $X_i$. 
\smallskip

\medskip\noindent
{\bf 5.5 Proposition.} {\it Let $ X$ be a finite topological graph whose
natural stratification contains no loops and no vertices of degree 1, 
and let $\cal R$ be any 
inverse sequence 
satisfying conditions (i1)--(i4), i.e. a null and dense inverse sequence
of $X$-graphs and $X$-blow-ups. Then the inverse limit
$\lim_\leftarrow{\cal R}$ is homeomorphic to the reflection tree
${\cal X}^r( X)$.}

\medskip
Before starting the proof of Proposition 5.5, we
present some fairly general result which describes
a method of modifying an inverse sequence of metric compacta without
affecting its limit. The modification consists of quotiening the spaces
in the sequence subject to their appropriate partitions. We call 
the sequence resulting from any such modification a {\it good quotient}
of the initial sequence. We will use a modification of this kind as a tool in our proof of Proposition 5.5.

Recall that an {\it inverse sequence of compact metric spaces}
is a tuple ${\cal P}=(\{ P_i \}, \{ \pi_i \})$, where each $P_i$ is 
a compact metric  space and each 
$\pi_i:P_{i+1}\to P_i$ is a continuous map. For each $i>j$ we denote
by $\pi_{i,j}:P_i\to P_j$ the composition 
$\pi_{i,j}:=\pi_j\circ\dots\circ\pi_{i-1}$. We also use the convention
that for each $i$ the map $\pi_{i,i}$ is the identity of $P_i$.
Recall  that the limit $\lim_\leftarrow{\cal P}$ is, by definition,
a subspace 
in the product $\prod_{i}P_i$, with the induced product topology, 
consisting of the sequences
$(x_i)\in\prod_iP_i$ such that $\pi_i(x_{i+1})=x_i$ for each $i$.
Sequences $(x_i)$ as in the previous sentence are called {\it threads} of $\cal P$. 

We refer the reader to [Dav] for the introduction and basic results
concerning partitions of the spaces (called also decompositions)
and their associated quotients.

\medskip\noindent
{\bf 5.6 Lemma.} {\it Let ${\cal P}=(\{ P_i \}_{i\ge1}, \{ \pi_i \}_{i\ge1})$
be an inverse sequence of compact metric spaces, and for each $i\ge1$
let ${\cal Q}_i$ be an upper-semicontinuous partition of $P_i$. Suppose that the 
following further conditions are satisfied:}

\itemitem{(p1)} {\it (compatibility) for each $i\ge1$ and for each 
$Q\in{\cal Q}_{i+1}$
there is $Q'\in{\cal Q}_i$ such that $\pi_i(Q)\subset Q'$;}

\itemitem{(p2)} {\it (fineness) for any sequence $(Q_i)_{i\ge k}$ of subsets
$Q_i\in{\cal Q}_i$ such that $\pi_i(Q_{i+1})\subset Q_i$ for each $i\ge k$,
the diameters of the images $\pi_{i,k}(Q_i):i\ge k$ converge to $0$.}


\noindent
{\it Then the induced inverse sequence of quotients
$$
{\cal P^*}=( \{ P_i/{\cal Q}_i \}_{i\ge1}, 
\{ \pi^*_i:P_{i+1}/{\cal Q}_{i+1}\to P_i/{\cal Q}_i \}_{i\ge1})
$$
is well defined, consists of compact metric spaces and continuous maps, and
$$
\lim_{\longleftarrow}{\cal P^*}\cong\lim_{\longleftarrow}{\cal P}.
$$}

\medskip\noindent
{\bf Proof:}
Note that, by upper semi-continuity assumption, each decomposition space
$P_i^*:=P_i/{\cal Q}_i$ is a compact metric space 
(see Theorem 2 on page 13 in [Dav]).
Denote by $\rho_i:P_i\to P_i^*$ the quotient map. Define 
$\pi_i^*:P_{i+1}^*\to P_i^*$ to be the unique map for which
$\pi_i^*\rho_{i+1}=\rho_i\pi_i$. It is well defined and continuous due to
the compatibility assumption (p1).

Observe further, that the family $(\rho_i)$ of the quotient maps yields
a morphism ${\cal P}\to {\cal P^*}$ of the inverse sequences, due to the commutativity
relations $\pi_i^*\rho_{i+1}=\rho_i\pi_i$ mentioned in the previous paragraph.
This morphism induces then a continuous map 
$$\rho:\lim_{\longleftarrow}{\cal P}\to\lim_{\longleftarrow}{\cal P^*}$$
given by $\rho((x_i)):=(\rho_i(x_i))$.
We will show that $\rho$ is a homeomorphism. Since both limits are compact,
it suffices to show that $\rho$ is injective and surjective.

To prove injectivity of $\rho$, consider $z,w\in\lim_{\leftarrow}{\cal P}$,
$z\ne w$, $z=(z_i)$ and $w=(w_i)$. 
We claim that there is $k$ such that $z_k$ and $w_k$ are not contained both
in the same $Q\in{\cal Q}_k$. Suppose this is not true. Then for each $i\ge1$
there is a set $Q_i\in{\cal Q}_i$ with $z_i,w_i\in Q_i$.
By compatibility condition (p1), we then have the inclusions $\pi_i(Q_{i+1})\subset Q_i$
for all $i$. By condition (p2), for each $k\ge1$ we have
$$
\lim_{i\to\infty}\hbox{diam}(\pi_{i,k}(Q_i))=0. 
\leqno(+)
$$
It follows from our a contrario assumption that for each $i\ge k$ both
$z_k=\pi_{i,k}(z_i)$ and $w_k=\pi_{i,k}(w_i)$ belong to the set $\pi_{i,k}(Q_i)$.
As a consequence of $(+)$ we then have $z_k=w_k$ for all $k$,
contradicting the assumption that $z\ne w$. Thus our claim follows.
Note however that the assertion of this claim directly means that 
$\rho_k(z_k)\ne\rho_k(w_k)$, and this obviously implies that $\rho(z)\ne\rho(w)$,
as required.

To prove surjectivity of $\rho$, consider any 
$x^*=(x_i^*)\in\lim_{\leftarrow}{\cal P^*}$.
We need to construct a thread $(x_i)$ for $\cal P$ such that 
$\rho_i(x_i)=x_i^*$ for each $i$. For any $i\ge1$ put 
$Q_i=\rho_i^{-1}(x_i^*)$, and note that $Q_i\in{\cal Q}_i$
and $\pi_i(Q_{i+1})\subset Q_i$ for all $i$. Note also that, 
for each $k\ge1$ the sequence of subsets $(\pi_{i,k}(Q_i))_{i\ge k}$ of $P_k$ 
is nested, i.e. $\pi_{i+1,k}(Q_{i+1})\subset\pi_{i,k}(Q_i)$, and, 
by fineness condition (p2), the diameters of these subsets converge to 0.
Hence, the intersection $\bigcap_{i\ge k}(\pi_{i,k}(Q_i))$ is a singleton,
and we take it as $x_k$. It is then straightforward to check that 
so defined $(x_i)$
is a thread, and that $\rho_i(x_i)=x_i^*$.

This completes the proof of the lemma.

\medskip\noindent
{\bf 5.7 Definition.} A {\it good quotient} of an inverse system $\cal P$ is
any inverse sequence ${\cal P^*}$ obtained from $\cal P$ as in Lemma 5.6.

\medskip\noindent
{\bf Proof of Proposition 5.5:} To prove the proposition, we will show that 
some good quotient $\cal R^*$ of the sequence $\cal R$
is isomorphic, as inverse sequence of topological spaces, to a cofinal
subsequence $\cal P$ in some violated
reflection inverse system ${\cal S}^v_X$ for $ X$.
In view of Proposition 4.1 and Lemma 5.6,  
this will imply that
$$
\lim_{\longleftarrow}{\cal R}\cong \lim_{\longleftarrow}{\cal R^*}\cong
\lim_{\longleftarrow}{\cal P}=\lim_{\longleftarrow}{\cal S}^v_X\cong
{\cal X}^r( X).
$$
To clarify the exposition, we split this rather long proof into steps.

\medskip\noindent{\bf Step 1.}

We start with some useful terminology and observations.
\smallskip
Recall that ${\cal R}=(\{ X_i \}_{i\ge1}, \{ \pi_i \}_{i\ge1})$
is a null and dense inverse sequence of $X$-graphs and $X$ blow-ups.
For any $i\ge1$, we call a point $x\in X_i$ a {\it direct blow-up point} of $X_i$
if $x$ is a blow-up point of $\pi_i$ or if for some $n\ge i+1$ 
the point $x$ is unaffected by the blow-ups $\pi_i,\pi_{i+1},\dots,\pi_{n-1}$ 
and,
when viewed as a point of $X_n$, $x$ is a blow-up point 
of $\pi_{n}$.
We define similarly a {\it direct blow-up segment} in $X_i$
as a closed segment $L$ contained in the interior of some edge of $X_i$
such that, either $L$ is the blow-up segment of $\pi_i$, 
or for some $n\ge i+1$  the blow-ups
$\pi_i,\pi_{i+1},\dots,\pi_{n-1}$ do not affect $L$ (so that it can be viewed as 
a segment contained in the interior
of some edge of $X_n$), and $L$ is a blow-up segment of $\pi_n$.
We use the term {\it direct blow-up locus} to denote both any
direct blow-up point and any direct blow-up segment.
It is not hard to make the following observations,
which hold for each $i\ge1$:
\smallskip
\itemitem{(d1)} if $\Lambda$ is the blow-up lucus of $\pi_i$,
then its preimage $\pi_i^{-1}(\Lambda)$ in $X_{i+1}$, and hence also
its preimage in any $X_j$ with $j>i$, is a subgraph of $X_j$
(i.e. the union of some family of closed essential edges of $X_j$ and some
essential vertices of $X_j$);
\smallskip
\itemitem{(d2)}  if $\Lambda$ is the blow-up lucus of $\pi_i$,
then the image  in $X_i$ (through the appropriate
composition of the maps $\pi_k$) of the blow-up locus of any map $\pi_j$
with $j>i$ is either contained in or disjoint with $\Lambda$;
\smallskip
\itemitem{(d3)} the direct blow-up loci in $X_i$ 
are pairwise disjoint
(this  can be deduced e.g. from property (d2));
\smallskip
\itemitem{(d4)} the image in $X_i$ (through the appropriate
composition of the maps $\pi_k$) of the blow-up locus of any map $\pi_j$
with $j>i$, is contained in some direct blow-up locus of $X_i$
(this can be deduced e.g. from properties (d2) and (d3));
\smallskip
\itemitem{(d5)} each essential vertex of $X_i$ is a direct
blow-up point (this follows from condition (i2)), and there are no more
direct blow-up points in $X_i$;
\smallskip
\itemitem{(d6)} the union of all direct blow-up loci 
in $X_i$ is a dense subset of $X_i$ (this follows from condition (i4) in
view of the above observation (d4));
\smallskip
\itemitem{(d7)} each point in $x\in X_i$ which is not an essential vertex 
and does not belong to a direct blow-up segment of $X_i$, is not affected
by any of the $ X$-blow-ups $\pi_k$ with $k\ge i$, and thus it can be canonically
identified with a point in the limit $\lim_{\leftarrow}{\cal R}$;
we call any such point $x$ a {\it stable} point of $X_i$;
\smallskip
\itemitem{(d8)} the essential vertices of $X_i$ and the ends
of all direct blow-up segments of $X_i$ belong to the closure
of the set of all stable points of $X_i$; moreover, each stable point
is an accumulation point of the set of all stable points,
and it can be approached by the stable points from each its side in the edge
of $X_i$ in the interior of which it is contained.
\smallskip

\noindent
For each $k\ge1$, denote by $\Lambda_k$
the blow-up locus of the map $\pi_k:X_{k+1}\to X_k$ in the sequence $\cal R$.
We make the following easy observation, the proof of which we skip.

\smallskip\noindent
{\bf 5.10 Fact.}
{\it For each $k\ge1$, if $\Lambda$ is a direct blow-up locus in $X_k$
distinct from $\Lambda_k$, then}
\item{(1)} {\it the restriction $\pi_k|_{\pi_k^{-1}(\Lambda)}: 
\pi_k^{-1}(\Lambda) \to \Lambda$ is a bijection, so that we can view
$\Lambda$ as a subset in $X_{k+1}$ (by identifying it with its preimage
via the above mentioned bijection);}
\item{(2)} {\it the preimage $\pi_k^{-1}(\Lambda)$ is a direct blow-up locus
in $X_{k+1}$.}

\smallskip
We say that a direct blow-up locus of $X_{k+1}$ is {\it inherited} from $X_k$
if it has a form $\pi_k^{-1}(\Lambda)$ for some direct blow-up locus $\Lambda$
of $X_k$ distinct from $\Lambda_k$. Note that, in each of the graphs $X_{k+1}$
we have two kinds of blow-up loci $\Lambda$: 
\item{$\bullet$} the ones 
inherited from $X_k$, and
\item{$\bullet$} those for which $\pi_k(\Lambda)\subset\Lambda_k$. 

\noindent
The latter will be called {\it non-inherited} blow-up loci of $X_{k+1}$.
Locus $\Lambda_{k+1}$
can be of any of the above two kinds.
We artificially view all the blow-up loci of $X_1$ as inherited.

\medskip\noindent{\bf Step 2.}

In this step we describe some good quotient $\cal R^*$ of the sequence $\cal R$,
by choosing appropriate partitions ${\cal Q}_i$ of the spaces $X_i$
in $\cal R$. We also verify the assumptions (p1) and (p2) of Lemma 5.6
for the chosen partitions.

We start with choosing a partition ${\cal Q}_1$ of $X_1$.
Let ${\cal J}_1$ be a family consisting of closed normal neighbourhoods
of vertices of $X_1$ and of closed segments contained in the interiors of
edges of $X_1$, such that
\smallskip
\itemitem{(j1)} the elements in ${\cal J}_1$ are pairwise disjoint,
\smallskip
\itemitem{(j2)} for any direct blow-up locus $\Lambda$ of $X_1$ there is 
$J\in{\cal J}_1$ such that $\Lambda\subset\hbox{int}(J)$.
\smallskip
\itemitem{(j3)} for each 
$J\in{\cal J}_1$ the boundary $\partial J$ consists of stable points of $X_1$.
\smallskip

\noindent
Such a family ${\cal J}_1$ can be constructed recursively, as follows:
\smallskip
\item{$\bullet$} 
put $J_{1}$ to be some closed 
normal neighbourhood in $X_1$
of the blow-up locus $\Lambda_1$  (the blow-up locus of $\pi_1$), 
such that the boundary $\partial J_{1}$
consists of stable points of $X_1$; 
\smallskip
\item{$\bullet$} having already defined the neighbourhoods
$J_{1},\dots,J_{n}$, let $i_n$ 
be the smallest integer such that
$\Lambda_{i_n}$ (the blow-uo locusof the map $\pi_{i_n}$) 
is a direct blow-up locus of $X_1$
not contained in $J_{1}\cup\dots\cup J_{n}$; note that, since the boundaries
of $J_1,\dots,J_n$ consist of stable points, they are disjoint with 
$\Lambda_{i_n}$, and since the latter is connected and not contained in any 
$J_i$, it is disjoint with all these $J_i$; by this and by property (d8),
one may choose 
some closed normal neighbourhood 
of $\Lambda_{i_n}$ in $X_1$ disjoint with
$J_{1}\cup\dots\cup J_{n}$ and such that $\partial J_{n+1}$
consists of stable points of $X_1$; we take this neigbourhood as $J_{n+1}$. 
\smallskip

\noindent
We take as ${\cal J}_1$ the family $\{ J_i:i\ge1 \}$, 
noting that it indeed satisfies
conditions (j1)-(j3).
%
%
We then take as 
the partition ${\cal Q}_1$ of $X_1$ all the sets of ${\cal J}_1$
and all singletons in the complement of their union. Note that, since 
the diameters of the sets in the family ${\cal J}_1$ obviously converge to $0$,
the partition ${\cal Q}_1$ is upper semi-continuous 
(see Proposition 3 on page 14 in [Dav]).

\smallskip
Suppose know, that for some $k\ge1$, and for each $i\le k$, 
we have already chosen the 
family ${\cal J}_i$ of subsets $J\subset X_i$, such that the following
conditions are satisfied:

\smallskip
\itemitem{(jj0)} each $J\in {\cal J}_i$
is either a closed normal neighbourhood of a vertex of $X_i$,
or a segment contained in the interior of an edge of $X_i$, and its boundary
$\partial J$ consists of stable points of $X_i$,

\smallskip
\itemitem{(jj1)} the sets of the family ${\cal J}_i$ are pairwise disjoint,

\smallskip
\itemitem{(jj2)} any direct blow-up locus of $X_i$ is contained in one of the
sets of ${\cal J}_i$,

\smallskip
\itemitem{(jj3)} if $i\ge2$, denote by $J(\Lambda_{i-1})$ this set of the family
${\cal J}_{i-1}$ which contains $\Lambda_{i-1}$; then for any $J\in{\cal J}_i$
either $\pi_{i-1}(J)\subset \hbox{int}(J(\Lambda_{i-1}))$ or 
$\pi_{i-1}(J)\cap J(\Lambda_{i-1})=\emptyset$, and in the latter case
we have $\pi_{i-1}(J)\in{\cal J}_{i-1}$;

\smallskip
\itemitem{(jj4)} if $i\ge2$, $J\in{\cal J}_i$, 
and $\pi_{i-1}(J)\subset J(\Lambda_{i-1})$, denote by $\Lambda_J$
this direct blow-up locus of $X_i$ contained in $J$ which is the blow-up
locus of $\pi_k$ with the smallest possible $k$; then for any $j<i$ 
we have
$$
\hbox{diam}(\pi_{i,j}(J))<\max[2\cdot\hbox{diam}(\pi_{i,j}(\Lambda_J)),1/i],
$$
where $\pi_{i,j}=\pi_j\circ\dots\circ\pi_{i-1}$. 
\smallskip

\noindent
We describe a family ${\cal J}_{k+1}$ of subsets of $X_{k+1}$ which also
satisfies conditions (jj0)-(jj4), with $k+1$ substituted for $i$.
It will be a (disjoint) union of two subfamilies ${\cal J}'_{k+1}$ and 
${\cal J}''_{k+1}$ described as follows. 
As ${\cal J}'_{k+1}$ take the family of preimeges $\pi_k^{-1}(J)$ for
all $J\in{\cal J}_k$, $J\ne J(\Lambda_k)$.
Note that, since the map $\pi_k:X_{k+1}\to X_k$ is injective on the preimege
$\pi_k^{-1}(X_k\setminus\Lambda_k)$ (because $X_k\setminus\Lambda_k$
is viewed as part of $X_{k+1}$ and $\pi_k$
is viewed as the identity map on this part), this map is also injective
on any preimage 
$\pi_k^{-1}(J)$ for any $J\in{\cal J}_k\setminus\{ J(\Lambda_k) \}$.
Thus, we can identify each set $\pi_k^{-1}(J)$ 
as above with the corresponding $J$. From this point of view, ${\cal J}'_{k+1}$
can be identified with the family ${\cal J}_k\setminus\{ J(\Lambda_k) \}$.
To describe ${\cal J}''_{k+1}$, we perform a recursive construction similar
to that in the description of the family ${\cal J}_1$ above.
More precisely, the construction consists of the following steps:

\smallskip
\item{$\bullet$}
let $i_0\ge k+1$ be the smallest integer such that $\Lambda_{i_0}$
is contained in $X_{k+1}\setminus(X_k\setminus J(\Lambda_k))$; 
note that then $\Lambda_{i_0}$ is a direct blow-up locus in $X_{k+1}$;
put $J''_1$ to be a closed normal neighbourhood of $\Lambda_{i_0}$ 
in $X_{k+1}$, with $\partial J''_1$ consisting of stable points of $X_{k+1}$,
and so close to $\Lambda_{i_0}$ that
\smallskip
\itemitem{(1)} $J''_1\subset X_{k+1}\setminus(X_k\setminus J(\Lambda_k))$,
and
\smallskip
\itemitem{(2)} for any $j<k+1$ we have
$$
\hbox{diam}(\pi_{k+1,j}(J''_1))<\max[2\cdot\hbox{diam}(\pi_{k+1,j}(\Lambda_{i_0})),1/k+1];
$$

\item{$\bullet$}
having already defined the neighbourhoods
$J''_{1},\dots,J''_{n}$, let $i_n$ 
be the smallest integer such that
$\Lambda_{i_n}$ is a direct blow-up locus of $X_k+1$ 
contained in $X_{k+1}\setminus(X_k\setminus J(\Lambda_k))$, but
not contained in $J''_{1}\cup\dots\cup J''_{n}$;
then $\Lambda_{i_n}$ is actually disjoint with $J''_{1}\cup\dots\cup J''_{n}$;
put $J''_{n+1}$ to be a closed normal neighbourhood of $\Lambda_{i_n}$ 
in $X_{k+1}$, with $\partial J''_{n+1}$ consisting of stable points of $X_{k+1}$,
and so close to $\Lambda_{i_n}$ that
\smallskip
\itemitem{(1)} $J''_{n+1}\subset X_{k+1}\setminus(X_k\setminus J(\Lambda_k))$
and $J''_{n+1}\cap (J''_{1}\cup\dots\cup J''_{n})=\emptyset$, and
\smallskip
\itemitem{(2)} for any $j<k+1$ we have
$$
\hbox{diam}(\pi_{k+1,j}(J''_{n+1}))<\max[2\cdot\hbox{diam}(\pi_{k+1,j}(\Lambda_{i_n})),1/k+1].
$$

\noindent
We take as ${\cal J}''_{k+1}$ the family $\{ J''_i:i\ge1 \}$.
We skip further details.

Now, for each $i\ge2$, we define the partition ${\cal Q}_i$ of the space $X_i$
as consisting of all sets from the family ${\cal J}_i$, and all singletons
from the complement of the union of these sets. It follows fairly directly
from conditions (jj0) and (jj1) that diameters of the sets $J\in{\cal J}_i$
converge to 0 (i.e. the family ${\cal J}_i$ is null), and hence the partition
${\cal Q}_i$ is upper semi-continuous.

It follows directly from condition (jj3) that the family 
${\cal Q}_i:i\ge1$ of partitions satisfies the
compatibility condition
(p1) of Lemma 5.6.  We will show that it satisfies also condition (p2) 
of this lemma. 

\smallskip
To verify the fineness condition (p2), let $(Q_i)_{i\ge k}$ be any sequence 
as in this condition. 
It follows from the incusions $\pi_i(Q_{i+1})\subset Q_i$ that the sequence
of diameters $(\hbox{diam}(\pi_{i,k}(Q_i)))_{i\ge k}$ is nonincreasing.
Note that, if some of the sets $Q_i\in{\cal Q}_i$
is a singleton, the condition (p2) is satisfied trivially. So we assume that for each 
$i\ge k$ we have $Q_i\in{\cal J}_i$.
Under this assumption, we have the following claim.

\smallskip\noindent
{\bf Claim.}
{\it For infinitely many $i\ge k$ we have $Q_i=J(\Lambda_i)$.}

\smallskip
To prove the claim, we note first that, as a consequence of observation (d6),
for each $i\ge k$ the neighbourhood $Q_i$ contains at least one direct
blow-up locus of $X_i$. Let $n_i$ be the smallest among the numbers $n\ge i$
such that $\Lambda_n$ is a direct blow-up locus of $X_i$ contained in $Q_i$.
Then the points of $Q_i$ are successively unaffected by the maps 
$\pi_i,\pi_{i+1},\dots,\pi_{n_i-1}$, and hence we view $Q_i$ also as a subset
of $X_{n_i}$. Moreover, by the construction of the families ${\cal J}_n$,
we have $Q_i\in{\cal J}_{n_i}$ (actually, $Q_i\in{\cal J}'_{n_i}$), 
and hence $J(\Lambda_{n_i})=Q_i$.
It is also not hard to observe that $Q_i$, as a subset of $X_{n_i}$,
coincides with the preimage $\pi_{n_i,i}^{-1}(Q_i)$. It follows that
for any $J\in{\cal J}_{n_i}$, $J\ne Q_i$, the preimage $\pi_{n_i,i}^{-1}(J)$
is disjooint with $Q_i\subset X_i$. Since we have 
$\pi_{n_i,i}(Q_{n_i})\subset Q_i$, we get that $Q_{n_i}=Q_i$,
and consequently $Q_{n_i}=J(\Lambda_{n_i})$.
Thus, for each $i\ge k$ there is $n\ge i$ such that $Q_n=J(\Lambda_n)$.
This obviously implies the claim.
 
\smallskip
Denote  by $(i_m)_{m\ge1}$ the infinite increasing sequence consisting of 
all $i\ge k+1$
such that $Q_{i-1}=\Lambda_{i-1}$. Obviously, for those $i$ we also have
$\pi_{i-1}(Q_i)\subset J(\Lambda_{i-1})$.
By condition (jj4), for each $m$ we have
$$
\hbox{diam}(\pi_{i_m,k}(Q_{i_m}))<
\max[2\cdot\hbox{diam}(\pi_{i_m,k}(\Lambda_{i_m})),1/{i_m}].
\leqno{(\star)}
$$
Since, by condition (i3) in Definition 5.4, we have 
$\lim_{m\to\infty}\hbox{diam}(\pi_{i_m,k}(\Lambda_{i_m}))=0$,
we deduce from 
the estimate~$(\star)$, and from the fact that the full sequence 
$(\hbox{diam}(\pi_{i,k}(Q_i)))_{i\ge k}$
is nonincreasing, that
$\lim_{i\to\infty}\hbox{diam}(\pi_{i,k}(Q_i))=0$ as well.
Thus, the partitions ${\cal Q}_i$ satisfy the fineness condition (p2).

\smallskip
Denote by $\pi^*_i:X_{i+1}/{\cal Q}_{i+1}\to X_i/{\cal Q}_i$
the maps naturally induced by the maps $\pi_i$ on the quotients
(which are well defined due to compatibility condition (p1)).
Put $\cal R^*$ to be the good quotient of $\cal R$ given by
$$
{\cal R^*}:=(\{ X_i/{\cal Q}_i \}_{i\ge1}, \{ \pi^*_i \}_{i\ge1}).
$$

\break

\medskip\noindent
{\bf Step 3.} 

For the remaining part of the proof of Proposition 5.5,
fix some countable dense subset $D\subset X$ containing all essential
vertices of $X$. Recall that $V$ and $E$ denote the sets of all essential
vertices and essential edges of $X$, respectively.

In this step of the proof we investigate the quotient spaces $X_i/{\cal Q}_i$
appearing in the sequence $\cal R^*$ constructed in the previous step.
We start with the following observation concerning 
$X_1/{\cal Q}_1$. In its statement, we denote by 
$\rho:X_1\to X_1/{\cal Q}_1$ the quotient map.

\medskip\noindent
{\bf 5.11 Lemma.}
{\it There is a homeomorphism $h_1:X_1/{\cal Q}_1\to X$ satisfying the following
conditions:}
\item{(1)} {\it $h_1$ respects labels of vertices and edges, 
in the
following sense: if $v_1$ is a vertex of $X_1$ labelled with $v\in V$ then
$h_1\rho(v_1)=v$; similarly, if $x$ is a point contained in an edge $e_1$
of $X_1$ labelled with $e\in E$ then $h_1\rho(x)$ is contained in the edge
$e$ of $X$;}
\smallskip
\item{(2)} {\it if we denote by $D_1\subset X_1/{\cal Q}_1$
the subset consisting of all points corresponding to the shrinked subsets
$J\in{\cal J}_1$, then $h_1(D_1)=D$.}

\medskip\noindent
{\bf Proof:}
We view existence of a homeomorphism $h_1$ satsfying (1) as rather obvious.
It easily follows from the fact that, after taking the quotient of a segment
given by shrinking to points its subsegments from any pairwise disjoint
family, we get a space homeomorphic to a segment. Additional property (2)
can be achieved by the following observations. First, note that the union of the
family ${\cal J}_1$ is dense in $X_1$. This is a consequence of
observation (d6) and condition (j2). 
It follows that $D_1$ is a countable dense subset of
$X_1/{\cal Q}_1$. Moreover,by condition (i2) in Definition 5.4, 
$D_1$ contains all vertices of $X_1/{\cal Q}_1$.
The fact that $h_1$ can be chosen to satisfy additionally condition (2)
follows then from Lemma $\Delta.1$. 

\medskip
Recall that we denote by $\Lambda_n$ the blow-up locus of the bonding map
$\pi_n:X_{n+1}\to X_n$ from the inverse sequence $\cal R$.
It follows from the description of $X$-blow-up maps in Definitions 5.2 and 5.3
that we can view $X_n\setminus\Lambda_n$ as a subset of $X_{n+1}$.
Under this perspective, we also view $X_n\setminus\hbox{int}(J(\Lambda_n))$
as a subset of $X_{n+1}$. Furthermore, 
it follows from our previous considerations
that any set $J\in{\cal J}_{n+1}$ is contained either in 
$X_n\setminus J(\Lambda_n)$ (viewed as a subset of $X_{n+1}$) or in
$X_{n+1}\setminus[X_n\setminus\hbox{int}(J(\Lambda_n))]$.
Actually, the sets $J\subset X_n\setminus J(\Lambda_n)$ form a subfamily in
${\cal J}_{n+1}$ which coincides with the family ${\cal J}'_{n+1}$
in our earlier recursive description of the family ${\cal J}_{n+1}$ in Step 2, 
while the remaining subsets $J$
constitute the subfamily ${\cal J}''_{n+1}$.
In particular, there is a well defined quotient space
$(X_{n+1}\setminus(X_n\setminus J(\Lambda_n))) / {\cal J}''_{n+1}$,
and we denote by 
$$
\rho_{n+1}: X_{n+1}\setminus(X_n\setminus J(\Lambda_n)) \to
(X_{n+1}\setminus(X_n\setminus J(\Lambda_n))) / {\cal J}''_{n+1}   
$$
the quotient map.
Under the just recalled notations and conventions, and
by the arguments similar as in the proof of Lemma 5.11, we get
the following result, the proof of which we skip.

\medskip\noindent
{\bf 5.12 Lemma.}

\item{(0)}
{\it For each $n\ge1$ the quotient space $X_n/{\cal Q}_n$ is homeomorphic
to $X_n$, and it naturally inherits from $X_n$ the structure of an $X$-graph.}

\smallskip
\item{(1)}
{\it For each $n\ge1$, view $X_n\setminus\hbox{int}(J(\Lambda_n))$
as a subspace of $X_{n+1}$, and the family ${\cal J}'_{n+1}$ as the partition
of this subspace (which consists of the sets from the family ${\cal J}'_{n+1}$
and of the singletons). Denote also by $p$ the point in the quotient space
$X_n/{\cal Q}_n$ corresponding to the shrinked subset $J(\Lambda_n)$,
and let $(X_n/{\cal Q}_n)^\#(p)$ be the blow-up of $X_n/{\cal Q}_n$ at $p$,
as described in Section 2. Then the map 
$$
g_n:(X_n\setminus\hbox{int}(J(\Lambda_n))/{\cal J}'_{n+1}\to 
(X_n/{\cal Q}_n)^\#(p)
$$ 
induced by the identity map of $X_n$ is a homeomorphism.}

\smallskip
\item{(2)} {\it Suppose that $\Lambda_n$ is a blow-up segment. 
Denote by $e'$ this edge of $X_n$ which contains $\Lambda_n$,
and by $e\in E$ its label (in the $X$-graph structure of $X_n$). 
Let $p\in D\subset X$ be any point from the interior of the edge $e$ of $X$,
and let $X^{\#}(p)$ be the blow-up of $X$ at $p$, as described 
in Section 2. Then 
${\cal J}''_{n+1}$, viewed as a partition of the subspace
$X_{n+1}\setminus(X_n\setminus\ J(\Lambda_n))$
consisting of the subsets from ${\cal J}''_{n+1}$ and of singletons,
is upper semi-continuous, and
there is a homeomorphism
$$
h_{n+1}:
(X_{n+1}\setminus(X_n\setminus J(\Lambda_n))) / {\cal J}''_{n+1}
\to  X^\#(p)
$$
satisfying the following conditions:}
\smallskip
\itemitem{(a)} {\it $h_{n+1}$ respects the 
labels of vertices and edges in the
following sense: if $v'$ is a vertex of $X_{n+1}$ labelled with $v\in V$,
and contained in
$X_{n+1}\setminus(X_n\setminus J(\Lambda_n))$,
then $h_{n+1}\rho_{n+1}(v')=v$; similarly, if $x$ is a point contained 
in an edge $\varepsilon'$
of $X_{n+1}$ labelled with $\varepsilon\in E$, and such that
$x\in X_{n+1}\setminus(X_n\setminus J(\Lambda_n))$,
then $h_{n+1}\rho_{n+1}(x)$ is contained in the edge
$\varepsilon$ of $X^\#(p)$ (or in the blown-up edge $e^\#(p)$ 
if $\varepsilon=e$);}
\smallskip
\itemitem{(b)} {\it if we denote by 
$
D_{n+1}\subset 
(X_{n+1}\setminus(X_n\setminus J(\Lambda_n))) / {\cal J}''_{n+1}
$
the subset consisting of all points corresponding to the shrinked subsets
$J\in{\cal J}''_{n+1}$, then 
$h_{n+1}(D_{n+1})=D\setminus\{ p \}$, where $D\setminus\{ p \}$
is naturally understood as a subset in $X^\#(p)$.}
\smallskip

\item{(3)}
{\it Suppose $\Lambda_n$ is a blow-up point, i.e. $\Lambda_n=\{ v' \}$
for some vertex $v'$ of $X_n$. Let $v\in V$ be the label of $v'$
(in the $X$-graph structure of $X_n$), and let $X^\#(v)$ be the blow-up of $X$
at $v$, as described in Section 2. Then 
${\cal J}''_{n+1}$, viewed as a partition of the subspace
$X_{n+1}\setminus(X_n\setminus J(\Lambda_n))$
consisting of the subsets from ${\cal J}''_{n+1}$ and of singletons,
is upper semi-continuous, and
there is a homeomorphism
$$
h_{n+1}:
(X_{n+1}\setminus(X_n\setminus J(\Lambda_n))) / {\cal J}''_{n+1}
\to  X^\#(v)
$$
satisfying the following conditions:}
\smallskip
\itemitem{(a)} {\it $h_{n+1}$ respects the 
labels of vertices and edges in the
following sense: if $w'$ is a vertex of $X_{n+1}$ labelled with $w\in V$,
and contained in $X_{n+1}\setminus(X_n\setminus J(\Lambda_n))$,
then $h_{n+1}\rho_{n+1}(w')=w$; similarly, if $x$ is a point contained 
in an edge $\varepsilon'$
of $X_{n+1}$ labelled with $\varepsilon\in E$, and such that
$x\in X_{n+1}\setminus(X_n\setminus J(\Lambda_n))$,
then $h_{n+1}\rho_{n+1}(x)$ is contained in the edge
$\varepsilon$ of $X^\#(v)$ (under the obvious interpretation of edges of $X$
as edges of $X^\#(v)$);}
\smallskip
\itemitem{(b)} {\it if we denote by 
$
D_{n+1}\subset 
(X_{n+1}\setminus(X_n\setminus J(\Lambda_n))) / {\cal J}''_{n+1}
$
the subset consisting of all points corresponding to the shrinked subsets
$J\in{\cal J}''_{n+1}$, then 
$h_{n+1}(D_{n+1})=D\setminus\{ v \}$, where $D\setminus\{ v \}$
is naturally understood as a subset in $X^\#(v)$.}

\break

\medskip\noindent
{\bf Step 4.}

In this last step of the proof of Proposition 5.5 
we give a description of a violated
reflection inverse system ${\cal S}_X^v$, its cofinal subsequence $\cal P$,
and an isomorphism $\phi:{\cal R^*}\to{\cal P}$ of inverse sequences.
A subsequence $\cal P$ of ${\cal S}^v_X$ will have the following form
(under notation as in Section 4):
for some increasing sequence $H_n:n\ge1$ of finite subtrees of $T$ such that 
each $H_n$ has exactly $n$ vertices, and such that 
$\bigcup_{n\ge1}H_n=T$, 
it holds that
${\cal P}=(\{ Z_{H_n} \}_{n\ge1}, \{ \pi^v_{H_{n+1}H_n} \}_{n\ge1})$.
An isomorphism $\phi:{\cal R^*}\to{\cal P}$ will be described as a sequence
$\phi_n:X_n/{\cal Q}_n\to Z_{H_n}$ of homeomorphisms which commute
with the bonding maps of the sequences, i.e. satisfy the equations 
$\phi_n\pi_n^*=\pi^v_{H_{n+1}H_n}\phi_{n+1}$.

We pass to the description of the tree $T$ and the sequence $H_n$
of its subtrees. For each $n\ge1$, denote by $k_n$ the smallest number $k$
such that the blow-up locus $\Lambda_n$ of the map $\pi_n$
is a direct blow-up locus in $X_k$. Obviously, we have $1\le k_n\le n$.
For the vertex set $V_T$ of $T$ take a countable infinite set
$t_1,t_2,\dots$ ordered into a sequence. For the edge set $E_T$
take the sequence $e_n:n\ge1$, where each $e_n$ connects $t_{n+1}$
with $t_{k_n}$. For each $n\ge 1$, take as $H_n$ the subtree of $T$
spanned on the vertex set $\{ t_1,\dots,t_n \}$; 
note that then  $\{ t_1,\dots,t_n \}$
is the vertex set of $H_n$, and $\{ e_1,\dots,e_{n-1} \}$ is its edge set.

Recall that at the beginning of Step 3 of this proof we have fixed
a countable dense subset $D\subset X$ containig all vertices of $X$.
Recall also that, in order to describe a violated inverse sysem ${\cal S}^v_X$,
we need:
\item{(a)} to associate to any edge $e_n$ of the above described
tree $T$  a label $\lambda(e_n)$ from the set $D$, so that the label function
$\lambda:E_T\to D$ satisfies the requirements described in Section 4, and
\item{(b)} for each edge $e_n$ whose label $\lambda(e_n)$ is not
a vertex of $X$, to decide whether the corresponding gluing map 
$\beta_\epsilon:P_{\lambda(e_n)}\to P_{\lambda(e_n)}$
(and $\beta_{\bar\epsilon}$, which is the inverse of $\beta_\epsilon$),
where $|\epsilon|=e_n$, is the identity or the transposition.

\noindent
We will describe the data (a) and (b) as above recursively, together with the
description of the homeomorphisms $\phi_n$ constituting an isomorphism $\phi$.

Start with a choice of any homeomorphism $h_1$ as in Lemma 5.11.
Recall that $J(\Lambda_1)$ is this set from the family ${\cal J}_1$
which contains the blow-up locus $\Lambda_1$ of the map $\pi_1$.
Label the edge $e_1$ of $T$ with the element $h_1([J(\Lambda_1)])\in D_1$,
where $[J(\Lambda_1)]$ is the point in $X_1/{\cal Q}_1$ corresponding
to the shrinked set $J(\Lambda_1)$. Since, by the description of Section 4,
the space $Z_{H_1}$ is canonically equal to $X$, we can take as
$\phi_1:X_1/{\cal Q}_1\to Z_{H_1}$ simply the homeomorphism $h_1$.

Now, assume that for some $n$ we have already associated the labels 
$\lambda(e_i)$
to the edges $e_1,\dots,e_{n-1}$, and that we have already described the gluing
maps $\beta_\varepsilon$ related to the same edges $e_1,\dots,e_{n-1}$.
This portion of data determines uniquely the spaces $Z_{H_i}$ for all $i\le n$,
and the maps $\pi^v_{H_{i+1}H_i}$ for all $i\le n-1$.
Assume also that we have already defined the homeomorphisms
$\phi_i:X_i/{\cal Q}_i\to Z_{H_i}$ for all $i\le n$, and that they satisfy
the appropriate earlier mentioned commutativity equations.
Finally, suppose that for some appropriately chosen homeomorphisms
$h_1$ as in Lemma 5.11 and $h_i:2\le i\le n$ as 
in parts (2) or (3) of Lemma 5.12, the homeomorphisms $\phi_i$ and
the labels $\lambda(e_i)$ above satisfy the following:
\smallskip
\itemitem{(f1)} under the natural identification of $Z_{H_1}$ with $X$,
we have $\phi_1=h_1$;
\smallskip
\itemitem{(f2)} for all $1\le i\le n-1$, we have 
$\lambda(e_i)=h_{k_i}([J(\Lambda_{i})])$ and 
$h_{i+1}(D_{i+1})=D\setminus \{ \lambda(e_i) \}$;
\smallskip
\itemitem{(f3)} if $2\le i\le n$, then under the natural identification of
$X_{t_i}$ with $X$, consider the subset 
$X^\#_{t_i}(\lambda(e_{i-1}))\subset Z_{H_i}$; consider also the subset
$(X_i\setminus(X_{i-1}\setminus J(\Lambda_{i-1})))/ {\cal J}''_i 
\subset X_i/{\cal Q}_i$; then the image of the restriction of $\phi_i$ to
$(X_i\setminus(X_{i-1}\setminus J(\Lambda_{i-1})))/ {\cal J}''_i$
coincides with $X^\#_{t_i}(\lambda(e_{i-1}))$, and the restriction itself
coincides with $h_i$, i.e.
$$
\phi_i|_{(X_i\setminus(X_{i-1}\setminus J(\Lambda_{i-1})))/ {\cal J}''_i}=h_i.
$$

\noindent
Note that all the above assumptions are indeed satisfied in the already described
case of $n=1$, mostly trivially. 
To accomplish an inductive step of our construction,
we need to determine the label $\lambda(e_{n})$, 
the form of the gluing map $\beta_\epsilon$ for $|\epsilon|=e_n$,
and the shape of the homeomorphism $\phi_{n+1}$.

To determine the label $\lambda(e_n)$, 
observe that, by definition of the number $k_n$, the set $\Lambda_n$
is a direct blow-up locus in $X_{k_n}$. Consequently, the set $J(\Lambda_n)$
is also a subset of $X_{k_n}$. Moreover, by minimality of $k_n$ in the
definition of $k_n$, we have that $J(\Lambda_n)\in{\cal J}''_{k_n}$
(where we use the convention that ${\cal J}''_1:={\cal J}_1$).
We also have 
$J(\Lambda_n)\subset X_{k_n}\setminus (X_{k_n-1}\setminus J(\Lambda_{k_n-1}))$ if $k_n\ge2$.
We put $\lambda(e_n)=h_{k_n}([J(\Lambda_n)])$, where $[J(\Lambda_n)]$
denotes the point corresponding to the shrinked subset $J(\Lambda_n)$ 
in $X_1/{\cal Q}_1$ if $k_n=1$, and in
$(X_{k_n}\setminus (X_{k_n-1}\setminus J(\Lambda_{k_n-1}))) /
{\cal J}''_{k_n}$ if $k_n\ge2$.

Note that if $\Lambda_n$ is not a segment then $\lambda(e_n)$ is not
a vertex of $X$. In this case we need to determine the form of the gluing map
$\beta_\epsilon:P(\lambda(e_n))\to P(\lambda(e_n))$ for $|\epsilon|=e_n$.
Denote by $e'$ the edge of $X_{k_n}$ which contains $\Lambda_n$,
and let $e\in E$ be the label of $e'$ in the $X$-graph structure of $X_{k_n}$.
Denote by $u,v\in V$ the endpoints of $e$ (which are the vertices of $X$),
and note that the endpoints of $e'$ (which are the vertices of $X_{k_n}$)
are labeled with $u$ and $v$. Denote by $u',v'$ the endpoints of $e'$ 
labeled with $u,v$, respectively. Denote also by $a_u$ this endpoint
of $\Lambda_n$ which is closer to $u'$ in $e'$ than the other endpoint.
Referring to the description of $X_{n+1}$ as obtained from 
$X_n\setminus\hbox{int}(\Lambda_n)$ and $X\setminus\hbox{int}(e)$
by gluing the endpoints of $\Lambda_n$ to the endpoints of $e$
(as in Definition 5.3),
for $|\epsilon|=e_n$, put $\beta_\epsilon$ to be
\itemitem{$\bullet$} the identity of $P(\lambda(e_n))$, if $a_u$ is identified 
in $X_{n+1}$ with $u$;
\itemitem{$\bullet$} the transposition of $P(\lambda(e_n))$, if $a_u$
is identified with $v$.

\smallskip
To describe a homeomorphism 
$\phi_{n+1}:X_{n+1}/{\cal Q}_{n+1}\to Z_{H_{n+1}}$ as required,
choose any homeomorphism
$h_{n+1}:(X_{n+1}\setminus(X_n\setminus J(\Lambda_n) )) / {\cal J}''_{n+1}
\to X^\#(\lambda(e_n))$ as in part (2) or (3) of Lemma 5.12,
where $\lambda(e_n)$ plays the role of $p$ or $v$.
Viewing $(X_n\setminus J(\Lambda_n)) / {\cal J}'_n$ naturally as a subspace
in both $X_{n+1}/{\cal Q}_{n+1}$ and $X_n/{\cal Q}_n$,
and under the identification of $X$ with $X_{t_{n+1}}$
(and hence also of $X^\#(\lambda(e_n))$ with 
$X^\#_{t_{n+1}}(\lambda(e_n))$), we put
\itemitem{$\bullet$} $\phi_{n+1}(x)=\phi_n(x)$ if 
$x\in(X_n\setminus J(\Lambda_n)) / {\cal J}'_n $, and
\itemitem{$\bullet$} $\phi_{n+1}(x)=h_{n+1}(x)$ if 
$x\in (X_{n+1}\setminus(X_n\setminus J(\Lambda_n) )) / {\cal J}''_{n+1}$.

\noindent
We skip the straightforward verification that 
\item{$\bullet$}
$\phi_{n+1}:X_{n+1}/{\cal Q}_{n+1}\to Z_{H_{n+1}}$ as above
is a well defined homeomorphism, and that it satisfies the commutativity
equation $\phi_n\pi_n^*=\pi^v_{H_{n+1}H_n}\phi_{n+1}$, and
\item{$\bullet$} the data $\phi_i:1\le i\le n+1$, $\lambda(e_i):1\le i\le n$,
and $\beta_\epsilon$ for $|\epsilon|=e_i$ and $1\le i\le n$,
satisfy conditions (f1)-(f3), with $n+1$ substitutes for $n$.

\noindent
We also skip the straightforward verification of other details allowing
to conclude that,
by iterating a step of the inductive construction as above,
we obtain a well defined  violated reflection system ${\cal S}^v_X$,
its cofinal subsequence $\cal P$, and an isomorphism of inverse sequences
$\phi:{\cal R}^*\to{\cal P}$.

The proof of Propsition 5.5 is thus completed.


\magnification1200

\bigskip\noindent
{\bf 6. The Coxeter-Davis complexes 
of systems $(W_\Gamma,S_\Gamma)$.}

\smallskip
In this section we describe the Coxeter-Davis complexes $ \Sigma_\Gamma$
associated to the right-angled Coxeter systems 
$(W_\Gamma,S_\Gamma)$ whose
nerves are simplicial graphs $\Gamma$. 
We then focus on a class of graphs $\Gamma$ 
satisfying conditions ($*$1) and ($*$2)
below, and for this class we introduce some useful stratification
of complexes $ \Sigma_\Gamma$ into strata called {\it sectors}
and {\it branch components}. This stratification will be
essentially used in Section P, to analyze the geodesic inverse
system in the complex $ \Sigma_\Gamma$ for such graphs $\Gamma$ 
(the limit of this system is,
by definition, the visual boundary $\partial_\infty(W_\Gamma,S_\Gamma)$ 
of the corresponding Coxeter system).

\medskip\noindent
{\it The Coxeter-Davis complex $ \Sigma_\Gamma$.}

For this part of exposition, let $\Gamma$ be any finite simplicial 
graph which is flag as simplicial complex (i.e. does not contain
cycles of length 3). Denote by $V_\Gamma,E_\Gamma$ the sets of { all}
vertices and edges of $\Gamma$, respectively. 
Let $(W_\Gamma,S_\Gamma)$ be the right angled Coxeter system with
nerve $\Gamma$. 
This means that $W_\Gamma$ is a group and $S_\Gamma$ is a generating set
for $W_\Gamma$
such that
\item{$\bullet$} the generators in $S_\Gamma$ are identified with
the vertices of $\Gamma$;
\item{$\bullet$} $W_\Gamma$ has the presentation
$
\langle  S_\Gamma\, |\, \{ s^2:s\in S_\Gamma \}\cup
\{ (ss')^2:[s,s'] \hbox{ is an edge of }\Gamma \} \rangle.
$

\noindent
The next fact provides description and basic 
properties of an object called {\it the Coxeter-Davis complex} for 
the system $(W_\Gamma,S_\Gamma)$.
We refer the reader to [D] for more detailed introduction
to Coxeter systems and their Coxeter-Davis complexes.
It is worth keeping in mind that for right angled Coxeter systems 
$(W_\Gamma, S_\Gamma)$ under our consideration 
the Coxeter-Davis complex of any such system coincides with
the Cayley complex 
of the group $W_\Gamma$, for the canonical presentation mentioned above.
(This is explained with more details in remark right after Fact 6.1.)

Recall that a {\it square complex} is a 2-dimensional cubical
complex.

\medskip\noindent
{\bf 6.1 Definition/Fact.} {\it Suppose that $\Gamma$ is a simplicial graph
which is flag, and let $\Sigma_\Gamma$ denote  
the Coxeter-Davis complex of the system $(W_\Gamma,S_\Gamma)$.}

\noindent
{\it (1) $\Sigma_\Gamma$
is the unique labelled square complex satisfying
the following properties:}

\itemitem{(cd1)} {\it squares of $ \Sigma_\Gamma$ are labelled
with the elements of $E_\Gamma$, and edges (1-cells) 
with elements of $V_\Gamma$;} 

\itemitem{(cd2)} {\it for any vertex $x$ of $ \Sigma_\Gamma$ the link
of $ \Sigma_\Gamma$ at $x$, denoted $( \Sigma_\Gamma)_x$, which we view as
a graph equipped with the induced labelling of vertices and edges 
(with elements of $V_\Gamma$ and $E_\Gamma$, respectively), is isomorphic,
as labelled graph, to $\Gamma$ equippped with the tautological
labellings of vertices and edges;} 

\itemitem{(cd3)} {\it $ \Sigma_\Gamma$ is simply connected.} 

\smallskip\noindent
{\it (2) $ \Sigma_\Gamma$ is a CAT(0) cubical complex. More precisely,
by this we mean that $\Sigma_\Gamma$ comes equipped with the so called 
standard piecewise euclidean metric (for which each cubical cell is
isometric with the unit euclidean cube of the corresponding dimension).
With this metric, $\Sigma_\Gamma$ is a geodesic metric space
satisfying the CAT(0) condition. We refer the reader to [BH] for an exhausting
exposition of the CAT(0) concept, and below in this section we recall
few of its aspects that will be used in this paper.}

\smallskip\noindent
{\it (3) The group $W_\Gamma$ acts on $ \Sigma_\Gamma$ by
combinatorial automorphisms preserving the labels. The action is
proper discontinuous and cocompact. In fact, $W_\Gamma$
coincides with the group of all label preserving combinatorial
automorphisms of $ \Sigma_\Gamma$, and its action is simply transitive
on the vertices.}

\medskip\noindent
{\bf Remark.}
The above description of $\Sigma_\Gamma$, as a labelled square complex,
can be given alternatively as follows. Let $C_\Gamma$ be the 
Cayley graph of the group $W_\Gamma$ with respect to the generating set
$S_\Gamma$. By the convention which we use in this paper,
$C_\Gamma$ has the vertex set which coincides with $W_\Gamma$,
and it has one unoriented edge for any pair $\{ g,gs \}$ of elements 
that differ by a multiplication on the right by a generator.
We label such an edge with the corresponding generator $s$.
For any edge $[s,s']\in E_\Gamma$, consider the special subgroup
$W_{\{ s,s' \}}<W_\Gamma$, and note that it is canonically isomorphic
to the dihedral group of order 4. To each left coset $gW_{\{ s,s' \}}$
in $W_\Gamma$ there is associated a subgraph in $C_\Gamma$
which is a cycle of length 4. For each such coset, we attach to $C_\Gamma$ a square 2-cell along the corresponding 4-cycle. We label this 2-cell with the
corresponding edge $[s,s']$. Doing this for all edges of 
$\Gamma$, we obtain a square 2-complex that we denote by $C^2_\Gamma$
and call the {\it Cayley complex} of $(W_\Gamma,S_\Gamma)$.
It is a well known fact that (under our assumptions on $\Gamma$), 
as a labelled complex, $C^2_\Gamma$
coincides with the Coxeter-Davis complex $\Sigma_\Gamma$,
as described in above Fact 6.1.

\medskip
Recall that {\it the visual boundary} of $(W_\Gamma,S_\Gamma)$ is, 
by definition, the visual boundary of the Coxeter-Davis 
complex $ \Sigma_\Gamma$. We donote it 
$\partial_\infty(W_\Gamma,S_\Gamma)$. 
If $W_\Gamma$ happens to be word-hyperbolic
(as it is the case  for all $\Gamma$ satisfying conditions ($*$1) and ($*$2)
below),
$\partial_\infty(W_\Gamma,S_\Gamma)$ coincides with the 
Gromov boundary of $W_\Gamma$ (see Proposition III.3.7(2) in [BH]).

\medskip\noindent
{\it The stratification of $ \Sigma_\Gamma$: sectors and branch components.}

Recall that for a simplicial graph $\Gamma$ we denote by $|\Gamma|$
the underlying topological graph.
In the remaining part of the section we restrict our attention to  the
class of graphs $\Gamma$ that are connected, not reduced
to a single vertex, and
satisfy the following properties:

\itemitem{($*$1)} the underlying topological graph $|\Gamma|$ 
contains no essential loop edges;

\itemitem{($*$2)} each essential edge of $|\Gamma|$ consists of at least 3 edges 
of the corresponding simplicial graph $\Gamma$.


\medskip\noindent
{\bf 6.2 Definition.}
Given an essential edge $\varepsilon$ of $|\Gamma|$
(which, by our assumptions ($*$1) and ($*$2), is not a loop
and consists of at least 3 edges of $\Gamma$), a {\it sector
of $ \Sigma_\Gamma$  labelled with $\varepsilon$} (shortly, an
{\it $\varepsilon$-sector}) is any maximal connected subcomplex
of $ \Sigma_\Gamma$ which is the union of squares labelled with
edges $e\in E_\Gamma$ contained in $\varepsilon$.

\medskip\noindent
{\bf Remark.}
Alternatively, $\varepsilon$-sectors can be described as follows.
Consider the special subgroup $W_\varepsilon<W_\Gamma$
with its standard generating set $S_\varepsilon=S_\Gamma\cap\varepsilon$.
Viewing $\Sigma_\Gamma$ as the Cayley complex of $(W_\Gamma,S_\Gamma)$,
we have a natural identification of the Cayley complex $C^2_\varepsilon$ of 
$(W_\varepsilon,S_\varepsilon)$ with a subcomplex in $\Sigma_\Gamma$.
The $\varepsilon$-sectors are then the translates in $\Sigma_\Gamma$,
under the action of $W_\Gamma$,
of the subcomplex $C^2_\varepsilon$. 

\medskip
The next fact describes the combinatorial types of sectors of 
$ \Sigma_\Gamma$, as square complexes.
We omit its straightforward proof.

\medskip\noindent
{\bf 6.3 Fact.}
{\it Let $k\ge3$ be the number of edges of $E_\Gamma$ contained in an
essential (non-loop) edge $\varepsilon$ of $|\Gamma|$. Then each 
$\varepsilon$-sector of $ \Sigma_\Gamma$ is isomorphic to the square
complex $\Omega_k$ uniquely determined by the following
properties:}

\itemitem{(k1)} {\it each vertex link of $\Omega_k$ is a polygonal
arc consisting of $k$ edges;} 

\itemitem{(k2)} {\it $\Omega_k$ is simply connected;}

\itemitem{(k3)} {\it for each square $Q$ in $\Omega_k$ each of
the two reflections in $Q$ through a midsegment of $Q$ parallel
to a pair of its opposite sides extends to a combinatorial automorphism
of the whole $\Omega_k$.}

\noindent
{\it It is a consequence of conditions (k1) and (k2) that
$\Omega_k$ is topologically a noncompact planar surface
with infinitely many boundary components which are all
noncompact. }

\medskip\noindent
{\bf Remark.} It is also not hard to realize that each $\varepsilon$-sector
$\Omega$ in $ \Sigma_\Gamma$
is combinatorially isomorphic, as labelled square complex
with labels inherited from $\Sigma_\Gamma$,
to the Coxeter-Davis complex $ \Sigma_\varepsilon$ of the
right angled Coxeter system $(W_\varepsilon,S_\varepsilon)$ 
corresponding
to the parabolic subgroup of $W_\Gamma$ spanned
by the standard generators corresponding to the vertices in
$\varepsilon$. 
The boundary $\partial\Omega$ is the subcomplex
of $\Omega$ consisting of all edges which are labelled with the
two essential vertices of $|\Gamma|$ contained in $\varepsilon$
(i.e. the endpoints of $\varepsilon$). Each component
$L$ of $\partial\Omega$ is an infinite polygonal line,
and the edges in $L$
are labelled alternately with the two endpoint vertices
of $\varepsilon$.

\medskip
Some useful properties of sectors, as subcomplexes of $ \Sigma_\Gamma$,
are gathered in the next fact. 
Recall that a subset $A$ of a geodesic metric space $X$ is {\it convex}
if for any two points of $A$, any geodesic in $X$ connecting these points
is contained in $A$. If $A$ is a subcomplex of a CAT(0) cubical complex $X$,
then convexity of $A$ has the following characterization: $A$ is convex iff
it is connected and for any vertex $v$ of $A$ the link of $A$ at $v$
is a full subcomplex in the link of $X$ at $v$. In the statement below we 
use the term {\it strict
convexity} for the following property, which is clearly stronger than convexity: 
a subcomplex $A$
of a CAT(0) square complex $X$ is {\it strictly convex}
if it is connected and for any vertex $v$ of $A$ the link
$A_v$ is {\it 3-convex} in the link $ \Sigma_v$. The latter means that
$A_v$ is a full subgraph of $ \Sigma_v$, and that any polygonal path
of $ \Sigma_v$ intersecting $A_v$ only at its endpoints has length
at least 3. 

\medskip\noindent
{\bf 6.4 Fact.}

\item{(1)}
{\it Each sector $\Omega$, as well as its any boundary 
component $L$, is a strictly convex subcomplex
of $ \Sigma_\Gamma$.}

\item{(2)}
{\it Each connected component of the complement
$ \Sigma_\Gamma\setminus\hbox{\rm int}(\Omega)$ is a strictly convex
subcomplex of $ \Sigma_\Gamma$.}

\item{(3)}
{\it The inclusion provides a bijective correspondence between
the boundary components of $\Omega$ and the connected components
of $ \Sigma_\Gamma\setminus\hbox{\rm int}(\Omega)$.}

\medskip\noindent 
{\bf Proof:}
Suppose that $\Omega$ is an $\varepsilon$-sector.

To prove part (1), note that for any vertex $v$ of $\Omega$
the pair of links $((\Sigma_\Gamma)_v, \Omega_v)$ is isomorphic
to the pair $(\Gamma,\varepsilon)$. Strict convexity of $\Omega$
follows then from condition ($*$2). Similarly, if $L$ is a boundary
component of $\Omega$, then for any vertex $v$ of $L$ the pair
$((\Sigma_\Gamma)_v, L_v)$ is isomorphic to the pair $(\Gamma,\{ a,b \})$,
where $a,b$ are the endpoints of $\varepsilon$.
Strict convexity of $L$ follows then again from condition ($*$2).

To prove part (2), consider a connected component $Q$ of 
$\Sigma_\Gamma\setminus\hbox{int}(\Omega)$, and note that it is
a subcomplex of $\Sigma_\Gamma$. Moreover, if $v$ is a vertex of $Q$
not contained in $\Omega$, we have $Q_v=(\Sigma_\Gamma)_v$,
and if $v$ is contained in $\Omega$ then the pair of links
$((\Sigma_\Gamma)_v,Q_v)$ is isomorphic to the pair
$(\Gamma,\Gamma\setminus\hbox{int}(\varepsilon))$.
Again, strict convexity of $Q$ follows easily from condition ($*$2).

To get part (3), we need to show two things about any
connected component $Q$
of $\Sigma_\Gamma\setminus\hbox{int}(\Omega)$: first, that $Q$ intersects $\Omega$,
and second, that $Q$ does not intersect more than one boundary
component of $\Omega$.
To see the first assertion above, note that if
it were not true, $Q$ would be a connected
component of $\Sigma_\Gamma$, which contradicts connectedness
of the latter. The second assertion follows immediately from
convexity of $Q$ and of $\Omega$.

This completes the proof of Fact 6.4.

\medskip
We pass to the discussion of objects that we call {\it branch components} 
of $\Sigma_\Gamma$.
Denote by $V$ the set of essential vertices of the underlying
topological graph $|\Gamma|$.

\medskip\noindent
{\bf 6.5 Definition.} A {\it branch component} of $ \Sigma_\Gamma$
is any maximal connected subcomplex of $ \Sigma_\Gamma$
which is the union of edges of $ \Sigma_\Gamma$ labelled
with the essential vertices of $|\Gamma|$ (i.e. with 
the elements of $V$). A {\it branch locus} of $ \Sigma_\Gamma$
is the union of all branch components, and we denote it
by $C_\Gamma$.

\medskip\noindent
{\bf Remark.} An alternative description of branch components can be given as follows. Consider the special subgroup $W_V<W_\Gamma$ spanned on
those generators from $S_\Gamma$ that correspond to all essential vertices
of $|\Gamma|$. View the Cayley complex $C^2_V$ (which actually coincides
with the Cayley graph $C_V$) naturally as the subcomplex of the Cayley complex
$C^2_\Gamma$. After identifying the latter with $\Sigma_\Gamma$,
the branch components of $\Sigma_\Gamma$ are the translates
under the action of $W_\Gamma$
of the subcomplex $C^2_V$.

\medskip\noindent
{\bf 6.6 Fact.}
\item{(1)}
{\it Each branch component $C$ of $ \Sigma_\Gamma$ is isomorphic
to the regular tree of degree equal to the cardinality of the set $V$.
Moreover, for each vertex $x\in C$ the labelling restricted to the
edges of $C$ issuing from $x$ is a bijection on $V$.}

\item{(2)}
{\it Every edge of a branch component $C$ is contained either
in exactly one or in at least three squares of $ \Sigma_\Gamma$.}

\item{(3)}
{\it Each boundary component of any sector in $ \Sigma_\Gamma$
is contained in some branch component. }

\item{(4)}
{\it Each branch component $C$ is a strictly convex subcomplex
of $ \Sigma_\Gamma$.}

\item{(5)}
{\it The connected componets of the complement 
$\Sigma_\Gamma\setminus C$ of any branch component $C$ are in the natural
bijective correspondence (via inclusion) with (the interiors of) the sectors
of $\Sigma_\Gamma$ adjecent to $C$.}

\medskip\noindent
{\bf Proof:} Note that for each vertex $v$ of any branch component $C$
the vertex link $C_v$, as a subcomplex in the vertex link 
$(\Sigma_\Gamma)_v$, corresponds to the set $V$ viewed as 
a subcomplex of $\Gamma$. In view of the condition ($*$2),
it follows that $C$ is a strictly convex subcomplex of $\Sigma_\Gamma$.
In particular, it is a convex subcomplex, and therefore it must be
simply connected, and hence a tree. The assertions (2) and (3) follow
even more directly from the above observation concerning links.

To prove (5) note that, since $\Sigma_\Gamma$ is connected, the clusure
of each connected component $U$ of $\Sigma_\Gamma\setminus C$
intersects $C$. It follows that $U$ contains the interior of at least one
sector $\Omega$ adjacent to $C$. To see that there is only one such sector
$\Omega$, consider the metric completion $\overline U$ of $U$, which can be viewed as obtained by attaching to $U$, disjointly (!), the boundary components
$L^\Omega$ contained in $C$ of all sectors $\Omega$ as above.
The natural map $\bar\iota:\overline U\to\Sigma_\Gamma$ induced by
the inclusion $\iota:U\to\Sigma_\Omega$ is then easily seen to be the local
isometry, since all the vertex links of $\overline U$ embed onto full subcomplexes 
in the corresponding vertex links of $\Sigma_\Gamma$.
Suppose that we have at least two sectors $\Omega$ as above, say 
$\Omega_1$ and $\Omega_2$. Consider any points 
$p_i\in L^{\Omega_i}\subset\overline U$, for $i=1,2$,
and let $\Gamma$ be a geodesic in $\overline U$ connecting $p_1$ with $p_2$.
Since the subcomplexes $L^{\Omega_i}$ of $\overline U$ are disjoint,
$\gamma$ necessarily passes through $U$. 
On the other hand, since $\bar\iota$ is a local isometry, the image 
$\bar\iota(\gamma)$ is a geodesic in $\Sigma_\Gamma$, and it connects 
the points $\bar\iota(p_i)$, which both belong to $C$.
Since $C$ is convex in $\Sigma_\Gamma$, the geodesic $\bar\iota(\gamma)$
is contained in $C$, which contradicts the earlier observation that 
$\Gamma$ passes through $U$. This completes the proof.

\medskip
The separation properties of the sectors and the branch components
in $ \Sigma_\Gamma$ are nicely represented by a dual object,
the {\it adjacency graph} ${\cal A}_\Gamma$, which we now
introduce.

\medskip\noindent
{\bf 6.7 Definition.}
\item{(1)} A sector $\Omega$ is {\it adjacent} to a branch component
$C$ of $ \Sigma_\Gamma$ if some boundary component of $\Omega$
is contained in $C$.
\item{(2)} An {\it adjacency graph} ${\cal A}_\Gamma$
is the bipartite graph whose vertices represent all sectors and all branch 
components of $ \Sigma_\Gamma$, and whose edges correspond to
the adjacency relation defined in (1).

\medskip\noindent
{\bf 6.8 Fact.} {\it For any $\Gamma$ as above, the adjacency
graph ${\cal A}_\Gamma$ is a tree.}
 
\medskip\noindent
{\bf Proof:} If $|\Gamma|$ is a single essential edge
then $\Sigma_\Gamma$ reduces to a single sector, and in this case
the assertion is obvious. Otherwise, 
it is still a fairly straightforward observation that the graph 
${\cal A}_\Gamma$ is connected. Thus,
it is sufficient to show that each vertex of ${\cal A}_\Gamma$
separates this graph. For vertices represented by sectors of 
$\Sigma_\Gamma$ this follows easily from Fact 6.4(3),
while for vertices represented by branch components - from Fact 6.6(5).

\medskip\noindent
{\bf Remark.} In view of Fact 6.8, we will call ${\cal A}_\Gamma$
the {\it adjacency tree} of $ \Sigma_\Gamma$. Due to Fact 6.4(3),
this tree has infinite
order at every vertex corresponding to a sector. 
It also has infinite order at every vertex corresponding to a branch
component,
except the case when $|\Gamma|$
is a single essential edge.

\medskip
An argument similar as in the proof of Fact 6.8 yields also the following.

\medskip\noindent
{\bf 6.9 Fact.}
{\it Suppose that each of the three distinct objects $S,S_1,S_2$ is
either the interior of some sector or a branch
component of $ \Sigma_\Gamma$, and let $s,s_1,s_2$ be the vertices
in the adjacency tree ${\cal A}_\Gamma$ corresponding to
$S,S_1,S_2$, respectively. Then $S$ separates $S_1$ from $S_2$
in $ \Sigma_\Gamma$ if and only if $s$ separates $s_1$ from $s_2$
in ${\cal A}_\Gamma$.}


\magnification1200

\noindent
{\bf 7. Approximation Lemma and the proof of Main Theorem.}

\smallskip
In this section we formulate a lemma 
(Lemma 7.11 below), which we call Approximation Lemma, 
and which, together with Proposition 5.5, 
immediately verifies Theorem 1.1 of the introduction for
topological graphs $X=|\Gamma|$ that are connected, have 
no esential loop edges and no vertices of degree 1.
We also state the main result of this paper in its full generality, as Theorem 7.12 below
(which slightly extends Theorem 1.1 of the introduction). 
We show how Theorem 7.12 follows from Approximation Lemma, in view of the results
of Sections 3 and 5. The proof of Approximation Lemma is postponed until
the next section (Section 8).

\medskip
We start with describing the setting, and with
some terminological preparations. 
Let $\Gamma$ be a finite connected
simplicial graph which is flag, and 
which satisfies the assumptions ($*$1) and ($*$2) from the previous section 
(i.e. the underlying topological graph $|\Gamma|$ has no essential loop edge
and each essential edge of $|\Gamma|$ consists of at least 3 edges of $\Gamma$),
and assume additionally that $\Gamma$ has no vertices of degree 1. 
Let $(W_\Gamma,S_\Gamma)$ be the right angled Coxeter
system with nerve $\Gamma$, and let  $\Sigma_\Gamma$ be the associated 
Coxeter-Davis complex, equipped with the standard piecewise euclidean
metric. Recall that, with this metric $\Sigma_\Gamma$ is a $CAT(0)$
geodesic metric space (in particular, any two points of $\Sigma_\Gamma$
are connected by a unique geodesic).
Let $x_0$ be any vertex of  $\Sigma_\Gamma$. Denote by $S_r$ the sphere
of radius $r$ in $\Sigma_\Gamma$ centered at $x_0$, and view it as
a topological space with the induced 
topology. For $r<r'$ let $g_{r'r}:S_{r'}\to S_r$ be the
{\it geodesic projection}, i.e. the map which to each point
$x'\in S_{r'}$ associates the unique point $x\in S_r$
lying on the geodesic $[x_0,x']$. The {\it geodesic inverse system}
in $(\Sigma_\Gamma,x_0)$ is the system $(\{ S_r \},\{ g_{r'r} \})$
of spheres and geodesic projections. By definition, we have 
$$
\partial_{\infty}(W_\Gamma,S_\Gamma):=\partial_{\infty}\Sigma_\Gamma=
\lim_{\longleftarrow}(\{ S_r \},\{ g_{r'r} \})
$$
(see II.8.5 in [BH], or the comments after Definition (2b.1) in [DJ]).

\medskip\noindent
{\bf 7.1 Definition.}
A {\it strict star domain} in $(\Sigma_\Gamma,x_0)$ is any 
closed bounded subset $D$ of $\Sigma_\Gamma$ 
distinct from the singleton $\{ x_0 \}$ 
such that for each $x\in D\setminus\{ x_0 \}$ 
the geodesic segment $[x_0,x]$ is contained in $D$, and its
part $[x_0,x)=[x_0,x]\setminus\{x\}$ is contained
in the interior $\hbox{int}(D)$ of $D$.

\medskip
Note that balls centered at $x_0$ are strict star domains
in $(\Sigma_\Gamma,x_0)$. We mention without proof the following 
obvious fact.

\medskip\noindent
{\bf 7.2 Fact.} {\it Let $D$ be a strict star domain in $(\Sigma_\Gamma,x_0)$.
Then for any point $y\in \Sigma_\Gamma\setminus\hbox{\rm int}(D)$ the geodesic
segment $[x_0,y]$ intersects the boundary $\partial D$ in
precisely one point.}

\medskip\noindent
{\bf 7.3 Definition.} Given two strict star domains $D\subset D'$
in $(\Sigma_\Gamma,x_0)$ the {\it geodesic projection}
$g:\partial D'\to\partial D$ is the map which to any point
$x'\in\partial D'$ associates the unique point $x\in\partial D$
lying on the geodesic segment $[x_0,x']$.

\medskip
The above notion of geodesic projection generalizes the one
for concentric spheres. It also has the following property.

\medskip\noindent
{\bf 7.4 Lemma}
{\it For any strict star domains $D\subset D'$
in $(\Sigma_\Gamma,x_0)$ the {geodesic projection}
$g:\partial D'\to\partial D$ is continuous with respect to the topologies
induced from $\Sigma_\Gamma$.}

\medskip\noindent
{\bf Proof:}
Fix any $x'\in\partial D'$, and any sequence $(x_n')$ of points of
$\partial D'$ converging to $x'$. Put $x=g(x')$, $x_n=g(x_n')$, and suppose
that $(x_n)$ does not converge to $x$. Since $\partial D$ is compact,
there is a subsequence of $(x_n)$ which converges to some $y\in\partial D$
distinct from $x$. Since geodesics in CAT(0) spaces depend continuously
on their endpoints (see Proposition 2.2 on p. 176 in [BH]),
the point $y$ must belong to the geodesic $[x_0,x']$. Since $x$ also
belongs to this geodesic, and since both $x$ and $y$ are contained
in $\partial D$, we get a contradiction with the definition of a strict star
domain, and hence $(x_n)$ converges to $x$.
This finishes the proof.

\medskip\noindent
{\bf 7.5 Definition.}
An {\it exhausting sequence} in $(\Sigma_\Gamma,x_0)$
is a sequence $D_n:n\ge1$ of strict star domains 
in $(\Sigma_\Gamma,x_0)$ such that

\itemitem{(es1)} for each $n\ge1$ we have $D_n\subset D_{n+1}$;
\itemitem{(es2)} for each $r>0$ there is $n$ such that the ball
of radius $r$ in $\Sigma_\Gamma$ centered at $x_0$ is contained
in $D_n$ (equivalently, the union 
$\bigcup_{n=1}^\infty D_n=\Sigma_\Gamma$). 

\medskip
To any exhausting sequence $(D_n)$ in $(\Sigma_\Gamma,x_0)$
there is associated an inverse sequence $(\{\partial D_n\},\{g_n\})$,
where $g_n:\partial D_{n+1}\to \partial D_{n}$ is the geodesic
projection.
We mention without proving the following easy observation.

\medskip\noindent
{\bf 7.6 Fact.}
{\it For any exhausting sequence $(D_n)$ in $(\Sigma_\Gamma,x_0)$
and the associated sequence of geodesic projections 
$g_n:\partial D_{n+1}\to\partial D_{n}$ we have}
$$
\lim_{\longleftarrow}(\{ \partial D_n \},\{ g_n \})=
\lim_{\longleftarrow}(\{S_r\},\{g_{r'r}\})=
\partial_{\infty}(W_\Gamma,S_\Gamma).
$$

\medskip\noindent
{\bf 7.7 Definition.}
A strict star domain $D$ in $(\Sigma_\Gamma,x_0)$ is {\it regular}
if the following conditions are satisfied

\itemitem{(r1)} the boundary $\partial D$ intersects the branch locus 
of $\Sigma_\Gamma$ at finitely many points;

\itemitem{(r2)} if $x\in\partial D$ belongs to the interior of an edge $e$
of $\Sigma_\Gamma$ contained in the branch locus then $\partial D$
is transversal to $e$ at $x$;

\itemitem{(r3)} if $x\in\partial D$ is a vertex of $\Sigma_\Gamma$,
then for some (small) open neighbourhood $U$ of $x$ in $\Sigma_\Gamma$
the intersection with $U$ of the interior of precisely one edge of the
branch locus that issues from $x$ is contained in  $\hbox{int}(D)$, and
intersections with $U$ of the interiors of the remaining
such edges are contained in $\Sigma_\Gamma\setminus D$. 

\medskip
Observe that the boundary $\partial D$ of any regular strict star
domain is naturally a finite topological graph whose vertex set coincides
with the intersection of $\partial D$ with the branch locus, 
and whose every essential edge coincides with a component of the intersection of $\partial D$
with some sector of $\Sigma_\Gamma$.
Moreover, we equip $\partial D$ (viewed as a topological graph)
with the labelling of its essential vertices and edges
described as follows:

\item{$\bullet$} if $u$ is a vertex of $\partial D$ which  
lies in the interior of some branch edge $e$ of $\Sigma_\Gamma$, 
we associate to $u$ the label $v\in V$
equal to the label of $e$;

\item{$\bullet$} if $u$ is a vertex of $\partial D$ which is
a vertex of $\Sigma_\Gamma$, we associate
to it the label $v$ of this unique branch edge issuing from $u$
whose interior intersects $\hbox{int}(D)$ arbitrarily close to $u$
(existence and uniqueness of such an edge follows from condition
(r3) in Definition 7.7);

\item{$\bullet$} each edge $d$ of $\partial D$ is contained in
a unique sector $\Omega$ of $\Sigma_\Gamma$, and we  
label $d$ with this $\varepsilon\in E$ for which
$\Omega$ is an $\varepsilon$-sector.

\smallskip
\noindent
We make a record of the following easy observation.

\medskip\noindent
{\bf 7.8 Lemma.}
{\it Let $D$ be a regular strict star domain in $(\Sigma_\Gamma,x_0)$, and view 
$\partial D$ as a topological graph equipped with the labelling of its essential vertices and edges, as described above. Then $\partial D$ is an $X$-graph.}

\medskip


Given any two continuous maps $f,g:Y\to Z$ between compact metric spaces, denote by
$|| f-g ||$ their uniform distance, i.e. the number
$|| f-g ||:=\sup_{y\in Y}d_Z(f(y),g(y))$.
In order to formulate Approximation Lemma 
(the main technical result of this section), we need to recall a result of
M. Brown (see Theorem 2 in [Br]).

\medskip\noindent
{\bf 7.9 Lemma (Brown's Lemma).}
{\it Let ${\cal S}=(\{ X_n \},\{ h_n \})$ be an inverse sequence of compact 
metric spaces. Then for any $k\ge0$ and for any sequence $s_k=(f_1,\dots,f_k)$
of continuous maps $f_i:X_{i+1}\to X_i$ (where by $s_0$ we mean the
empty sequence) there is a positive real number 
$\epsilon(s_k)=\epsilon(f_1,\dots,f_k)>0$ such that, if $(h_i')_{i\ge1}$
is a sequence of continuous maps $h_i':X_{i+1}\to X_i$ satisfying
the estimates
$$
|| h_{i+1}'-h_{i+1} ||<\epsilon(h_1',\dots,h_i') \quad 
\hbox{for all $i\ge0$}
\leqno{(B)}
$$
then for the modified inverse sequence 
${\cal S}'=(\{ X_n \},\{ h_n' \})$ we have
$$
\lim_{\longleftarrow}{\cal S}'=\lim_{\longleftarrow}{\cal S}.
$$
}

We will call any inverse sequence ${\cal S}'$ obtained from a sequence $\cal S$
subject to the condition $(B)$ as in the above lemma
a {\it Brown's approximation} of $\cal S$.



\medskip\noindent
{\bf 7.10 Remarks.}
\item{(1)} It is not hard to
observe that
a Brown's approximation can be always constructed recursively so that it 
additionally satisfies the following condition:

\item{} {\it for any $n\ge2$, let $\delta_n>0$ be a uniform continuity constant for $\epsilon_n={1\over n}$, common for all maps 
$g_i'\circ g_{i+1}'\circ\dots\circ g_{n-1}'$ with $i<n$; then for each $m\ge n$
we have
$$
||  g_n'\circ g_{n+1}'\circ\dots\circ g_{m}' - 
g_n\circ g_{n+1}\circ\dots\circ g_m  ||<\delta_n.
\leqno{(A)}
$$
}
\item{} We will demand and use the above condition (A) in our later arguments
in Section 8.

\item{(2)}
The proof of Brown's Lemma in [Br] gives also the following 
property, to which we will refer in our later arguments in Section 8.
{\it For any
approximation ${\cal S}'$ satisfying condition (B),
if $x\in\lim_{\leftarrow}{\cal S}'$ is represented by a thread $(x_n')$
then for each $k$ the limit
$$
y_k:=\lim_{n\to\infty}g_k\dots\circ g_{n-1}(x_n')
$$
exists, and the sequence $(y_k)$ is a thread of the inverse sequence $\cal S$.
Moreover, if we denote by $y\in\lim_\leftarrow{\cal S}$ the element
induced by $(y_k)$, then the map given by $x\mapsto y$ describes
a homeomorphism between $\lim_\leftarrow{\cal S}'$
and $\lim_\leftarrow{\cal S}$.}

\medskip\noindent
{\bf 7.11 Lemma (Approximation Lemma).}
{\it  Let $\Gamma$ be a finite connected simplicial graph which is flag,
not reduced to a single vertex,
which satisfies the assumptions ($*$1) and ($*$2) from Section 6
(i.e. $|\Gamma|$ has no essential loop edge, and each essential edge
of $|\Gamma|$ consists of at least 3 edges of $\Gamma$),
and which has no vertex of degree 1.
Then
there exists an exhausting 
sequence $D_n:n\ge1$ of regular strict star domains
in $(\Sigma_\Gamma,x_0)$ satisfying the following property:
the inverse sequence $( \{ \partial D_n \},\{ g_n \} )$
(where $g_n:\partial D_n\to\partial D_{n-1}$ are the geodesic projections,
and where we view each $\partial D_n$ as a $|\Gamma|$-graph)
admits a Brown's approximation $(\{ \partial D_n \},\{ g_n' \})$ 
which is isomorphic to some null and dense inverse sequence $\cal R$
of $|\Gamma|$-graphs and $|\Gamma|$-blow-ups.}

\medskip
We postpone the proof of the above lemma until the next section,
and we present now how this lemma implies the main result
of the paper stated below in its full generality, as Theorem 7.12.

Recall that a topological graph $X$ is {\it non-separable} if it is connected
and has no
separating essential vertex. (In particular, the circle is non-separable.)
A {\it block}
of $X$
is any subgraph of $X$ (for the natural cell structure consisting 
of essential vertices and edges)
which is maximal for the inclusion in the family of all non-separable
subgraphs of $X$, and which is not (homeomorphic to) a segment
or a singleton. We view each block as a topological
graph, and we note that its natural cell structure (consisting of its own essential
vertices and edges) is not necessarily induced from the natural cell structure
of $X$, and in general it might be ``coarser'' than the structure induced from $X$. 
Note that any block has no vertices of degree 1 and, if it is not a circle,
it also has no loop edges.

\medskip\noindent
{\bf 7.12 Theorem (Main Theorem).}
{\it Let $\Gamma$ be a finite simplicial graph which is flag, 
and denote by $|\Gamma|$
the underlying topological graph. Suppose that 
for any block $Y$ of $|\Gamma|$ which is not a circle,
each essential edge of $Y$
(which automatically is not a loop) is the union of at least 3 edges of $\Gamma$.
Let $(W_\Gamma,S_\Gamma)$ be the right angled Coxeter system 
with nerve $\Gamma$.
Then the visual boundary $\partial_\infty(W_\Gamma,S_\Gamma)$ 
is homeomorphic
to the reflection tree of graphs $|\Gamma|$, 
i.e. $\partial_\infty(W_\Gamma,S_\Gamma)\cong{\cal X}^r(|\Gamma|)$.}

\medskip\noindent
{\bf Proof:}   
Suppose first that
$\Gamma$ is connected, distinct from a singleton, 
satisfies assumptions ($*$1)
and ($*$2) from Section 6, and contains no vertex of degree 1.
Let $(D_n)$ be an exhausting sequence for $(\Sigma_\Gamma,x_0)$
as asserted by Approximation Lemma (Lemma 7.11).
By Lemma 7.6, the visual boundary $\partial_\infty(W_\Gamma,S_\Gamma)$
is homeomorphic to the inverse limit 
$\lim_{\longleftarrow}(\{ \partial D_n \},\{ g_n \})$.
It follows from Approximation Lemma that the inverse sequence
$(\{ \partial D_n \},\{ g_n \})$ has a Brown's approximation
$(\{ \partial D_n \},\{ g_n' \})$ which is isomorphic to some
null and dense inverse sequence $\cal R$ of $X$-graphs and $X$-blow-ups.
By Proposition 5.5 we get
$$
\partial_{\infty}(W_\Gamma,S_\Gamma)\cong
\lim_{\longleftarrow}(\{ \partial D_n \},\{ g_n \})\cong
\lim_{\longleftarrow}(\{ \partial D_n \},\{ g_n' \})\cong
\lim_{\longleftarrow}{\cal R}\cong
{\cal X}^r(|\Gamma|) \leqno{(7.12.1)}
$$ 
for graphs $\Gamma$ considered in this case.

Now, let $\Gamma$ be any graph satisfying the assumptions of the theorem, 
and let $Y_1,\dots,Y_k$ be the family
of factors of some terminal split decompositions of all connected components
of $|\Gamma|$. Suppose that the first $m$ of these factors $Y_i$
form the set of all blocks of $|\Gamma|$,
while the others are all trivial. For $i=1,\dots,k$ denote by
$(W_i,S_i)$ the Coxeter system of the special subgroup of $W_\Gamma$
spanned on the subset $S_i$ consisting of those generators of $S_\Gamma$
which correspond to the vertices of $\Gamma$ contained in $Y_i$.

If $m=0$ (i.e. $|\Gamma|$ has no blocks) then $\Gamma$ is
a tree or a disjoint union of trees. If $\Gamma$ is 
a singleton or a doubleton,  the theorem holds true by a direct inspection.
Otherwise, it is easy to verify (e.g. by referring to [D], Proposition 8.8.2 and
Theorem 8.7.4) that $W_\Gamma$
is then virtually free non-abelian, and thus its visual boundary is homeomorphic
to the Cantor set. On the other hand, by Lemma 2.4(2) the space
${\cal X}^r(|\Gamma|)$ is then also homeomorphic to the Cantor set.
This verifies the assertion of the theorem in the case $m=0$.

If $m=k=1$ and $|\Gamma|=Y_1$ is a circle $S^1$, it is well known that
$W_\Gamma$ is a cocompact reflection group in the euclidean or 
the hyperbolic plane, see Example 14.2.2 in [D]. 
As a consequence, we have 
$\partial_\infty(W_\Gamma,S_\Gamma)\cong S^1$.
Since by Lemma 2.4(3) we have also 
${\cal X}^r(|\Gamma|)={\cal X}^r(S^1)\cong S^1$,
the assertion follows in this case. Furthermore,
if $m=k=1$ and $|\Gamma|=Y_1$ is not a circle, we are in the setting of the
first paragraph of the proof (because any block distinct from a circle 
has no loop edges and no vertices of degree 1), and by (7.12.1) we get the assertion as well.

Finally, we are left with the case when $m\ge1$ and $k\ge 2$.
In this case $W_\Gamma$ has infinitely many ends (see [D], Theorem 8.7.4).
Note that the subgroups $W_i$ 
form a family of special subgroups
of $W_\Gamma$ corresponding to 
a split decomposition of the graph $\Gamma$ (i.e. a decomposition as
described in Definition 3.3), 
and this coincides in our case with a more general decomposition of simplicial complexes
as desribed in Subsection 7.1
of the paper [Amalgam]. It is shown in [Amalgam] that if a Coxeter group is not
2-ended then its visual boundary is homeomorphic to the dense amalgam
of the visual boundaries of its special subgroups corresponding 
to any decomposition
of its nerve (see Proposition 7.3.2 in that paper). In particular,
in our setting we get that
$$
\partial_\infty(W_\Gamma,S_\Gamma)\cong\tilde\sqcup
(\partial_\infty(W_1,S_1),\dots,\partial_\infty(W_k,S_k)).
$$
Observe that for $i>m$ each of the groups $W_i$ is either finite or virtually
free (including the 2-ended case), and consequently its visual boundary
is either empty, or homeomorphic to a doubleton or to the Cantor set. 
It follows then from Proposition 3.7(4), and from the conventions concerning
the dense amalgam of a tuple containing the empty space, that
$$
\partial_\infty(W_\Gamma,S_\Gamma)\cong\tilde\sqcup
(\partial_\infty(W_1,S_1),\dots,\partial_\infty(W_m,S_m)).
$$
Since by (7.12.1) we have that
$\partial_\infty(W_i,S_i)\cong{\cal X}^r(Y_i)$ for all $i\le m$,
denoting by $Y_1',\dots,Y_p'$ the homeomorphism types of the blocks
$Y_1,\dots,Y_m$ and
applying Propositions 3.7(3) and 3.8, we get
$$
\partial_\infty(W_\Gamma,S_\Gamma)\cong\tilde\sqcup
({\cal X}^r(Y_1),\dots,{\cal X}^r(Y_m))\cong
\tilde\sqcup
({\cal X}^r(Y_1'),\dots,{\cal X}^r(Y_p'))
\cong{\cal X}^r(|\Gamma|),
$$
which completes the proof.


\magnification1200

\bigskip\noindent
{\bf 8. The proof of Approximation Lemma.}

\smallskip
This section is devoted entirely to the proof of Approximation Lemma
(Lemma 7.11).
The proof is long and requires many observations and partial
constructions.

\medskip\noindent
{\it Geodesic rays and the branch locus of $\Sigma_\Gamma$.}

We begin with some observations concerning the 
interactions between branch components and geodesic
rays in $\Sigma_\Gamma$ started at the base vertex $x_0$. 
Recall that all branch components of $\Sigma_\Gamma$
are trees (see Fact 6.6(1)), and that they are strictly
convex in $\Sigma_\Gamma$ (Fact 6.6(4)). 
Denote by $C_0$ the branch component of $\Sigma_\Gamma$
containing $x_0$. The next claim follows directly from convexity
of $C_0$ in $\Sigma_\Gamma$.

\medskip\noindent
{\bf 8.1 Claim.} {\it Each ray in the tree $C_0$ issuing from $x_0$
is a geodesic ray of $\Sigma_\Gamma$.}

\medskip
Recall that for any closed convex subset $A$ of a CAT(0) space $Z$ 
and any point $x\in Z\setminus A$ there is a unique point
$x'\in A$ which is closest to $x$ (see [BH], Proposition II.2.4). 
This point is called the
{\it projection} of $x$ to $A$. Since branch components are closed 
and convex in $\Sigma_\Gamma$, we may speak of projections of the point
$x_0$ to them, and we have the following property of such projections.

\medskip\noindent
{\bf 8.2 Claim.} {\it For each branch component $C\ne C_0$ the
projection of $x_0$ to $C$, denoted $x_C$, is a vertex of $\Sigma_\Gamma$.}

\medskip\noindent
{\bf Proof:} If $x_C$ were not a vertex, it would be an interior
point of some edge $e$ of $C$. By the property of projection,
the geodesic segment $[x_0,x_C]$ would be then orthogonal to $e$ 
at $x_C$. It is not hard to see that in any square complex a geodesic
hitting some edge orthogonally at an interior point
does not pass throug any vertex of the complex.
This contradicts the fact that $x_0$ is a vertex, thus proving 
the claim.

\medskip
For each branch component $C\ne C_0$ denote by 
$\Omega_C$ this sector of $\Sigma_\Gamma$ from which the geodesic
$[x_0,x_C]$ approaches $x_C$. Denote also by $e_{C,1}$ and
$e_{C,2}$ the two edges in $C$ issuing from $x_C$ which are
adjacent to (i.e. contained in the boundary of) $\Omega_C$.
Denote by $x_{C,1},x_{C,2}$, respectively, the endpoints of the edges
$e_{C,1},e_{C,2}$ other than $x_C$.

\medskip\noindent
{\bf 8.3 Claim.} {\it If $C\ne C_0$ is a branched component of $\Sigma_\Gamma$,
then}
\item{(1)}
{\it any ray in $C$ started at $x_C$ and not passing along
$e_{C,1}$ or $e_{C,2}$ is a part of a geodesic ray of $\Sigma_\Gamma$
started at $x_0$;}
\item{(2)}
{\it for $i=1,2$, any ray in $C$ started at $x_{C,i}$ and not passing along
$e_{C,i}$ is a part of a geodesic ray of $\Sigma_\Gamma$ started at $x_0$.}

\medskip\noindent
{\bf Proof:} We sketch the proof of part (2), skipping that of
part (1), which is analogous. 
It is not hard to deduce from Claim 6.4 that
the interior 
$\hbox{int}(\Omega)$ separates $x_0$ from $C$.
It follows that, for $i=1,2,$ the
geodesic $[x_0,x_{C,i}]$ approaches $x_{C,i}$ either along the
edge $E_{C,i}$ or directly from the interior of the sector $\Omega_C$. 
Moreover, in the latter case the angle at $x_{C,i}$ between
$[x_0,x_{C,i}]$ and $e_{C,i}$ is strictly less than $\pi/2$.
As a consequence, by strict convexity of $C$ in $\Sigma_\Gamma$,
the angle at $x_{C,i}$ between 
the geodesic $[x_0,x_{C,i}]$ and any edge $e$ of $C$
distinct from $e_{C,i}$ is strictly larger than $\pi$.
It follows that $[x_0,x_{C,i}]$ extends as a geodesic along $e$.
By this, and by convexity of $C$ in $\Sigma_\Gamma$, 
any extension of $[x_0,x_{C,i}]$ by a ray as in the statement
yields a geodesic ray in $\Sigma_\Gamma$. This completes the proof.

\medskip\noindent
{\it Bifurcations and shadows of geodesics.}

Next observations deal with bifurcations of geodesics in $\Sigma_\Gamma$,
i.e. the phenomenon when a geodesic has a non-unique extension.
Recall that in a piecewise euclidean CAT(0) complex $K$ the possible
local extensions of a geodesic $\alpha=[x,y]$ behind its end $y$
are described in terms of a {\it shadow} of $\alpha$ at $y$,
denoted $\hbox{Sh}_y\alpha$. By definition, this is the subset of all
points in the spherical link $\hbox{Lk}_yK$ lying at distance $\ge\pi$
from the point $a\in\hbox{Lk}_yK$ representing $\alpha$.
Each $b\in\hbox{Sh}_y\alpha$ represents uniquely a germ
of a geodesic $\beta=[y,z]$ which yields an extension of $\alpha$.
(We refer the reader to [BH], Section I.7.14, for the description
of spherical links in piecewise euclidean complexes, and to Lemma (2d.1)
in [DJ] for the proof of the above mentioned characterization of geodesic
extensions in terms of shadows.)

An endpoint $y$ of a geodesic segment $\alpha$ is a {\it bifurcation point}
if the cardinality of the shadow $\hbox{Sh}_y\alpha$ is greater than 1.
It is not hard to observe that, under our assumptions on $\Gamma$,
the shadow of any geodesic segment in $\Sigma_\Gamma$ at its any endpoint
is nonempty. For geodesics with an endpoint at a vertex $v$ of $\Sigma_\Gamma$
(at which the spherical link $\hbox{Lk}_v\Sigma_\Gamma$ is isometric
with $\Gamma_{\pi/2}$), this follows e.g. by observimng that any point
in $\Gamma_{\pi/2}$ lies on a cycle of perimeter $\ge2\pi$ . For geodesics with
other endpoints, this is even more obvious. Thus, the endpoint of any geodesic
segment in $\Sigma_\Gamma$ is either a point of unique local extension,
or a bifurcation point. As a consequence, any geodesic segment in 
$\Sigma_\Gamma$ can be extended to a geodesic ray 
(i.e. $\Sigma_\Gamma$ is geodesically complete).

\medskip\noindent
{\bf 8.4 Claim.}
{\it Geodesics in $\Sigma_\Gamma$ have no
bifurcations at interior points of the sectors.}

\medskip\noindent
{\bf Proof:} At each interior point $y$ of a sector in $\Sigma_\Gamma$
the spherical link $\hbox{Lk}_y\Sigma_\Gamma$ is the standard circle 
(of perimeter $2\pi$). Consequently, any shadow is a singleton,
and hence geodesics have unique extensions.

\medskip\noindent
{\bf 8.5 Claim.}
{\it Geodesics in $\Sigma_\Gamma$ started at $x_0$ have no
bifurcations at interior points of all edges $e$ of the following two kinds:}
\item{(1)} {\it edges contained in the branch
component $C_0$ containing $X_0$;}
\item{(2)} {\it edges distinct from $e_{C,1}$ 
and $e_{C,2}$ contained in any other branch component $C$.}

\medskip\noindent
{\bf Proof:} Let $y$ be an interior point in any edge $e$ as in the 
statement. By Claim 8.3, the geodesic $\alpha=[x_0,y]$ approaches $y$
along $e$ (more precisely, along this part of $e$ which is
closer to $x_0$ or $x_C$ in the corresponding branch component). 
In the spherical link $\hbox{Lk}_y\Sigma_\Gamma$
the only point at distance $\ge\pi$ from the point $a$ representing
$\alpha$ is the point $b$ representing the geodesic issuing from $y$
and going along the other part of $e$ (more distant from $x_0$ 
or $x_C$). This completes the proof.

\medskip\noindent
{\it Shadows at vertices of $\Sigma_\Gamma$.}

Recall that link of $\Sigma_\Gamma$ at any vertex $x$ is,
by definition of $\Sigma_\Gamma$, canonically
isomorphic to the graph $\Gamma$. When viewed as spherical
link $\hbox{Lk}_x\Sigma_\Gamma$, it is canonically {\it isometric}
to $\Gamma$ equipped with the length metric for which each
edge has length $\pi/2$. We denote the resulting metric space
by $\Gamma_{\pi/2}$, to emphasise the above mentioned metric 
with which it is equipped. 
We will view the spherical vertex links of $\Sigma_\Gamma$
(under their identifications with $\Gamma_{\pi/2}$) naturally as
$|\Gamma|$-graphs, for the tautological labellings of their essential vertices
and edges. We will also view as $|\Gamma|$-graphs, with the induced
labellings, various subspaces of the vertex links, notably the shadows
of geodesics terminating at these vertices. 

We now describe the shadows of geodesics $[x_0,x]$, as subsets of 
$\Gamma_{\pi/2}=\hbox{Lk}_x\Sigma_\Gamma$,
for all vertices $x\ne x_0$ in $\Sigma_\Gamma$.
We use the following notation: if $a\in\Gamma_{\pi/2}$
then $$\hbox{Sh}(a):=\{ b\in\Gamma_{\pi/2}:d(a,b)\ge\pi \},$$
where $d$ is the metric in $\Gamma_{\pi/2}$.
The following claim is obvious.

\medskip\noindent
{\bf 8.6 Claim.}
{\it  Let $x\ne x_0$ be any vertex in $\Sigma_\Gamma$ 
and let $C$ be the component of the branch locus containing $x$.
Suppose that the geodesic
$[x_0,x]$ terminates with an edge $e$ contained in $C$. 
Let $v\in V$ be the label of $e$
(and at the same time the point of $\Gamma_{\pi/2}$).
Then $\hbox{Sh}_x[x_0,x]=\hbox{Sh}(v)$ and thus, as 
a $|\Gamma|$-graph, this shadow is isomorphic to
$|\Gamma|\setminus U_v$, where $U_v$ is some open normal
neighbourhood of $v$ in $|\Gamma|$.}

\medskip\noindent
{\bf 8.7 Remark.}
Note that, in view of Fact 8.3, Claim 8.6 applies to
\item{(1)} all vertices $x\ne x_0$ contained in $C_0$, and
\item{(2)} all vertices in any other branch component $C$
which are distinct from $x_C,x_{C,1}$ and $x_{C,2}$.

\medskip
Slightly less obvious is the following observation concerning
vertices $x=x_{C,i}$.

\medskip\noindent
{\bf 8.8 Claim.}
{\it  Let $x$ be any of the vertices of form $x_{C,i}$, for any
branch component $C\ne C_0$, and let $v\in V$
be the label of the edge $e_{C,i}$
(and at the same time the point of $\Gamma_{\pi/2}$). Then the shadow 
$\hbox{Sh}_x[x_0,x]\subset \Gamma_{\pi/2}$, viewed as a
$|\Gamma|$-graph, is isomorphic to
$|\Gamma|\setminus U_v$, where $U_v$ is some open normal
neighbourhood of $v$ in $|\Gamma|$.}

\medskip\noindent
{\bf Proof:} Let $a_0$ be the point in the spherical link
$\hbox{Lk}_x\Sigma_\Gamma=\Gamma_{\pi/2}$ representing
the geodesic $[x_0,x]$, 
and let $\varepsilon\in E$ be such that the sector $\Omega_C$
(as defined right before Claim 8.3)
is an $\varepsilon$-sector.
As it was shown in the first part of proof of Claim 8.3, 
$a_0$ either coincides 
with $v$, or it is an interior point of the essential edge 
$\varepsilon$ lying at distance less than $\pi/2$ from $v$.
In any case the shadow 
$\hbox{Sh}_x[x_0,x]=\hbox{Sh}(a_0)$ has the form
$|\Gamma|\setminus U_v$, as required, which completes
the proof.

\medskip
Finally, we deal with the shadows of the geodesics $[x_0,x_C]$
at their endpoints $x_C$. Recall that $\Omega _C$ is the sector of
$\Sigma_\Gamma$ as described right before Claim 8.3.

\medskip\noindent
{\bf 8.9 Claim.}
{\it  Let $x$ be any of the vertices of form $x_C$, for any
branch component $C\ne C_0$, and let $\varepsilon\in E$ 
be such that the sector $\Omega_C$
is an $\varepsilon$-sector. Then the shadow 
$\hbox{Sh}_x[x_0,x]\subset \Gamma_{\pi/2}$, viewed as a
$|\Gamma|$-graph, is isomorphic to
$|\Gamma|\setminus U$, where $U\subset|\Gamma|$ 
has one of the following forms:}
\item{(1)} {\it $U$ is an open segment contained in the interior of
$\varepsilon$ (in fact, this case consists of three subcases:
0, 1 or 2 endpoints of $U$ are the endpoints of $\varepsilon$);}
\item{(2)} {\it $U$ is an open normal neighbourhood of
one of the endpoints of $\varepsilon$ in $|\Gamma|$;}
\item{(3)} {\it $U$ is the union of the interior of $\varepsilon$ 
and an open normal neighbourhood of
one of the endpoints of $\varepsilon$ in $|\Gamma|$;}
\item{(4)} {\it $U$ is an open normal neighbourhood of $\varepsilon$
in $|\Gamma|$.}

\medskip\noindent
{\bf Proof:} Let $a_0$ be the point in the spherical link
$\hbox{Lk}_x\Sigma_\Gamma=\Gamma_{\pi/2}$ representing
the geodesic $[x_0,x]$, and let $v_i\in V$, for $i=1,2,$ 
be the labels of the edges
$e_{C,i}$ (and at the same time the endpoints of 
the essential edge $\varepsilon$ in $|\Gamma|$).
Since $x$ is the projection of $x_0$ to $C$, the distances
in $\Gamma_{\pi/2}$ between $a_0$ and $v_i$ are both $\ge\pi/2$.
The assertion follows then directly by considering the cases
whether these distances are smaller, equal or larger than $\pi$. 
We omit the details.

\medskip\noindent
{\it A filtration of the branch locus $C_\Gamma$.}

In order to describe an appropriate sequence $D_n$ of strict star domains
in $(\Sigma_\Gamma,x_0)$, as asserted in Approximation Lemma,
we use certain auxilliary object that we call a {\it filtration} of
$C_\Gamma$, which we now introduce. 
To do this, we need more terminology, notation
and some observations.
For any sector $\Omega$ of $\Sigma_\Gamma$ let 
$\partial_{x_0}\Omega$ be this component of $\partial\Omega$
which either contains $x_0$ or separates $\hbox{int}(\Omega)$
from $x_0$. Existence and uniqueness of such component
of $\partial\Omega$ follows from Facts 6.4(3) and 6.9.
Next, denote by $H_\Omega$ the closure 
of this connected component of $\Sigma_\Gamma\setminus\partial_{x_0}\Omega$
which contains $\hbox{int}(\Omega)$. Finally, let
$g_\Omega:H_\Omega\to\partial_{x_0}\Omega$ be the geodesic
projection (i.e. the map which, to any point $y\in H_\Omega$
associates the point closest to $y$ in the intersection of the 
geodesic $[x_0,y]$ with $\partial_{x_0}\Omega$).
In general, $g_\Omega$ is not continuous, but it satisfies
the following.

\medskip\noindent
{\bf 8.10 Claim.}
{\it  Let $\Omega$ be a sector of $\Sigma_\Gamma$ and let $C$
be any branch component adjacent to $\Omega$ and
not containing $\partial_{x_0}\Omega$. Then 
$g_\Omega(C)=g_\Omega(e_{C,1}\cup e_{C,2})$.
Moreover, $g_\Omega(C)$ is a bounded subset of
$\partial_{x_0}\Omega$.}

\medskip\noindent
{\bf Proof:} The first assertion follows easily from Claim 8.3.
In view of the first assertion, the second one follows by observing 
that the union $\bigcup\{ [x_0,y]:y\in e_{C,1}\cup e_{C,2} \}$
is a bounded subset of $\Sigma_\Gamma$, and the image
$g_\Omega(C)$ is its subset.

\medskip\noindent
{\bf 8.11 Definition.}
A {\it filtration} of the branch locus $C_\Gamma$
is any sequence $F_n:n\ge1$ of finite subgraphs of $C_\Gamma$
satisfying the following conditions:
\item{(1)} $F_1=\{ x_0 \}$, $F_{n+1}$ has exactly one more 
vertex than $F_n$ for all $n\ge1$, and $\bigcup F_n=C_\Gamma$;
\item{(2)} for each $n$ the intersection of $F_n$ with any branch
component $C$ is either empty or a subtree containing
the vertex $x_C$ (we use the convention that $x_{C_0}=x_0$);
\item{(3)} for any branch component $C\ne C_0$, if $i$ is
the smallest
 index for which $x_C\in F_i$, then 
$g_{\Omega_C}(C)\subset F_{i-1}$ (here again $\Omega_C$ denotes
the sector of $\Sigma_\Gamma$ described right before Claim 8.3).

\medskip
We skip an easy argument for
showing that a filtration of $C_\Gamma$
always exists (in fact, there are lots of choices and plenty
of filtrations). We only note that, in order to construct a filtration, 
it is convenient to order the components of the branch locus $C_\Gamma$
into a sequence $C_0,C_1,\dots$ such that:
\item{$\bullet$} $C_0$ is the component containing $x_0$;
\item{$\bullet$} for each $n\ge0$, the subtree of the adjacency tree
${\cal A}_\Gamma$ (described in Definition 6.7 and in Fact 6.9)
spanned on the vertices corresponding to the branch components
$C_0,\dots,C_n$ contains no other vertices corresponding to branch components.

\noindent
It is also helpful to note that if $\partial_{x_0}\Omega_{C_n}\subset C_i$
then $i<n$, and to refer to Claim 8.10. We omit further details.

\medskip\noindent
{\it Description of the domains $D_n$ and the inverse sequence 
induced by a filtration.}

We now describe an exhausting sequence of regular strict star
domains $(D_n)$ in $(\Sigma_\Gamma,x_0)$, as required in 
Approximation Lemma. It will be determined essentially uniquely by
a choice of a filtration of the branch locus $C_\Gamma$.

Given a filtration $(F_n)$ of $C_\Gamma$, for each $n$ consider
a regular strict star domain $D_n$ satisfying the following
properties:

\itemitem{(d1)} for each branch component $C$, $C$ intersects $\hbox{int}(D_n)$
iff  $C\cap F_n\ne\emptyset$, and if this is the case then $C\cap\hbox{int}(D_n)$
is a connected subset of $C$ containing $C\cap F_n$ and containing no
other vertices of $C$ than those of $C\cap F_n$;

\itemitem{(d2)} if for some $C\ne C_0$ we have 
$F_n\cap C=\{ x_C \}$ then $x_{C,i}\in\partial D_n$ for $i=1,2$.

\noindent
It is fairly obvious that domains $D_n$ as above exist.
Moreover, we have the following uniqueness result.

\medskip\noindent
{\bf 8.12 Claim.}
{\it If $D_n\subset D_n'$ are two regular strict star domains 
satisfying conditions (d1)-(d2) then the geodesic projection
$g:\partial D_n'\to\partial D_n$ is a homeomorphism
which is also an isomorphism of $|\Gamma|$-graphs.
As a consequence, any two regular strict star domains $D_n$ 
satisfying conditions (d1)-(d2) are isomorphic as $|\Gamma|$-graphs.}

\medskip\noindent
{\bf Proof:} The first assertion follows from the fact that
geodesics started at $x_0$ have no bifurcations in the set
$\hbox{int}(D_n')\setminus\hbox{int}(D_n)$
(this is a direct consequence of Claims 8.4 and 8.5),
and from the fact (following from Claims 8.1, 8.3 and 8.5)
that for each $y'\in\partial D_n'$ its image $y=g(y')$
satisfies the following:
\item{$\bullet$} if $y'\in\hbox{int}(\Omega)$ for some sector 
$\Omega$ then $y\in\hbox{int}(\Omega)$; 
\item{$\bullet$} if $y'$ lies in the interior of some branch edge
of $\Sigma_\Gamma$ then $y$ belongs to the interior of the same edge;
\item{$\bullet$} if $y'$ is a vertex of $\Sigma_\Gamma$ then $y=y'$.

\noindent
To prove the second assertion, note that for any two regular strict
star domains $D_n,D_n'$ satisfying (d1)-(d2), their intersection
is a domain of the same kind. Thus the second assertion follows 
from the first one.

\medskip
We also obviously have the following.

\medskip\noindent
{\bf 8.13 Claim.} 
{\it The $|\Gamma|$-graph $\partial D_1$ is isomorphic to
$|\Gamma|$ (with the tautological labelling of its essential vertices
and edges).}

\medskip
We next observe that
a sequence of domains $(D_n)$ as above can easily be chosen
so that it forms an increasing exhausting sequence. It follows from
Claim 8.12 that the induced inverse sequence 
$(\{ \partial D_n \},\{ g_n \})$, where $g_n:\partial D_{n+1}\to D_{n}$
are the geodesic projections, is then unique up to isomorphism
of inverse sequences of $|\Gamma|$-graphs.
We will call this sequence {\it the inverse sequence induced
by the filtration $(F_n)$}.

To conclude the proof of Approximation Lemma, it remains to show
that the maps $g_n$ in the above inverse sequence 
can be approximated by appropriate $|\Gamma|$-blow-ups.
To do this, we will use, among others, the previous analysis of shadows.

\medskip\noindent
{\it Approximating geodesic projections by $|\Gamma|$-blow-ups.}

Given a filtration $(F_n)$ of the branch locus $C_\Gamma$,
let $(\{ \partial D_n \},\{ g_n \})$ be the inverse sequence of
$|\Gamma|$-graphs induced by this filtration. For each $n\ge2$,
denote by $x_n$ the unique vertex in $F_n\setminus F_{n-1}$.
We construct the required approximations of the maps $g_{n-1}$
separately in the three cases corresponding to the three types of the 
vertices $x_n$.

\medskip\noindent
{\bf Case 1:} $x_n\in C_0$ or 
$x_n\in C\setminus\{ x_C,x_{C,1},x_{C,2} \}$ for some 
$C\ne C_0$.

Note that, since $n\ge2$, we have $x_n\ne x_0$. Thus, in view
of Remark 8.7, $x_n$ satisfies the assumption of Claim 8.6. 
We pass to the notation as in Claim 8.6, with $x=x_n$ and with 
$e$ coinciding with the unique edge in $F_n$ adjacent to $x_n$.
Denote by $v_e$ the vertex in the $|\Gamma|$-graph 
$\partial D_{n-1}$ corresponding to the intersection point
$\partial D_{n-1}\cap e$.
It follows from Claim 8.6, and from the arguments as in the proof
of Claim 8.12, that the $|\Gamma|$-graph $\partial D_n$ 
and the map $g_{n-1}$ have the following form.
The graph $\partial D_n$ can be naturally viewed as obtained
from $\partial D_{n-1}$ by:
\item{$\bullet$} first, deleting $v_e$ and replacing it with as many
points as edges adjacent to $v_e$;
\item{$\bullet$} then, gluing to the so obtained graph
the shadow $\hbox{Sh}(v)=\Gamma\setminus U_v$,
through the map which associates to the points which 
have replaced $v_e$
the points of $\partial U_v$, in the way respecting labels
of the edges containing these points.

\noindent
Under this perspective on $\partial D_n$, the map $g_{n-1}$
has the following description:
\item{$\bullet$} the points of $\partial D_n$ corresponding to
$\partial D_{n-1}\setminus\{ v_e \}$ are mapped identically
to the same points in $\partial D_{n-1}$;
\item{$\bullet$} the points of $\partial D_n$ corresponding to
the shadow $\Gamma\setminus U_v$ are all mapped to $v_e$.

\noindent
The map $g_{n-1}$
can be approximated, arbitrarily close with respect to
the uniform distance, by maps $g_{n-1}'$ of the following form:
\item{$\bullet$} for any edge $\eta$ of $\partial D_{n-1}$
issuing from $v_e$ choose a point $p_\eta$ on $\eta$ sufficiently close
to $v_e$; denote by $q_\eta$ the new point at the end of $\eta$
that has replaced $v_e$ in the description of $\partial D_n$ above;
let $q_\eta'$ be the point in $\Gamma\setminus U_v$
to which $q_\eta$ is glued, and let $v_\eta$ be the endpoint
other than $v$ in the edge of $\Gamma$ containing $q_\eta'$;
\item{$\bullet$} put $g_{n-1}'$ to coincide with $g_{n-1}$ on the part
of $\partial D_n$ corresponding to 
$\partial D_{n-1}\setminus\bigcup_\eta[p_\eta,v_e]$
and on the part corresponding to the complement of the open 
star of the vertex $v$ in $|\Gamma|$
(i.e. on the subgraph of $|\Gamma|$, for the natural stratification,
spanned on all essential vertices except $v$); 
moreover, for each $\eta$
as before put $g_{n-1}'$ to map the segment 
$[p_\eta,q_\eta]\cup[q_\eta',v_\eta]$ in $\partial D_n$ 
homeomorphically to the segment $[p_\eta,v_e]$ in
$\partial D_{n-1}$.

\noindent
Any map $g_{n-1}'$ described above clearly has the form of
a map associated to the $|\Gamma|$-blow-up of 
$\partial D_{n-1}$ at $v_e$, and $g_{n-1}$ can be approximated
by such a map as close as necessary.

\medskip\noindent
{\bf Case 2:} $x_n=x_{C,i}$ for some $C\ne C_0$ and some
$i\in\{ 1,2 \}$.

The argument in this case is the same as in Case 1, except that
the role of the vertex $v_e$ is now played by $x_n$ itself,
the roles of $e$ and $v$ are played by the edge $e_{C,i}$
and by its label $v\in V$, respectively,
and we use Claim 8.8 in place of Claim 8.6.
We omit further details.

\medskip\noindent
{\bf Case 3:} $x_n=x_C$ for some $C\ne C_0$. 

As we show below, in this case, which is slightly different from
the previous ones, the map $g_{n-1}$ can be approximated by
maps $g_{n-1}'$ corresponding to blow-ups at a segment rather 
than at a vertex. Accordingly with the variety of forms for
the set $U$ in Claim 8.9, this case splits into many subcases.
We deal in detail with only one of those subcases, leaving
the remaining ones (in which the argument is very similar)
to the reader.

We start with the setting for the case under consideration.
We use the notation from Claim 8.9 and from its proof.
Let $a_0$ be the point in the spherical link 
$\hbox{Lk}_{x_n}\Sigma_\Gamma=\Gamma_{\pi/2}$ representing
the geodesic $[x_0,x_n]$, and let $v_i$ be the the points
in this link representing the edges $e_{C,i}$. We assume that
in the considered subcase
the distance in the link from $a_0$ to $v_1$ is smaller than 
$\pi$, while the corresponding distance to $v_2$ is greater than
$\pi$.

Fix the following notation. Put $z_0:=[x_0,x_n]\cap\partial D_{n-1}$
and $z_1:=[x_0,x_{C,1}]\cap\partial D_{n-1}$ and note that,
since the distance from $a_0$ to $v_1$ in the link at $x_n$ 
is less than $\pi$, these two points are distinct. Assume that
the sector $\Omega_C$ is labelled with $\varepsilon\in E$.
Then the graph $\partial D_{n-1}$ contains an edge 
(labelled with $\varepsilon$) corresponding
to its part contained in $\Omega_C$, and we denote this edge
by $\varepsilon_{n-1}$. Note that both $z_0$ and $z_1$ belong
to the interior of this edge.
Further, denote by $a_2$ the point in the link
$\hbox{Lk}_{x_n}\Sigma_\Gamma$, lying at distance $\pi$ from $a_0$
on the segment $[a_0,v_2]$. This is one of the "extremal" points
in the shadow $\hbox{Sh}(a_0)$.
Consider the geodesic issuing from $x_n$ which corresponds to $a_2$
and note that it does not bifurcate before reaching $\partial D_n$;
denote by $y_2$ the point of intersection of this geodesic
with $\partial D_n$. Similarly, for each essential edge $\eta$ other than
$\varepsilon$ issuing from $v_1$ in 
$\hbox{Lk}_{x_n}\Sigma_\Gamma=\Gamma_{\pi/2}$, denote by $a_\eta$
the interior point of $\eta$ having distance $\pi$ from $a_0$.
Points of the form $a_\eta$ are the remaining "extremal" points
in the shadow $\hbox{Sh}(a_0)$. For each $\eta$, 
consider the geodesic issuing from $x_n$ which corresponds to $a_\eta$
and note that it does not bifurcate before reaching $\partial D_n$;
denote by $y_\eta$ the point of intersection of this geodesic
with $\partial D_n$.

We make the following observations, which fully describe the form
of the map $g_{n-1}$:
\item{(1)} points of the preimage $g_{n-1}^{-1}(z_0)$ are in the
natural bijective correspondence (which is also an isomorphism
of $|\Gamma|$-graphs) with the points of the shadow
$\hbox{Lk}_{x_n}\Sigma_\Gamma=\hbox{Sh}(a_0)$; this correspondence
is given by associating to a point $y\in g_{n-1}^{-1}(z_0)$ the point
in the link induced by the geodesic $[x_n,y]$;
\item{(2)} for each $\eta$ as above, there is a segment 
$[y_\eta,x_{C,1}]$ in the graph $\partial D_n$, and each such segmnet
is mapped by $g_n$ homeomorphically on the segment $[z_0,z_1]$
in $\partial D_{n-1}$;
\item{(3)} the complement in $\partial D_n$ of the set 
$g_{n-1}^{-1}(z_0)\cup\bigcup_\eta[y_\eta,x_{C,1}]$
is mapped by $g_n$ isomorphically (by isomorphism of 
$|\Gamma|$-graphs) on the complement 
$\partial D_{n-1}\setminus[x_0,x_1]$.

Referring to the above description, we can approximate the map $g_{n-1}$, arbitrarily close, by maps $g_{n-1}'$ of the following form:
\item{$\bullet$} consider the segment $[y_2,x_{C,2}]$ 
in the graph $\partial D_n$ and let $y_3\in\partial D_n$ 
be a point lying outside this segment as close to $y_2$ as necessary;
similarly, let $y_4\in\partial D_n$ be a point lying outside the segment
$[y_3, x_{C,2}]$ as close to $y_3$ as necessary, and put 
$z_3:=g_{n-1}(y_3)$ and
$z_4:=g_{n-1}(y_4)$; finally, let $z_2$ be a point in the interior
of the segment 
$[z_0,z_1]$ contained in the edge $\varepsilon_{n-1}$ of 
$\partial D_{n-1}$, and assume that it is as close to $z_0$ as necessary;
\item{$\bullet$} for any edge $\xi$ of $\partial D_n$ adjacent
to $x_{C,2}$ and not containing $y_2$, choose a point $y_\xi$
in the interior of $\xi$ as close to $x_{C,2}$ as necessary;
\item{$\bullet$} 
$g_{n-1}'$ maps the segment $[y_4,y_2]\subset\partial D_n$ homeomorphically
onto the segment $[z_4,z_3]\subset\partial D_{n-1}$, similarly it maps
$[y_2,x_{C,2}]$ onto $[z_3,z_0]$, $[x_{C,2},y_\xi]$ onto $[z_0,z_2]$
for all $\xi$ as above, $[y_\eta,x_{C,1}]$ onto $[z_2,z_1]$ for all $\eta$
as above, it also maps $g_{n-1}^{-1}(z_0)\setminus[y_2,x_{C,2}]$
to the point $z_2$, and it coincides with $g_{n-1}$ on the remaining part of
$\partial D_n$.


\noindent
It is not hard to note that
any map $g_{n-1}'$ above has the form of
a map associated to the $|\Gamma|$-blow-up of 
$\partial D_{n-1}$ at the segment $[z_0,z_1]$ contained
in the interior of the edge $\varepsilon_{n-1}$. Moreover, 
$g_{n-1}$ can be clearly approximated
by such a map as close as necessary.
In one of the later arguments we will also use the following property
of the above described map $g_{n-1}'$
(in the statement of which $\pi_{D_n}$ denotes the geodesic projection
towards $x_0$ from the complement of $\hbox{int}(D_n)$ to the
boundary $\partial D_n$):
$$
\hbox{the blow-up segment of }g_{n-1}' \hbox{ coincides with the image }
\pi_{D_{n-1}}(e_{C,1}\cup e_{C,2}).
\leqno(8.1)
$$
In all other subcases of Case 3
the construction of approximations $g_{n-1}'$ can be also performed so that
property (8.1) holds. We skip further details.

\medskip
We now construct recursively 
a Brown's approximation $(\{ \partial D_n \},\{ g_n' \})$
as required. Suppose that for some $n\ge0$ we have already chosen
appropriate modified maps $g_1', \dots,g_n'$. Choose 
$g_{n+1}':\partial D_{n+2}\to\partial D_{n+1}$ to be a map as described
in the appropriate corresponding case above such that
$$
|| g_{n+1}'-g_{n+1} ||<\epsilon(g_1',\dots,g_n'),
$$
where $\epsilon(g_1',\dots,g_n')$ is the appropriate constant from Brown's Lemma
(Lemma 7.9).

It is clear from
the description of the maps $g_i'$ given in the three cases above that
these maps are $|\Gamma|$-blow-ups.

\medskip\noindent
{\it An approximating inverse sequence 
$(\{ \partial D_n \},\{ g_n' \})$ is null.}

To conclude the proof of Approximation Lemma
we need to show that a Brown's approximation
$(\{ \partial D_n \},\{ g_n' \})$ as constructed above is null and dense 
(i.e. satisfies conditions (i3) and (i4) of Definition 5.4). 

We first deal with condition (i3), i.e. nullness. 
To verify it, we need the following. 

\medskip\noindent
{\bf 8.14 Claim.} 
{\it For arbitrary $i\ge1$ consider the family of images through geodesic
projections on $\partial D_i$ of the sets $e_{C,1}\cup e_{C,2}$,
for all branch components $C$ which are disjoint with $D_i$.
Then the diameters of these images converge to zero.}

\medskip
The claim follows easily from convexity of the CAT(0) metric in $\Sigma_\Gamma$
(see Proposition II.2.2 in [BH]), since the sets $e_{C,1}\cup e_{C,2}$
have diameters uniformly bounded by 2, and their distances from
$x_0$ (and from $\partial D_i$) diverge to infinity.
\medskip

Now, observe that for each $i$ the projections to $\partial D_i$
of the blow up segments of the sequence $(\{ \partial D_n \},\{ g_n' \})$
satisfy the following conditions:
\item{$\bullet$} they are all connected;
\item{$\bullet$} they form a {\it nested} family, i.e. any two sets are either
disjoint or one of them is contained in the other 
(see properties (d2)-(d4) in Step 1 of the proof of Proposition 5.5, in Section 5).

\noindent
Note that, given a countable nested family
$A_n:n\ge1$ of connected subsets
in a finite graph, the only reason for such a family not to be null is that
there is an increasing infinite sequence $n_k$ such that:
\item{$\bullet$} for all $k\ge1$ we have
$A_{n_{k+1}}\subset A_{n_k}$, 
and 
\item{$\bullet$}
$\lim_{k\to\infty}\hbox{diam}(A_{n_k})>0$. 

\noindent
To prove that the appropriate families of projections of blow-up segements,
as in condition (i3),
are null, it is thus sufficient to exclude the phenomenon as above.

Suppose on the contrary that for some $i_0\ge1$ 
we are given a sequence $L_k$ of blow-up
segments, with $L_k\subset\partial D_{n_k}$ and $n_k>i_0$, 
such that,
denoting by $L_k'$ the projections of $L_k$ onto $\partial D_{i_0}$
through appropriate compositions of the maps $g_i'$,
we have:
\item{(1)} $L_{k+1}'\subset L_k'$ for all $k\ge1$;
\item{(2)} $\lim_{k\to\infty}\hbox{diam}(L_k')>0$.

Since all the sets $L_k'$ are compact, it follows from property (2) above
that the intersection $\cap_{k\ge1}L_k'$ consists of at least two points. 
Let $x_0^1,x_0^2$ be any two distinct points of this intersection.
Note also that, without loss of generality, we can assume that
for each $k\ge1$ we have
$n_{k+1}>n_k$, and consequently 
$g_{n_k}'\circ\dots g_{n_{k+1}-1}'(L_{k+1})\subset L_k$.
It follows that there exist strings $(x_n^1)$ and $(x_n^2)$ for the inverse
sequence ${\cal S}'=(\{ \partial D_n \},\{ g_n' \})$ such that
for $j=1,2$ we have
$x_{i_0}^j=x_0^j$ and $x_{n_k}^j\in L_k$ for all $k\ge1$.

Recall that by (8.1), for each $k$ there is a branch component
$C_k$ in $\Sigma_\Gamma$ such that 
$L_k=\pi_{D_{n_k}}(e_{C_k,1}\cup e_{C_k,2})$.
Moreover, it is not hard to observe that the components $C_k$
are then pairwise distinct. For each $k\ge1$ and for $j=1,2$ choose 
$z_{k}^j\in e_{C_k,1}\cup e_{C_k,2}$ such that 
$\pi_{D_{n_k}}(z_k^j)=x_{n_k}^j$.
Note also that, by Remark 7.10(2),
the limits $y^j=\lim_{k\to\infty}g_{k_0}\circ\dots\circ g_{n_k-1}(x_{n_k}^j)$
exist and are distinct. On the other hand, for each $k$
we have that $y^j=g_{k_0}\circ\dots\circ g_{n_k-1}\circ\pi_{D_{n_k}}(z_k^j)$,
which simply means that $y^j=\pi_{D_{k_0}}(z_k^j)$.
Since it follows then from Claim 8.14
that, as $k\to\infty$, the distance between the projections 
$\pi_{D_{k_0}}(z_k^j)$ converges to 0, 
we conclude that $y^1=y^2$,
contradicting an earlier observation that these points are distinct.

This completes the proof of nullness.

\break

\medskip\noindent
{\it An approximating inverse sequence 
$(\{ \partial D_n \},\{ g_n' \})$ is dense.}

\smallskip
To prove denseness (condition (i4)), we will need the following.

\medskip\noindent
{\bf 8.15 Claim.} 
{\it The set $M=\{ x_C:C\ne C_0 \}$ is a net in $\Sigma_\Gamma$,
i.e. there is $r>0$ such that each ball of radius $r$ in $\Sigma_\Gamma$
intersects $M$.}

\medskip
To prove the claim, we will show that for any sector $\Omega$
in $\Sigma_\Gamma$ the intersection $\Omega\cap M$ is a net
in $\Omega$. Denote by $L_0$ this boundary component of
$\Omega$ through which geodesics started at $x_0$ enter $\Omega$.
Note that each point of $\Omega$ lies at distance at most 2
from some square $Q$ of $\Omega$ disjoint with $L_0$.
Furthermore, any such square $Q$ separates $L_0$ from
some other boundary componet $L\subset\partial\Omega$
disjoint from $Q$ and lying at distance at most 2 
from $Q$. It is also not hard to see that, if $C$ denotes the branch
component in $\Sigma_\Gamma$ containing $L$ then $x_C\in L$ and
the distance from $x_C$ to $Q$ is also $\le2$.
Since $\hbox{diam}(Q)\le2$, the number $r=6$ is as required, which
completes the proof of Claim 8.15.

\medskip
As a consequence of Claim 8.15, 
we get the following.

\medskip\noindent
{\bf 8.16 Corollary.}
{\it For any $n\ge1$, let 
$\pi_{D_n}:\Sigma_\Gamma\setminus\hbox{int}(D_n)\to\partial D_n$
be the geodesic projection towards the base vertex $x_0$.
Then the set
$$
\pi_{D_n}(\{ x_C:\hbox{\rm $C$ is a branch component such that }C\cap D_n=\emptyset \})
$$
is dense in $\partial D_n$.}

\medskip
To justify the corollary, recall that under our assumtions on $\Gamma$,
the Coxeter-Davis complex $\Sigma_\Gamma$ is geodesically complete,
and hence each of the maps $\pi_{D_n}$ is surjective. The corollary thus follows
from Claim 8.15 due to convexity of the CAT(0) metric in $\Sigma_\Gamma$.

\medskip
In order to show that our inverse sequence $(\{ \partial D_n \},\{ g_n' \})$
of $|\Gamma|$-graphs and $|\Gamma|$-blow-ups is dense, we assume,
without loss of generality, that it satisfies the additional condition
described in Remark 7.10(1), and in particular the inequalities (A) mentioned
there.

Suppose that $L\subset\partial D_i$ is a blow-up segment (which means that
the map $g_i':\partial D_{i+1}\to\partial D_i$ is a $|\Gamma|$-blow-up
at $L$). It follows from the description of $g_i'$ in Case 3 above, and more
precisely from (8.1), that then for some branch component $C$ we have
$x_C\in\partial D_{i+1}$ and $g_i'(x_C)\in L$. In view of this, to prove that
the inverse sequence $(\{ \partial D_n \},\{ g_n' \})$ is dense,
it is sufficient to show the following.

\medskip\noindent
{\bf 8.17 Claim.}
{\it For any branch component $C$ of $\Sigma_\Gamma$, denote by 
$n_C$ this number for which $x_C\in\partial D_{n_C}$. Then for each $i\ge1$,
the set
$$
\{ g_i'\circ g_{i+1}'\circ\dots\circ g_{n_C-1}'(x_C):C\cap D_i=\emptyset  \}
$$
is dense in $\partial D_i$.}

\medskip
To prove Claim 8.17, fix any point $p\in\partial D_i$, and any $\varepsilon>0$.
We will show that there is $C$ such that the distance in $\partial D_i$
from $p$ to the point $g_i'\circ g_{i+1}'\dots\circ g_{n_C-1}'(x_C)$
is less than $2\varepsilon$.
Pick $n>i$ such that ${1\over n}<\varepsilon$, and let $\delta_n$ be
a constant as in Remark 7.10(1).
Let $q\in\partial D_n$ be any point such that
$g_i'\circ g_{i+1}'\circ\dots\circ g_{n-1}'(q)=p$
(such $q$ exists since all the maps $g_i'$ are surjective,
which follows from the form of $|\Gamma|$-blow-up maps in view
of the assumption that $\Gamma$ has no vertices of degree 1).
By Corollary 8.16, there is $m>n$ and a point $x_C\in\partial D_m$
such that
$$
d_{\partial D_n}(q, g_n\circ g_{n+1}\circ\dots\circ g_{m-1}(x_C))<\delta_n,
$$
where $d_{\partial D_n}$ is the metric in $\partial D_n$
(this is true because 
$g_n\circ g_{n+1}\circ\dots\circ g_{m-1}(x_C)=\pi_{D_n}(x_C)$).
By condition (A) of Remark 7.10(1), we have
$$
d_{\partial D_n}(g_n'\circ g_{n+1}'\circ\dots\circ g_{m-1}'(x_C),
g_n\circ g_{n+1}\circ\dots\circ g_{m-1}(x_C))<\delta_n.
$$
Since $\delta_n$ is a uniform continuity constant for $\epsilon_n={1\over n}$
(for the map $g_i'\circ g_{i+1}'\circ\dots\circ g_{n-1}'$), we get
$$\matrix{
d_{\partial D_i}(p,g_i'\circ g_{i+1}'\circ\dots\circ g_{n-1}'\circ g_n'\circ\dots
g_{m-1}'(x_C))= \cr
=d_{\partial D_i}(g_i'\circ g_{i+1}'\circ\dots\circ g_{n-1}'(q),g_i'\circ g_{i+1}'\circ\dots\circ g_{n-1}'\circ g_n'\circ\dots
g_{m-1}'(x_C))\le \cr
=d_{\partial D_i}(g_i'\circ\dots\circ g_{n-1}'(q),g_i'\circ\dots\circ g_{n-1}'\circ g_n\circ\dots
g_{m-1}(x_C))+ \cr
+ d_{\partial D_i}(g_i'\circ\dots\circ g_{n-1}'\circ g_n\circ\dots
g_{m-1}(x_C),g_i'\circ\dots\circ g_{n-1}'\circ g_n'\circ\dots
g_{m-1}'(x_C))< \cr
<{1\over n}+{1\over n}<2\varepsilon.
}
$$
This finishes the proof of the claim, and thus also of Approximation Lemma.


\magnification1200

\bigskip
\centerline{\bf References}
\medskip

\itemitem{[AS]} F. Ancel, L. Siebenmann,
{\it The construction of homogeneous homology manifolds},
Abstracts Amer. Math. Soc. 6 (1985), 92.

\itemitem{[Amalgam]} J. \'Swi\c atkowski,
{\it The dense amalgam of metric compacta and topological characterization of boundaries of free products of groups}, Groups, Geometry, and Dynamics 10
(2016),  407--471.

\itemitem{[BH]} M. Bridson, A. Haefliger, Metric Spaces of Non-Positive
Curvature, Grundlehren der mathematischen Wissenschaften 319, Springer,
1999.

\itemitem{[Br]} M. Brown, {\it Some applications of an appropximation
theorem for inverse limits}, Proc. Amer. Math. Soc. 11 (1960), 478--481. 

\itemitem{[Capel]}  C, E, Capel, {\it Inverse limit spaces,} Duke Math. J. 21 (1954), 
233--245.

\itemitem{[Dav]} R.J. Daverman, {Decompositions of Manifolds}, 
Academic Press, 1986.

\itemitem{[D]} M.W. Davis, 
{The geometry and topology of Coxeter groups},
London Mathematical Society Monographs Series, vol. 32,
Princeton University Press, 2008.

\itemitem{[DJ]} M. Davis, T. Januszkiewicz, 
{\it Hyperbolization of polyhedra}, J. Differential Geometry 34 (1991),
347--388.

\itemitem{[DT]} P. Dani, A. Thomas, 
{\it Bowditch's JSJ tree and the quasi-isometry classification of certain Coxeter groups}, Journal of Topology 10 (2017), 1066--1106.

\itemitem{[En]} R. Engelking, General Topology, Sigma series in pure mathematics, Vol. 6, Heldermann, Berlin, 1989.

\itemitem{[Fi]} H. Fischer, {\it Boundaries of right--angled Coxeter 
groups with manifold nerves},
Topology 42 (2003), 423--446.

\itemitem{[J]} W. Jakobsche, {\it Homogeneous cohomology manifolds which are inverse limits},
Fundamenta Mathematicae 137 (1991), 81--95.

\itemitem{[KK]} M. Kapovich, B. Kleiner, 
{\it Hyperbolic groups with low-dimensional boundary}, 
Ann. Sci. ENS 33 (2000) 647--669.

\itemitem{[PS]} P. Przytycki, J. \'Swi\c atkowski,  
{\it Flag-no-square triangulations and Gromov boundaries
in dimension 3}, Groups, Geometry \& Dynamics  3 (2009), 453--468.

\itemitem{[St]} P.R. Stallings, {\it An extension of Jakobsches construction
of $n$--ho\-mo\-geneous continua to the
nonorientable case}, in Continua (with the Houston Problem Book),
ed. H. Cook, W.T. Ingram, K. Kuperberg, A. Lelek, P. Minc,
Lect. Notes in Pure and Appl. Math. vol. 170 (1995), 347--361.

\itemitem{[Tr-metr]} J. \'Swi\c atkowski, {\it Trees of metric compacta and trees of manifolds,} preprint, 2013,  arXiv:1304.5064v2.  

\itemitem{[Tr-mfld]} J. \'Swi\c atkowski, 
{\it Trees of manifolds as boundaries of spaces and groups,}
preprint, 2013,  arXiv:1304.5067v2.


\bye